\documentclass[twoside,11pt]{article}

\usepackage[top=1in,bottom=1in,left=1in,right=1in]{geometry}

\usepackage{enumitem}
\setlist[itemize]{label=--}

\RequirePackage{amsthm,amsmath,amsfonts,amssymb}
\RequirePackage[numbers]{natbib}
\RequirePackage[colorlinks,citecolor=blue,urlcolor=blue]{hyperref}
\RequirePackage{graphicx}
\usepackage{subcaption}
\usepackage{bbm}

\usepackage{textgreek}
\theoremstyle{plain}

\newtheorem{theorem}{Theorem}[section]
\newtheorem{lemma}[theorem]{Lemma}
\newtheorem{proposition}[theorem]{Proposition}
\newtheorem{corollary}[theorem]{Corollary}
\theoremstyle{remark}
\newtheorem{definition}[theorem]{Definition}
\newtheorem{rem}[theorem]{Remark}
\newtheorem{example}[theorem]{Example}

\newtheorem{ass}{Assumption}

\def\beq{\begin{equation}} 
\def\eeq{\end{equation}}
\def\beqn{\begin{eqnarray*}}
\def\eeqn{\end{eqnarray*}}
\def\Bal{\begin{align}}
\def\Eal{\end{align}}
\def\Bitem{\begin{itemize}\setlength{\itemsep}{.2in}}
\def\bitem{\begin{itemize}\setlength{\itemsep}{.05in}}
\def\eitem{\end{itemize}}
\def\blatin{\begin{enumerate}\setlength{\itemsep}{.05in}\renewcommand{\labelenumi}{\roman{enumi}.}}
\def\elatin{\end{enumerate}}
\def\Benum{\begin{enumerate}\setlength{\itemsep}{.2in}}
\def\benum{\begin{enumerate}\setlength{\itemsep}{.05in}}
\def\eenum{\end{enumerate}}
\def\bmult{\begin{multline*}}
\def\emult{\end{multline*}}
\def\bcenter{\begin{center}}
\def\ecenter{\end{center}}
\def\bframe{\begin{frame}}
\def\eframe{\end{frame}}

\newcommand{\thmref}[1]{Theorem~\ref{thm:#1}}
\newcommand{\prpref}[1]{Proposition~\ref{prp:#1}}
\newcommand{\corref}[1]{Corollary~\ref{cor:#1}}
\newcommand{\lemref}[1]{Lemma~\ref{lem:#1}}
\newcommand{\secref}[1]{Section~\ref{sec:#1}}
\newcommand{\figref}[1]{Figure~\ref{fig:#1}}

\newcommand{\assref}[1]{Assumption~\ref{ass:#1}}

\def\cA{\mathcal{A}}

\def\cC{\mathcal{C}}

\def\cL{\mathcal{L}}

\def\cP{\mathcal{P}}

\def\cR{\mathcal{R}}
\def\cS{\mathcal{S}}

\def\bL{\mathbf{L}}

\def\bbE{\mathbb{E}}

\def\bbR{\mathbb{R}}

\def\bbX{\mathbb{X}}

\def\ind{\mathbbm{1}}

\DeclareMathOperator{\diam}{diam}

\DeclareMathOperator{\tv}{TV}
\DeclareMathOperator{\TV}{TV}

\DeclareMathOperator{\Id}{Id}

\DeclareMathOperator{\imag}{Im}
\renewcommand{\Im}{\imag}
\DeclareMathOperator{\II}{I\!I}

\DeclareMathOperator{\pr}{\pi}

\DeclareMathOperator{\vol}{vol}

\DeclareMathOperator{\rad}{Rad}

\DeclareMathOperator{\conv}{conv}

\let\lac\{
\let\rac\}
\renewcommand{\{}{\left\lac}
\renewcommand{\}}{\right\rac}
\newcommand{\inner}[2]{\langle #1, #2 \rangle}

\def\({\left(}
\def\){\right)}

\def\implies{\ \Rightarrow \ }

\newcommand{\ve}{\varepsilon}
\newcommand{\vp}{\varphi}
\newcommand*\diff{\mathop{}\!\mathrm{d}}
\newcommand\wt{\widetilde}
\newcommand\wh{\widehat}
\renewcommand{\leq}{\leqslant}
\renewcommand{\geq}{\geqslant}
\newcommand*\ball{\mathop{}\mathsf{B}}

\DeclareMathOperator{\sdr}{sdr}
\newcommand{\s}{\mathrm{s}}
\DeclareMathOperator{\rch}{rch}
\DeclareMathOperator{\wfs}{wfs}
\newcommand{\dsr}{\d_{\cS(r)}}
\DeclareMathOperator{\D}{D}

\newcommand{\R}{\mathbb{R}}
\newcommand{\E}{\mathbb{E}}
\newcommand{\X}{\mathbbm X}

\newcommand{\1}{\mathbbm 1}

\renewcommand{\dh}{\mathrm{d_H}}
\renewcommand{\d}{\mathrm{d}}
\newcommand{\Med}{\mathrm{Med}}

\newcommand{\tens}[2]{#1^{\otimes #2}}
\newcommand{\norm}[1]{\left\Vert #1\right\Vert}

\DeclareMathOperator{\support}{Support}

\newcommand{\manif}{\mathcal{C}}
\newcommand{\manifolds}[3]{\manif^{{#1}}_{{#2},{#3}}}

\newcommand{\distrib}{\mathcal{P}}
\newcommand{\distributions}[5]{\distrib^{{#1}}_{{#2},{#3}}({#4}, {#5})}

\newcommand{\circlett}{\mathsf{C}}

\usepackage{xcolor}
\definecolor{darkred}{RGB}{100,0,0}
\definecolor{darkgreen}{RGB}{0,100,0}
\definecolor{darkblue}{RGB}{0,0,150}
\definecolor{ccol}{RGB}{20, 143, 119 }

\usepackage{hyperref}
\hypersetup{colorlinks=true, linkcolor=darkred, citecolor=darkgreen, urlcolor=darkblue}
\usepackage{url}

\begin{document}

\thispagestyle{empty}

\title{
Optimal Reach Estimation and Metric Learning
}

\author{
	Eddie Aamari%
	\footnote{CNRS \& U. Paris \& Sorbonne U., Paris, France (\url{https://www.lpsm.paris/\~aamari/})}
	\and
	Cl\'ement Berenfeld%
	\footnote{CNRS \& U. Paris-Dauphine \& PSL, Paris, France (\url{https://www.ceremade.dauphine.fr/\~berenfeld/})}
	\and
	Cl\'ement Levrard%
	\footnote{CNRS \& U. Paris \& Sorbonne U., Paris, France (\url{http://www.normalesup.org/\~levrard/})}
}

\date{}
\maketitle

\begin{abstract}
We study the estimation of the \emph{reach}, an ubiquitous regularity parameter in manifold estimation and geometric data analysis.
Given an i.i.d. sample over an unknown $d$-dimensional $\mathcal{C}^k$-smooth submanifold of $\mathbb{R}^D$, we provide optimal nonasymptotic bounds for the estimation of its reach.
We build upon a formulation of the reach in terms of maximal curvature on one hand, and geodesic metric distortion on the other hand.
The derived rates are adaptive, with rates depending on whether the reach of $M$ arises from curvature or from a bottleneck structure.
In the process, we derive optimal geodesic metric estimation bounds.
\end{abstract}

\section{Introduction}
\subsection{Geometric Inference}

Topological data analysis and geometric methods now constitute a standard toolbox in statistics and machine learning~\cite{Wasserman18,Chazal21}.
In this family of methods, data $\bbX_n := \{X_1,\dots,X_n\}$ are usually seen as point clouds in high dimension, for which complex structural correlations give rise to an underlying structure that is neither full-dimensional, nor even linear.
Dealing with non-linearity is very well understood through the prism of \emph{non-parametric regression}.
However, in absence of distinguished ``covariate'' and ``response'' variables (i.e. coordinates), regression does not make sense anymore.
Hence, one needs to adopt a more global and coordinate-free approach: data are naturally viewed as lying on a submanifold $M \subset \R^D$ of dimension $d \ll D$, where $d$ corresponds to its \emph{true} number of degrees of freedom.

This approach opens the way to the estimation of numerous geometric and topological quantities to describe data.
Central to it is the manifold itself~\cite{Genovese12,Genovese12b,Kim15,Fefferman19,Divol21,Aizenbud21,Puchkin22}, where error is most commonly measured in Hausdorff distance.
Among many others, let us also mention the homology~\cite{Balakrishnan12}, persistent homology~\cite{Chazal14}, differential quantities~\cite{Aamari19}, intrinsic metric~\cite{Arias20} and regularity~\cite{Aamari19b}.

\subsection{Reach and Regularity}

Similarly to functional estimation, the theoretical study of nonparametric geometric problems naturally comes with regularity conditions.
By far, the most ubiquitous regularity and scale parameter in this context is the \emph{reach}.
First introduced by H. Federer's seminal paper~\cite{Federer59} on geometric measure theory, the reach $\rch(K) \in \R_+$ of a set $K \subset \R^D$ measures how far $K$ is from being convex~\cite{Attali13}. 
It hence provides a typical scale at which it shares most of the properties of a convex set.
These properties include -- among others -- uniqueness of the projection map, contractibility of balls, and explicit formulas for the volume of thickenings (see~\cite{Federer59}).
When $K = M$ is a submanifold, the reach also assesses quantitatively how it deviates from its tangent spaces. 
Therefore, the reach also provides an upper bound on curvature (that is, a bound in $\cC^2$) and a minimal scale of possible quasi self-intersections~\cite{Aamari19}.

For all these reasons, the reach practically appears in \emph{all} geometric inference methods as a natural scale parameter, which either drives a bandwidth used in a localization method~\cite{Genovese12,Aamari19b}, a minimal regularity scale in a minimax study~\cite{Kim15}, or a \emph{signal} part in a signal-to-noise ratio~\cite{Genovese12b,Fefferman19,Aizenbud21}. See~\cite{Aamari19b} for more examples of its use.
On the estimation side, the reach has already been studied under several angles.
\begin{itemize}
\item
The formulation of $\rch(M)$ in terms of deviation to tangent spaces from~\cite[Theorem~4.18]{Federer59} has been put to use through a plugin in~\cite{Aamari19}.
The authors derived non-matching upper and lower bounds for the estimation of $\rch(M)$ over $\cC^3$ submanifolds.
In addition to being suboptimal, the method of~\cite{Aamari19} requires the knowledge of tangent spaces, and is very sensitive to uncertainty on them (see~\cite[Section~6]{Aamari19}).
\item
Extending the minimax study of~\cite{Aamari19},~\cite{Berenfeld22} took advantage of the so-called \emph{convexity defect function} introduced by~\cite{Attali13} to propose another plugin strategy, with rates obtained over more general $\cC^k$-smooth manifold classes.
Despite still deriving non-matching upper and lower bounds, \cite{Berenfeld22} managed to exhibit two different estimation rates, depending on whether the reach testifies of a high curvature zone (the so-called local case, with slow rates) or of a narrow bottleneck structure (global case, with faster rates).
In this work, the derived rates are only suboptimal when the reach is achieved by curvature.
\item
More recently,~\cite[Theorem~1]{Boissonnat19} gave a new formulation of the reach in terms of geodesic distortion. Informally, they showed that $\rch(K)$ is the largest radius $r \geq 0$ for which the geodesic distance $\d_K$ is smaller than the geodesic distance $\dsr$ on a Euclidean ball of radius $r$.
Based on this purely metric statement, \cite{Cholaquidis21} proposed to plug-in a nearest-neighbor graph distance of the data in this formulation. This method provides a consistent estimator under very weak assumptions. Unfortunately, it fails to take advantage of high order regularity, when the reach is achieved by curvature (again).
\end{itemize}
With this analysis of possible estimation flaws in mind, this article proposes a two-step method. 
In short, we decouple the estimation of the local and global reaches \cite{Aamari19b}, and estimate them separately via max-curvature estimation and geodesic distance estimation respectively.

\subsection{Metric Learning}
In the data analysis area, \emph{metric learning} refers to the problem of finding a distance $\wh \d$ over the space of observations $\bbX_n \times \bbX_n$ that is relevant for a given task at stake~\cite{yang2006distance, suarez2021tutorial}. 
For instance, in a supervised framework where one is provided with tuples of allegedly \emph{similar} or \emph{dissimilar} observations, the goal is to find a distance that is small on the similar tuples and large on the dissimilar ones. 
There is a wide range of existing methods in the literature, ranging from parametric (LSI~\cite{xing2002distance}, MCML~\cite{globerson2005metric}, LDML~\cite{guillaumin2009you} among others) to nonparametric (DMLMJ~\cite{nguyen2017supervised}, kernel methods~\cite{kwok2003learning,chatpatanasiri2010new}, to cite a few).

In an unsupervised setting, metric learning aims at finding a metric that takes into account the underlying geometry of the data. 
That is, it amounts to estimating of shortest path (or \emph{geodesic}) distance.
Often, this is done via a dimension reduction technique: any low-dimensional embedding of the data gives rise to a new distance over the data in the embedded space.
Existing algorithms include PCA, t-SNE~\cite{hinton2002stochastic}, MDS~\cite{cox2008multidimensional}, Isomap~\cite{tenenbaum2000global}, or MVU~\cite{Arias14}.
See \cite{suarez2021tutorial} for a thorough overview of the field.

Astonishingly, despite the variety of existing methods, we are not aware of any general minimax study of geodesic metric learning. 
Though, two major theoretical references seem to stand out:
\begin{itemize}
\item
In \cite{Trillos2019}, the authors use a neighborhood graph to estimate distances and derive convergence rates in the $\cC^2$ case, but only for nearby points.
\item
In~\cite{Arias20}, estimation rates of geodesic distances are derived in the $\cC^2$ case using a reconstructing mesh.
Lower bounds are also obtained, but in a fixed-design setting only.
\end{itemize}
We propose a simple plugin method, and show that estimating the geodesic metric is no harder than estimating the manifold itself in Hausdorff distance.
This general strategy is also supported by a matching minimax lower bound.

\subsection{Contribution and Outline}

This article deals with the framework where data lies on an unknown $d$-dimensional $\cC^k$-submanifold of $\R^D$ (Section~\ref{sec:models}).
The main contribution consists of nearly-tight minimax bounds for reach estimation (Section~\ref{sec:optimal_reach_estimation}).
Along the way, three major building blocks, interesting in their own rights, are developed thoroughly:
\begin{itemize}
\item[(Section~\ref{sec:reach})]
We propose a general plug-in strategy for estimating the reach of a manifold. It is based on curvature estimation on one hand, and on the estimation of an intermediate scale (framed between the reach and the weak feature size) on the other hand.

\item[(Section~\ref{sec:sdr})]
We define the so-called \emph{spherical distortion radius} at scale $\delta > 0$ and study its estimation.
From the metric characterization of the reach from~\cite{Boissonnat19},
we notice that this purely metric quantity can be used to play the role of an intermediate scale for reach estimation.
We show that its stability properties make it well-suited to play the role of the intermediate scale of Section~\ref{sec:reach}.
\item[(Section~\ref{sec:metric})]
We propose a general plugin strategy for metric learning, and derive optimal geodesic metric estimation upper and lower bounds.
\end{itemize}
The proofs and the most technical points are deferred to the Appendix.

\subsection{General Notation}

In what follows, $\R^D$ ($D\geq 2$) is endowed with the Euclidean norm $\norm{\cdot}$.
The closed ball of radius $r\geq 0$ centered at $x \in \R^D$ is denoted by $\ball(x,r)$. If $x \in T \subset \R^D$ is a linear subspace, we write $\ball_T(x,r) := T \cap \ball(x,r)$ for the same ball in $T$.
Throughout, $c_\square,c'_\square,C_\square,C'_\square \geq 0$ denote generic constants that depend on $\square$, and that shall change from line to line to shorten notation. 
Similarly, universal constants shall generically be denoted by $c,c',C,C' \geq 0$.

\section{Geometric and Statistical Model}
\label{sec:models}
Let us first present the models in which we will work throughout.
As will be defined and discussed at length in Section~\ref{sec:reach_wfs_defi}, we let $\rch(K)$ denote the \emph{reach} of a subset $K \subset \bbR^D$ of the Euclidean space.

Building upon the standard regression setup, the following class is a good analog of Hölder classes of order $k \geq 2$, that is well adapted to submanifolds for stability reasons (see~\cite[Proposition~1]{Aamari19}).
Here, the analogy is to be understood as $T_p M$ being the (local) covariate space, and $\Psi_p$ being the regression function.

\begin{definition}[{\cite[Definition~1]{Aamari19b}}]
\label{def:geometric-model}
Let $k \geq 2$, $\rch_{\min}>0$, and $\mathbf{L} =  (L_2, L_3,\ldots,L_k)$.
We let $\manifolds{k}{\rch_{\min}}{\mathbf{L}}$ denote the set of $d$-dimensional compact connected submanifolds $M$ of $\R^D$ with $\rch(M) \geq \rch_{\min}$, such that for all $p \in M$, there exists a local one-to-one parametrization $\Psi_p$ of the form:
\begin{align*}
  \Psi_p \colon \ball_{T_pM}\left(0,r \right) &\longrightarrow M 
  \\
  v &\longmapsto p + v + \mathbf{N}_p(v)
\end{align*}
for some $r \geq \frac{1}{4 L_2}$,
with $\mathbf{N}_p \in \cC^{k}\left(\ball_{T_pM}\left(0, r \right) , \R^D\right)$ such that for all $\norm{v} \leq \frac{1}{4 L_2}$,
\begin{align*}
\mathbf{N}_p(0)=0,
\quad 
\diff_0 \mathbf{N}_p = 0
,\quad 
\text{and~}
\norm{\diff^j_v \mathbf{N}_p}_{\mathrm{op}} \leq L_j \text{ for all }j \in \{2,\ldots,k\},
\end{align*}
where $\diff^j_v \mathbf{N}_p$ stands for the $j$th differential of $\mathbf{N}_p$ at $v$, and $\norm{\cdot}_{\mathrm{op}}$ for the Euclidean operator norm over tensors.
\end{definition}

As explained in~\cite[Section~2.2]{Aamari19b}, radii $1/(4L_2)$ in local parametrizations have only been chosen for convenience. 
For $k=2$, the existence of parametrizations $\Psi_p$ is always guaranteed as soon as $\rch_{\min} > 0$ and $L_2 \geq 2/\rch_{\min}$ (see~\cite[Lemma~1]{Aamari19b}).

\begin{definition}
\label{def:statistical-model}
We let $\distributions{k}{\rch_{\min}}{\mathbf{L}}{f_{\min}}{f_{\max}}$ denote the set of Borel probability distributions $P$ on $\R^D$ satisfying:
\begin{itemize}
\item
Its support $M := \support(P)$ belongs to $\manifolds{k}{\rch_{\min}}{\mathbf{L}}$;
\item
It has a density $f$ with respect to the volume measure on $M$, such that 
$$
f_{\min} \leq f(x) \leq f_{\max}
~~\text{~for all }~
x \in M
.$$
\end{itemize}
\end{definition}

On the estimation side, the uniform smoothness of the parametrizations in Definition~\ref{def:geometric-model} allows for estimation of the manifold via local polynomial fitting around sample points in the models $\distributions{k}{\rch_{\min}}{\mathbf{L}}{f_{\min}}{f_{\max}}$.
Recall that the \emph{Hausdorff distance} between two compact subsets $K,K' \subset \R^D$ is defined by
\begin{align}
\label{eq:hausdorff}
\dh(K,K') 
:= 
\max\{
\sup_{x \in K} \d(x,K')
,
\sup_{x' \in K'} \d(x',K)
\},
\end{align}
where for all $u \in \R^D$,
\begin{align}
\label{eq:distance-function}
\d(u,K) := \min_{x \in K} \norm{x-u}
\end{align}
stands for the \emph{distance function} to $K$.
The estimation rates over the model $\distributions{k}{\rch_{\min}}{\mathbf{L}}{f_{\min}}{f_{\max}}$ have been studied in~\cite{Aamari19b}. A key result that we will use is the following.

\begin{theorem}[{\cite[Theorem~6]{Aamari19b}}] \label{thm:hatm}
There exists an estimator $\wh M$ such that for $n$ large enough,
\label{thm:hausdorff_estimation_Ck_expectation}
\begin{align*}
\sup_{ P \in \distributions{k}{\rch_{\min}}{\mathbf{L}}{f_{\min}}{f_{\max}}} 
\E_{P^n} [ \dh(\wh M,M) ]
\leq
C_{d,k,\rch_{\min},\mathbf{L},f_{\min},f_{\max}} 
\left ( \frac{\log n}{n} \right )^{k/d},
\end{align*}
where in the supremum, $M$ stands for $\support(P)$.
\end{theorem}
This rate is minimax optimal up to $\log n$ factors~\cite[Theorem~7]{Aamari19b}. 
It can be achieved by a local polynomial patch estimator $\wh M$ (see \eqref{eq:polynomial_expansion} below) that we will use as a preliminary step towards reach estimation. 
Let us also mention here that these fitted local polynomials also allow for estimation of differential quantities of $M$, such as tangent spaces and curvature at sample points, with (minimax) convergence rates of order $O(n^{-(k-1)/d})$ and $O(n^{-(k-2)/d})$ respectively (see~\cite[Theorems~2 to 5]{Aamari19b}). 
This fact will be of key importance in Section~\ref{sec:plug-in-for-reach}, where estimating the maximal curvature of $M$ will allow to estimate the so-called ``local reach''.

\section{Reach and Related Quantities}
\label{sec:reach}
\subsection{Characterizations and Relaxations of the Reach}\label{sec:reach_wfs_defi}

Let $K$ be a compact subset of $\R^D$. Following the original definition of~\cite{Federer59}, the \emph{reach} of $K$, denoted by $\rch(K)$, may be thought of as the largest radius of a neighborhood of $K$ onto which the projection map $\pi_K$ onto $K$ is well-defined. More formally, define the \emph{medial axis} of $K$ by
\begin{align*}
\Med(K) := \{ u \in \R^D \mid \exists x_1 \neq x_2 \in K \quad \|u - x_1\| = \|u-x_2\| = \d(u,K) \}. 
\end{align*}
The reach of $K$ is then defined as the smallest distance between $K$ and $\Med(K)$.
\begin{definition}\label{defi:reach}
For all closed $K \subset \R^D$, the \emph{reach} of $K$ is defined by
\begin{align*}
\rch(K) 
:= 
\min_{x \in K} \d(x,\Med(K)) 
= 
\inf_{u \in \Med(K)} \d(u,K)
.
\end{align*}
\end{definition}  
Note that in full generality, the medial axis might not be a closed set, so that the infimum in Definition~\ref{defi:reach} may not be attained (for instance in the case where $K$ is one-dimensional with a sharp edge). From a topological viewpoint, a key property of sets with positive reach is that the projection onto $K$ induces continuous retractions from the \emph{offset} $K^r := \{ u \in \R^D \mid \d(u,K) \leq r\}$ onto $K$, whenever $r < \rch(K)$~\cite[Theorem 4.8]{Federer59}. 
This property is at the core of topologically consistent reconstruction procedures such as that of~\cite{Boissonnat14}.  

Sets with positive reach can also been thought of as generalizations of convex sets, characterized by the smoothness of their distance function. Indeed, based on the remark that $x \mapsto \d(x,K)$ is $\cC^1$ on $\R^D \setminus K$ whenever $K$ is convex,~\cite{Clarke95} define $r$-\textit{proximally-smooth} sets as the sets $K$ such that $\d(\cdot,K)$ is $\cC^1$ over $\{u \in \R^D \mid 0 < \d(u,K) < r\}$. Interestingly, for subsets of $\R^D$, $r$-proximally smooth sets are exactly sets with reach $\rch(K) \geq r$~\cite{Poliquin00}, so that the reach may be alternatively defined in terms of gradients of the distance function. 
To this aim, following~\cite{Chazal05}, a generalized gradient function can be defined over $\R^D \setminus K$.
For all $x \in \R^D \setminus K$, we write
\begin{align}\label{eq:defi_gradient_distance}
\nabla \d(x,K) := \frac{x-c_K(x)}{\d(x,K)},
\end{align}
where $c_K(x)$ is the center of the smallest enclosing ball of the set $\pi_K(\{x\})$ of nearest neighbors of $x$ on $K$.
Since $c_K(x) = \pi_K(x)$ whenever $x \notin \Med(K)$, the medial axis can actually be characterized as
\begin{align*}
\Med(K) = \{ x \in \R^D \setminus K \mid \| \nabla \d(x,K) \| < 1 \},
\end{align*}
and the reach as 
$$
\rch(K)
=
\sup \{ r >0 \mid 0 < \d(x,K) < r \Rightarrow \| \nabla \d(x,K) \| =1 \}
.$$ 
This characterization of the reach allows for a straightforward relaxation. Namely, for a parameter $\mu \in [0,1]$, the seminal paper~\cite{Chazal06} introduces the so-called \emph{$\mu$-medial axis} as being 
\begin{align*}
\Med_\mu(K) := \{ x \in \R^D \setminus K \mid \| \nabla \d(x,K) \| \leq \mu \}, 
\end{align*}
and the $\mu$-reach as 
\begin{align}
\label{eq:mu-reach}
\rch_\mu(K) := \inf_{u \in \Med_\mu(K)} \d(u,K). 
\end{align}
It is clear that for all $\mu<1$, $\rch(K) \leq \rch_\mu(K)$, with $\rch(K)$ corresponding to the limit $\rch_{1^-}(K)$. 
Furthermore, this relaxation of the reach still yields enough regularity guarantees that the offsets $K^r = \{ u \in \R^D \mid \d(u,K) \leq r\}$ are isotopic for all $r \in (0, \rch_\mu(K))$~\cite[Lemma~2.1]{Chazal06}. 
Hence, the condition that $\rch_\mu(K) > 0$ conveys enough regularity properties for many topological estimators to work~\cite{Chazal09}.

Through this lens, the largest radius that ensures the topological stability of the offsets is the $0$-reach, also called \textit{weak-feature size},
\begin{align}\label{eq:wfs_defi}
\wfs(K) := \inf_{ u \in \Med_0(K) } \d(u,K),
\end{align}
that is the distance from $K$ to the set of critical points of $\d(\cdot,K)$. As detailed in the following section, the weak-feature size plays a special role in the case where $K$ is a manifold.
Here come a few elementary properties of the weak feature size that we will use later on.
\begin{proposition}\label{prop:wfs_properties}
Let $K \subset \R^D$ be compact. 
\bitem
\item[(i)] If $K$ is a closed submanifold of $\R^D$, then $\wfs(K) < +\infty$;
\item[(ii)] If $\wfs(K) < +\infty$, then for all $\mu \in [0,1)$,
\begin{align*}
\rch(K) 
\leq 
\rch_\mu(K)
\leq
\wfs(K) 
\leq 
\sqrt{\frac{D}{2(D+1)}}
\diam(K)
.
\end{align*}
\eitem
\end{proposition}

A proof is given in Section \ref{sec:proof_of_prp_wfs_properties}. Proposition~\ref{prop:wfs_properties} thus ensures that $\wfs(M)$ is uniformly bounded over the classes $\manifolds{k}{\rch_{\min}} {\mathbf{L}}$ introduced in Section~\ref{sec:models}. 
Since $\wfs(K)$ and $\rch_\mu(K)$ both measure a typical scale for topological stability, estimating them from sample could be of practical interest for topological inference.
Unfortunately, the following negative result shows that this estimation problem is intractable, even over a well-behaved model of closed $\cC^k$-submanifolds such as $\distributions{k}{\rch_{\min}}{\mathbf{L}}{f_{\min}}{f_{\max}}$.

\begin{theorem}\label{thm:wfs_inconsistency}
\label{thm:mu-reach-inconsistency}
Assume that $f_{\min} \leq c_{d,k}/\rch_{\min}^d$ and $f_{\max} \geq C_{d,k}/\rch_{\min}^d$, and $L_j \geq C_{d,k}/\rch_{\min}^{j-1}$ for all $j \in \{2,\ldots,k\}$.
Then there exists $\tilde{c}_{d,k}>0$ such that for all $n \geq 1$ and $\mu \in [0,1)$,
\begin{align*}
\inf_{\widehat{r}_\mu} \sup_{P \in \distributions{k}{\rch_{\min}}{\mathbf{L}}{f_{\min}}{f_{\max}}} 
\E _{P^{\otimes n}}\left[
| \widehat{r}_\mu - \rch_\mu(M)| 
\right]
\geq
\tilde{c}_{d,k}\rch_{\min}
>
0
,
\end{align*}
where $\widehat{r}_\mu$ ranges among all the possible estimators based on $n$ samples.
\end{theorem}
An intuition behind Theorem \ref{thm:mu-reach-inconsistency} is that for all $\mu < 1$, the $\mu$-medial axis is an unstable structure.
For certain manifolds $M_0 \in \manifolds{k}{\rch_{\min}}{\mathbf{L}}$, one can find arbitrarily small perturbations of $M_0$ whose $\mu$-medial axes remain at a fixed Hausdorff distance from $\Med_\mu(M_0)$. 
See the proof of Theorem \ref{thm:mu-reach-inconsistency} in Section~\ref{sec:proof-of-lower-bound-mu-reach} for a precise statement of this intuition.

Despite the fact that $\rch(K) = \rch_{1^-}(K)$, this negative result indicates that we cannot leverage $\mu$-reach estimation to obtain quantitative bounds for reach estimation.
We shall hence turn towards other reach-related quantities. In fact, the particular case where $K= M$ is a manifold offers us several other characterizations of the reach, which suggest other estimation strategies.

\subsection{Reach of Submanifolds}
\label{sec:reach-of-submanifolds}

In what follows, $M$ stands for a $d$-dimensional closed submanifold of $\R^D$. 
Note that~\cite[Remarks~4.20 and~4.21]{Federer59} and~\cite{Boissonnat19} assert that a closed submanifold with positive reach is at least of regularity $\mathcal{C}^{1,1}$, so that geodesics and tangent spaces are always defined in the usual differential sense. 
For the manifold case, the intuition of $\rch(M)$ as a generalized convexity parameter is further backed by~\cite[Theorem~4.8]{Federer59}. 
Indeed, the inequality $\left\langle x - \pi_C(x), \pi_C(x) - c \right\rangle \geq 0$ valid for all $c \in C$ and $x \in \R^D$ whenever $C$ is convex, translates to $\left\langle x - \pi_M(x), \pi_M(x) - y \right\rangle \geq - \|\pi_M(x)-y\|^2 \|x-\pi_M(x)\|/(2\rch(M))$ being valid for all $y \in \R^D$ and $x \in \R^D$ such that $\d(x,M) < \rch(M)$. 
This leads to the following characterization of the reach, in the manifold case.

\begin{theorem}[{\cite[Theorem 4.18]{Federer59}}]\label{thm:Fed_reach_tangentspace}
For a submanifold $M \subset \R^D$ without boundary, 
\[
\rch(M) = \inf_{p \neq q \in M} \frac{\|p-q\|^2}{2 \d(q-p,T_p M)},
\]
where $T_p M$ denotes the tangent space of $M$ at $p$.
\end{theorem} 
This result provides a natural plugin estimator, proposed by \cite{Aamari19b}, which consists in replacing $M$ and $T_p M$ by suitable estimators of them.
A key result from~\cite{Aamari19b} is a description of how the infimum in Theorem~\ref{thm:Fed_reach_tangentspace} is achieved, possibly asymptotically.
\begin{theorem}[{\cite[Theorem 3.4]{Aamari19b}}]
\label{thm:reach_wfs_local}
Let $M \subset \R^D$ be a compact $\mathcal{C}^2$ submanifold without boundary. 
Then,
\begin{align*}
\rch(M) = \wfs(M) \wedge R_{\ell}(M),
\end{align*}
where denoting by $\II_p: T_p M \times T_p M \to T_p M^\perp$ the \emph{second fundamental form} of $M$ at $p \in M$,
$$R_{\ell}(M) := \min_{p \in M} \| \II_p\|_{\mathrm{op}}^{-1}
$$ stands for the minimal curvature radius of $M$.
\end{theorem}

This result conveys the following intuition in the manifold case: the infimum in the right-hand side of Theorem~\ref{thm:Fed_reach_tangentspace} may be attained:
\begin{itemize}
\item[(Local case)]
Asymptotically, for pairs of points $(p,q)$ converging to a maximal curvature point in some direction, so that $\rch(M) = R_\ell(M)$.
\item[(Global case)]
For a pair of points $(p,q)$ belonging to parallel areas of $M$, forming a bottleneck zone, so that $\rch(M) = \wfs(M)$.
\end{itemize}
This local/global dichotomy of the reach may also be retrieved in the recent characterization given by~\cite{Boissonnat19} in terms of metric distortion.

\begin{theorem}[{\cite[Theorem 1]{Boissonnat19}}]\label{thm:reach_characterization_metric_distortion}
Let $K \subset \R^D$ be a closed subset. 
Then
\begin{align*}
\rch(K) = \sup \left\lac {r >0 \mid \forall p,q \in K, \|p-q\| < 2r \Rightarrow \d_K(p,q) \leq 2r \arcsin \left ( \frac{\|p-q\|}{2r} \right )} \right\rac
,
\end{align*}
where $\d_K: K \times K \to \bar{\R}_+$ stands for the \emph{shortest-path} (or \emph{geodesic}) distance on $K$.

\end{theorem}
Recall that, for all $p,q \in K$, the distance $\d_K(p,q)$ is the infimum of the length of all the continuous path in $K$ between $p$ and $q$.
As will be detailed in Section~\ref{sec:sdr}, the above result allows to characterize the reach in terms of metric distortion with respect to metrics on spheres of radii $r$. 
In the same spirit as Theorem~\ref{thm:reach_wfs_local}, when $K=M$ is a submanifold, the configurations of $(p,q,r)
$ in the supremum of Theorem~\ref{thm:reach_characterization_metric_distortion} are limited by the same two local and global layouts: 
\begin{itemize}
\item[(Local case)]
When $p$ and $q$ tend to a maximal curvature point in some direction, the geodesic distance $\d_K$ behaves like that of a sphere of radius $R_{\ell}(M)$ at this point in this direction.
\item[(Global case)]
When $p$ and $q$ are in parallel areas, their geodesic distance must be larger than the spherical distance of radius $\|p-q\|/2$. 

\end{itemize}
\subsection{Plug-in Methods for Reach Estimation}
\label{sec:plug-in-for-reach}

The characterizations of the reach given in Section~\ref{sec:reach-of-submanifolds} all lead to their associate plug-in estimators:
\begin{itemize}
\item
Studying a $\cC^3$ model similar to $\distributions{3}{\rch_{\min}}{\mathbf{L}}{f_{\min}}{f_{\max}}$,~\cite{Aamari19b} took advantage of the characterization with tangent spaces (Theorem~\ref{thm:Fed_reach_tangentspace}) to conceive a reach estimator that converges at rate $O(n^{-2/(3d-1)})$ in the local case ($\rch(M) = R_{\ell}(M)$), and $O(n^{-1/d})$ in the global case ($\rch(M) = \wfs(M)$).
\item
Based on the metric distortion characterization of Theorem~\ref{thm:reach_characterization_metric_distortion},~\cite{Cholaquidis21} propose a reach estimator that is consistent whenever $M$ has positive reach.
\end{itemize}
In light of Theorem~\ref{thm:reach_wfs_local}, differences of convergence rates between the local and global case are to be expected.
To quantify this intuition,~\cite{Berenfeld22} introduces subclasses of the model $\distributions{k}{\rch_{\min}}{\mathbf{L}}{f_{\min}}{f_{\max}}$, parametrized by the gap between $R_\ell(M)$ and $\wfs(M)$. They obtain the following lower bounds.

\begin{theorem}[{\cite[Theorem 7.1]{Berenfeld22} and~\cite[Proposition 2.9]{Aamari19b}}]
\label{thm:lwr_bounds_clementb}
Let $\alpha \in \R$, $k \geq 2$, and write
\begin{align*}
\distributions{k}{\rch_{\min}}{\mathbf{L}, \alpha}{f_{\min}}{f_{\max}} 
:= 
\left\lac  P \in \distributions{k}{\rch_{\min}}{\mathbf{L}}{f_{\min}}{f_{\max}} \mid R_{\ell}(M) \geq \wfs(M) + \alpha \right\rac,
\end{align*}
where $M$ denotes $\support(P)$. Then, for all $\rch_{\min} > 0$ there exists small enough $f_{\min}$ and large enough $f_{\max}$, $\mathbf{L}$ such that
\begin{align*}
\inf_{\widehat{\rch}} \sup_{P \in \distributions{k}{\rch_{\min}}{\mathbf{L}, \alpha}{f_{\min}}{f_{\max}}} \E | \widehat{\rch} - \rch(M)| & \geq c_{\rch_{\min},d,k} \left(\frac{1}{n} \right)^{(k-2)/d}, \quad \mbox{\text{if} $\alpha \leq 0$}, \\
\inf_{\widehat{\rch}} \sup_{P \in \distributions{k}{\rch_{\min}}{\mathbf{L}, \alpha}{f_{\min}}{f_{\max}}} \E | \widehat{\rch} - \rch(M)| & \geq c_{\rch_{\min},d,k, \alpha} \left(\frac{1}{n} \right)^{k/d}, \quad \mbox{\text{if} $\alpha > 0$}.
\end{align*}
\end{theorem} 
These bounds indicate that estimating the reach is at least as hard as estimating the curvature in the local case ($\rch(M)=R_{\ell}(M)$), and at least as hard as estimating the manifold in the global case ($\rch(M) = \wfs(M)$).
We will prove in Section~\ref{sec:optimal_reach_estimation} that these rates are in fact minimax optimal up to $\log n$ factors.
This means that reducing reach estimation to curvature \emph{and} manifold estimation is a good way to go, as it leads to optimal rates.
To do so, following the idea behind Theorem~\ref{thm:reach_wfs_local}, estimating $R_{\ell}(M)$ -- or some notion of local reach -- and $\wfs(M)$ -- or some notion of global reach -- separately seems a sensible approach.

\subsubsection{Local Reach Estimation}
For (max-)curvature estimation, the strategy that we adopt follows from the polynomial patches estimator proposed in~\cite{Aamari19}.
Given a localization bandwidth $h>0$, and a parameter $t > 0$, for all $i \in \{1,\ldots,n\}$,
we let $\hat{\pi}_i : \R^D \to \R^D$ be an orthogonal projector of rank $d$ and $\hat{\mathbb{T}}^{(j)}_{i} : \bigl( \R^D \bigr)^{\otimes j} \to \R^D$ be symmetric tensors solutions of the least squares problem
\begin{align}\label{eq:defi_pol_fit}
\min_{
\substack{ 
\pi
\\
\max_{2\leq j \leq k-1} \|\mathbb{T}^{(j)}\|^{\frac{1}{j-1}} \leq t
}
} P_{n-1}^{(i)} \left [ \left \| x - \pi(x) - \sum_{j=2}^{k-1} \mathbb{T}^{(j)}(\pi(x)^{\otimes j}) \right \|^2 \1_{\ball(0,h)}(x) \right],
\end{align}
where $P_{n-1}^{(i)} := \frac{1}{n-1} \sum_{p \neq i} \delta_{X_p-X_i}$ denotes the empirical measure centered  at point $X_i$. 
Following~\cite[Section 3]{Aamari19}, if $h$ is taken to be of order $\Theta\( (\log n / n)^{1/d} \)$, that $t$ is chosen such that $t^k h \leq 1$, and that $\hat{T}_i := \Im(\hat{\pi}_i)$ denotes the image of $\hat{\pi}_i$  -- which is a $d$-dimensional vector space by construction --, then the local patches 
\begin{align}
\label{eq:polynomial_expansion}
  \widehat{\Psi}_i \colon \ball_{\hat{T}_i}(0,7h/8) &\longrightarrow \R^D \notag
  \\
  v &\longmapsto X_i + v + \sum_{j=2}^{k-1} \hat{\mathbb{T}}_{i}^{(j)}(\tens{v}{j})
\end{align}
are local $O(h^k)$ approximations of $M$ whenever $n$ is large enough. 
Furthermore, for $v \in \ball_{\hat{T}_i}(0,h/4)$, we can estimate the curvature tensor at $\pi_M(\widehat{\Psi}_i(v))$ via the second derivative of $\widehat{\Psi}_i$ at $v$, expressed in local coordinates around $\widehat{\Psi}_i(v)$ given by a basis of $\mathrm{Im}(\d_v \widehat{\Psi}_i)$. 
To summarize, for all $v \in \ball_{\hat{T}_i}(0,h/4)$,~\eqref{eq:polynomial_expansion} provides a $d$-dimensional space $\hat{T}_{i,v} := \mathrm{Im}(\d_v \widehat{\Psi}_i)$, as well as a symmetric bilinear map 
\begin{align*}
\hat{\mathbb{T}}_{i,v}^{(2)}: \hat{T}_{i,v} \times \hat{T}_{i,v} \rightarrow \hat{T}_{i,v}^\perp, 
\end{align*} 
that is provably close to $\II_{\pi_M(\widehat{\Psi}_i(v))}$. The precise definition of $\hat{\mathbb{T}}_{i,v}^{(2)}$ is given in Section \ref{sec:proof_thm_cvrates_curvature_max}. A minimal curvature radius (i.e. maximal curvature) estimator may then be computed as the minimal curvature radius of all the polynomial patches around sample points, that is
\begin{align}\label{eq:defi_curvature_radius_estimate}
\wh R_\ell := \min_{1 \leq i \leq n} \min_{ v \in \ball_{\hat{T}_i}(0,h/4)} \|\hat{\mathbb{T}}_{i,v}^{(2)}\|_{\mathrm{op}}^{-1}.
\end{align} 
Provided $M$ is uniformly well approximated by $\bigcup_{i=1}^n \widehat{\Psi}_i(\ball_{\hat{T}_i}(0,h/4))$, the convergence rate of $\wh R_\ell$ towards $R_\ell(M)$ will follow from uniform curvature bounds, similar to the pointwise ones from~\cite[Theorem~4]{Aamari19}.
We are able to prove the following.

\begin{theorem}\label{thm:cv_rates_curvature_max}
Let $k \geq 3$ and $P \in \distributions{k}{\rch_{\min}}{\mathbf{L}}{f_{\min}}{f_{\max}}$. Write $h = \left ( C_{d,k} \frac{f^2_{\max} \log n}{f^3_{\min}n} \right )^{1/d}$.
Then for $n$ large enough, with probability larger than $1- 2n^{-k/d}$, we have
\begin{align*}
\bigl| \wh R_\ell - R_{\ell}(M) \bigr| \leq C_{d,k,\mathbf{L}, \rch_{\min}} R_{\ell}^2(M) \sqrt{\frac{f_{\max}}{f_{\min}}} h^{k-2}.
\end{align*}
\end{theorem}
We refer to \secref{proof_thm_cvrates_curvature_max}  for a proof of this result. In particular, the estimator $\wh R_\ell$ achieves the rate of the lower bound from Theorem~\ref{thm:lwr_bounds_clementb} in the case where $\rch(M) = R_{\ell}(M)$ (i.e. $\alpha \leq 0$), up to $\log n$ factors. 

\subsubsection{Global Reach Estimation}
To complete the construction of an estimator of $\rch(M)$, building an estimator of $\wfs(M)$ could be a possibility.
However, Theorem~\ref{thm:wfs_inconsistency} shows that building an estimator of the weak feature size with a uniform convergence rates over $\distributions{k}{\rch_{\min}}{\mathbf{L}}{f_{\min}}{f_{\max}}$ is hopeless. 
Nonetheless, it is important to note that a uniform estimation rate of $\wfs(M)$ over $\mathcal{P}^k$ is not necessary to obtain uniform convergence rate for $\rch(M)$.
Indeed, an estimator $\widehat{\wfs}$ of $\wfs(M)$ that exhibits an optimal uniform convergence rate whenever $\wfs(M) \leq R_{\ell}(M)$, and that is provably larger than $R_{\ell}(M)$ otherwise, is enough to build an optimal reach estimator when combined with $\wh R_\ell$.
This is the case, for instance, of the weak feature size estimator of~\cite{Berenfeld22} based on the so-called \emph{convexity defect function}.

Based on this remark, we adopt a more general strategy, by seeking for an intermediate geometric scale $\theta(M)$ (or \emph{feature size}) such that for all $M \in \manifolds{k}{\rch_{\min}}{\mathbf{L}}$,
$$
\rch(M) \leq \theta(M) \leq \wfs(M)
.
$$
In such a case, Theorem~\ref{thm:reach_wfs_local} extends trivially, with $\wfs(M)$ replaced by $\theta(M)$.
\begin{proposition}
\label{prop:reach_intermediate_local}
Assume that $\theta : \cC^{2}_{\rch_{\min}}{} \to \R_+$ is such that $\rch(M) \leq \theta(M) \leq \wfs(M)$ for all $M \in \cC^{2}_{\rch_{\min}}{}$.
Then,
\begin{align*}
\rch(M) = \theta(M) \wedge R_{\ell}(M).
\end{align*}
\end{proposition}

Given such an intermediate scale parameter of interest $\theta(M)$, and assuming that a consistent estimator $\wh \theta $ of $\theta(M)$ is available, one can naturally consider the plugin $\widehat{\rch} := \wh{R}_\ell \wedge \wh{\theta}$. For free, Proposition~\ref{prop:reach_intermediate_local} yields that $\theta(M)\1_{R_{\ell}(M) > \rch(M)} = \rch(M)\1_{R_{\ell}(M) > \rch(M)}$, so that
\begin{align}
\label{eq:double-plugin-precision}
|\rch(M) - \widehat{\rch}| \leq | \wh{R}_\ell - R_{\ell}(M)| \1_{R_{\ell}(M) \leq \rch(M)} + | \wh{\theta} - \theta(M)| \1_{R_{\ell}(M) > \rch(M)}
,
\end{align}
as soon as $|R_{\ell}(M) - \wh{R}_\ell| + |\theta(M) - \wh{\theta}| \leq |R_{\ell}(M) - \theta(M)|$.
In addition, such a quantity would provide a local scale that is of interest for further topological inference, as exposed in Section~\ref{sec:reach_wfs_defi}.  

According to Theorem~\ref{thm:wfs_inconsistency}, taking $\theta(M)$ to be related to the medial axis characterization of the reach -- such as the $\mu$-reach, or the $\lambda$-reach defined in~\cite{Chazal05}) -- is likely to lead to an unsolvable statistical problem, because of the inherent instability of the medial axis.
Hence, we rather build upon the metric distortion characterization of the reach given by Theorem~\ref{thm:reach_characterization_metric_distortion}, and provide a better-behaved intermediate scale $\theta(M)$: the \emph{spherical distortion radius}.

\section{Spherical Distortion Radius} 
\label{sec:sdr}

\subsection{Motivation and Definition}
Based on Theorem~\ref{thm:reach_characterization_metric_distortion}, 
we now build a geometrically stable feature size that measures the maximum radius (or scale) at which the geodesic distance can be compared to the corresponding spherical distance.
To be more precise, for $x,y \in \bbR^D$ and $r > 0$, we define the \emph{spherical distance} $\dsr(x,y)$ -- or \emph{great-circle distance} -- as the distance between $x$ and $y$ when seen as both lying on a sphere of radius $r$. That is,
\[
\dsr(x,y) 
:= 
\begin{cases}
2 r \arcsin\left(\frac{\|x-y\|}{2r}\right)
&
\text{if } \norm{x-y} \leq 2r,
\\
+\infty
&
\text{otherwise}
\end{cases}
\]
Note that the map $r \mapsto \dsr(x,y)$ is decreasing on $[\|x-y\|/2,\infty)$ and that
$$
\dsr(x,y) = \frac12 \pi \|x-y\| ~~\text{for}~~r = \frac{\|x-y\|}{2}~~~\text{and}~~~\dsr(x,y) \xrightarrow[r \to \infty]{} \|x-y\|.
$$ 
Then, Theorem~\ref{thm:reach_characterization_metric_distortion} can be rewritten as 
\begin{align*}
\rch(K) = \sup \left \lac r >0 \mid \forall x,y \in K, \|x-y\| < 2r \Rightarrow \d_K(x,y) \leq \dsr(x,y) \right\rac.
\end{align*}
It should be noted that $\dsr$ is not formally a distance on $K$ (unless $K$ is a subset of a sphere of radius $r$), but this is of little importance in what follows.

Based on the same idea that motivates the introduction of the $\mu$-reach, we intend to discard curvature effects to obtain some notion of global reach.
In the metric characterization of the reach from Theorem~\ref{thm:reach_characterization_metric_distortion}, this can be done by the supremum restricting to points that are not too close.
\begin{definition}\label{defi:SDR}
Let $K$ be a compact subset of $\R^D$, $\d$ a distance on $K$ and $\delta >0$. The \emph{spherical distortion radius} of the metric space $(K,\d)$ at scale $\delta$ is defined by
$$ 
\sdr_\delta(K,\d) := \sup\{r > 0~\middle|~\forall x,y \in K,~ \delta \leq \|x-y\| < 2r \implies \d(x,y) \leq \dsr(x,y) \}.
$$
\end{definition}
In words, the spherical distortion radius at scale $\delta > 0$ is the largest radius $r$ for which the distance $\d$ is bounded above by the spherical distance at radius $r$, when restricted to points that are at least $\delta$-apart for the Euclidean distance.

\begin{figure}[!htbp]
\centering
	\includegraphics[width = 0.5\linewidth]{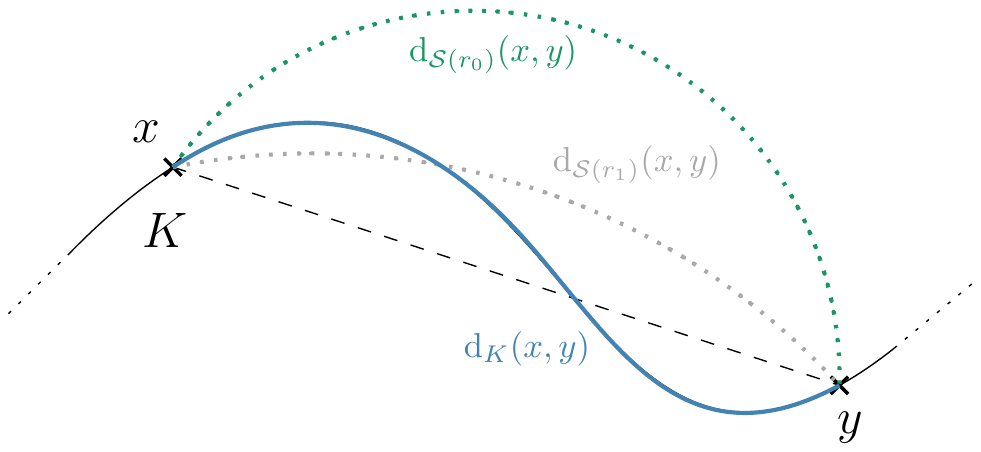}
\caption{
A curve $K$  in the plane. In blue is the shortest path between two points $x$ and $y$, whose length is $\d_K(x,y)$. 
In green (resp. grey) is the circle portion of radius $r_0$ (resp. $r_1$) going through $x$ and $y$. 
The layout is chosen so that $r_0 \leq r_1$ and $\d_{\cS(r_1)}(x,y) \leq \d_{K}(x,y) \leq \d_{\cS(r_0)}(x,y)$.
}
\label{fig:sdr}
\end{figure}

By construction, $\sdr_\delta(K,d) \geq\delta/2$ for all $\delta > 0$. 
Furthermore, whenever $\delta$ is strictly greater than $\diam K$, then no pairs of points in $x,y \in K$ satisfies $\|x-y\| \geq \delta$ so that $\sdr_\delta(K) = +\infty$. 
On the other hand, if $\delta = 0$, then the spherical distortion radius of $(K,\d_K)$, coincides with the reach of $K$ (Theorem~\ref{thm:reach_characterization_metric_distortion}). In fact, Proposition~\ref{prp:interpolate} below confirms that the spherical distortion radius interpolates between the reach and the weak feature size. 

\begin{proposition} \label{prp:interpolate}
For all closed $K \subset \bbR^D$ and all metric $\d$ on $K$, the map $\delta \mapsto \sdr_\delta(K,\d)$ is non-decreasing.
Furthermore, for  $\d = \d_K$,
$$\rch(K) \leq \sdr_\delta(K,\d_K) \leq \wfs(K) ~~~\text{for all} ~~~0 \leq \delta \leq \sqrt{\frac{2(D+1)}{D}}\wfs(K).$$ 
\end{proposition}
A proof of Proposition~\ref{prp:interpolate} is given in Appendix \ref{sec:proof-compare-reach-sdr}.

\begin{example}\label{ex:sdr_polygon} As a toy example, let us study the spherical distortion radius of the wedge shape $K_\alpha = \cL_1 \cup \cL_2$ where $\cL_1$ and $\cL_2$ are two half-line originated from a common point $z \in \bbR^D$ (see \figref{examplesdr}). We let $\alpha \in (0,\pi)$ be the angle between these two lines. In this context, we have $\rch(K_\alpha) = 0$, and it is easy to see that $\wfs(K_\alpha) = \infty$. Furthermore, the usual interpolations between the reach and the weak feature size exhibit a very degenerate behavior in the presence of an angular configuration such as this one, with for instance
$$
\rch_\mu(K_\alpha) = \begin{cases} 0~~&\text{if}~~ \mu \geq \sin(\alpha/2),\\
\infty ~~&\text{if}~~ \mu < \sin(\alpha/2).
\end{cases}
$$
On the contrary, we show hereafter that the spherical distortion radius interpolates non-trivially between $\rch(K_\alpha)$ and $\wfs(K_\alpha)$ in this case, giving rise to a new family of relevant characteristic scales even for non-smooth subsets $K_\alpha$. 

To see this, take $x \in \cL_1$ and $y \in \cL_2$, and denote by $a := \|x-z\|$ and $b := \|y-z\|$. The intrinsic distance $\d_{K_\alpha}(x,y)$ is given by $a+b$ while $\|x-y\|^2 = a^2 + b^2 - 2 ab \cos(\alpha)$. 
Now the solution of the minimization problem
$$
\min\{a^2+b^2-2ab \cos(\alpha)~|~ a + b = \d_{K_{\alpha}}(x,y)\} 
$$
is given by $a = b = \d_{K_{\alpha}}(x,y)/2$ and equals $\d^2_{K_{\alpha}}(x,y) \sin^2(\alpha/2)$. The spherical distortion radius of $K_\alpha$ at scale $\delta$ is thus the largest $r$ such that
\beq \label{exsdr}
\frac{\delta}{\sin(\alpha/2)} \leq 2r \arcsin\(\frac{\delta}{2r}\).
\eeq
Since the right-hand side above ranges between $\delta$ and $\delta\pi/2$, we distinguish two cases: 
\begin{itemize}
\item
If $\sin(\alpha/2) < 2/\pi$, then no $r$ can fulfill \eqref{exsdr}. Hence, $\sdr_\delta(K_\alpha,\d_{K_\alpha}) = \delta/2$.
\item
Otherwise $\sin(\alpha/2) \geq 2/\pi$, in which case the largest $r$ is given by the equality $\vp(2r/ \delta) = 1/\sin(\alpha/2)$, where $\vp(u) := u \arcsin(1/u)$ is a bijection between $[1,\infty)$ and $(1,\pi/2]$. 
\end{itemize}
All in all, it holds
$$
\sdr_{\delta}(K_\alpha,\d_{K_\alpha}) = \begin{cases}
\delta/2~~~~&\text{if}~~~\alpha < \alpha_* \\
(\delta/2) \vp^{-1}(1/\sin(\alpha/2)) ~~~~&\text{if}~~~\alpha \geq \alpha_*
\end{cases}
$$
where $\alpha_* = 2 \arcsin(2/\pi) < \pi/2$.
Note that compared to $\rch_\mu(K_\alpha)$, there is no discontinuity in $\sdr_\delta(K_\alpha,\d_{K_\alpha})$ as $\alpha$ varies.
\end{example}
\begin{figure}[!htbp]
\centering
\begin{subfigure}{0.49\textwidth}
	\centering
	\includegraphics[height = 0.5\linewidth]{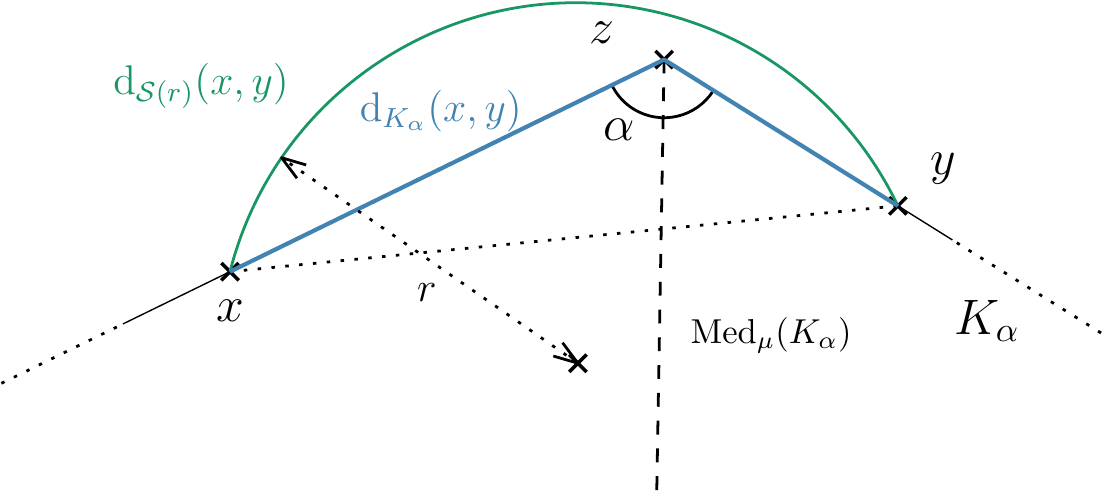}
	\caption{}
	\label{fig:examplesdr}
\end{subfigure}
\begin{subfigure}{0.49\textwidth}
	\centering
	\includegraphics[height = 0.5\linewidth]{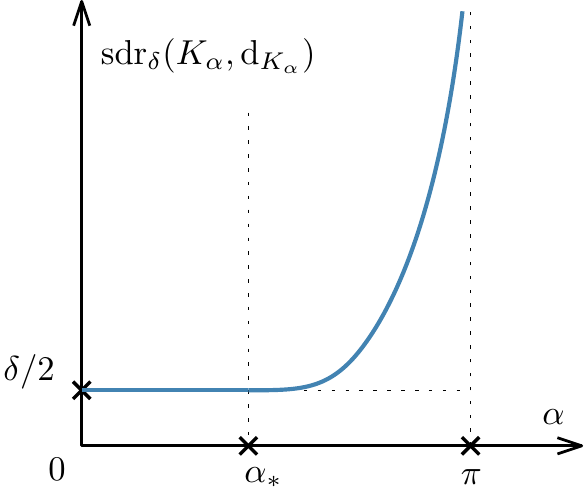}
	\caption{}
	\label{fig:examplesdr2}
\end{subfigure}
\caption{(\subref{fig:examplesdr}) Diagram of $K_\alpha = \cL_1 \cup \cL_2$ with an angle $\alpha$ between the two half-lines. The shortest path between $x$ and $y$ is drawn in blue. In dashed the $\mu$-medial axis for $\mu > \sin(\alpha/2)$, showing in particular that $\rch_\mu(K_\alpha) = 0$ in this case. 
(\subref{fig:examplesdr2}) Plot of the function $\alpha \mapsto \sdr_\delta(K_\alpha,\d_{K_\alpha})$, which operates a smooth interpolation between $\delta/2$ and $\infty$.}
\label{fig:examplesdr}
\end{figure}
Example~\ref{ex:sdr_polygon} above carries the intuition that the spherical distortion radius seems somehow stable with respect to Hausdorff perturbations, contrary to the $\mu$-reach. We quantify this intuition in the following section.

\subsection{Stability Properties}
In this section, we will be comparing different metric spaces on subsets of $\R^D$.
Let $K$ and $K'$ be two subsets of $\bbR^D$, endowed with distances $\d$ and $\d'$ respectively. 
We intend to prove that $\sdr_\delta(K,\d)$ and $\sdr_\delta(K',\d')$ are close whenever $(K,\d)$ and $(K',\d')$ are close, and that $(K,\d)$ has good properties.
The notion of proximity between $K$ and $K'$ will be measured in Hausdorff distance (see~\eqref{eq:hausdorff}). 
It remains to define a notion of proximity between $\d$ and $\d'$, which is called the \emph{mutual distortion}.

\begin{definition}\label{def:metric_distortion}
Let $(K,\d)$ and $(K',\d')$ be two metric subspaces of $\bbR^D$. The \emph{metric distortion} of $\d'$ relative to $\d$ at scale $\delta > 0$ is
 $$
 \D_\delta(\d' | \d) :=  \sup_{\substack{x',y' \in K' \\ \|x'-y'\| \geq \delta}} \frac{\d'(x',y')}{\d(\pr_K(\{x'\}), \pr_K(\{y'\}))}.
 $$
where $\pr_K$ is the (possibly multivalued) closest-point projection onto $K$ for the ambient Euclidean distance, and where  
$$
\d(\pr_K(\{x'\}), \pr_K(\{y'\})) 
:= 
\inf
\{\d(x,y)~|~x\in \pr_K(\{x'\}),~y \in \pr_K(\{y'\}) \}
.
$$ 
We adopt the convention $\D_\delta(\d'|d) = 0$ if $\delta > \diam(K')$. 
The \emph{mutual distortion} of $\d$ and $\d'$ is then defined as
\begin{align*}
\D_\delta(\d,\d') := \max\{\D_\delta(\d'|\d), \D_\delta(\d|\d')\}.
\end{align*}
\end{definition}
The mutual distortion defined above allows to compare distances on different spaces, while taking into account their respective embeddings in $\R^D$.
A small distortion $\D_\delta(\d,\d')$ means that, if $a,b \in K$ and $x,y \in K'$ are two couples of points that are $\delta$-separated and such that $x$ and $a$, and $y$ and $b$ are respectively close to each other, then $\d(a,b)$ and $\d'(x,y)$ should be close as well. This definition of \emph{mutual distortion} between metric subspaces of $\bbR^D$ is related to the existing notion \emph{metric distortion} of an embedding. 
See for instance~\cite{bourgain1985lipschitz} or more recently~\cite{chennuru2018measures} which deals with distortion measures in a statistical framework.
It is nonetheless significantly different, in particular because the usual notion of distortion is invariant through re-scaling of either $\d$ or $\d'$.  In our framework, invariance with respect to scaling is an undesirable property, since we want to estimate the reach, which is itself a scale factor (or feature size).
\begin{rem} \label{rem:bilip} When $K = K'$, the mutual distortion can be seen as the bi-Lipschitz coefficient of $\Id : (K,\d) \to (K,\d')$ at \emph{scale} $\delta$, meaning that for all $x,y \in K$
$$
\|x-y\| \geq \delta \implies  \frac1L \d'(x,y) \leq \d(x,y) \leq L \d'(x,y),
$$
where $L = \D_\delta(\d,\d')$. In particular, a mutual distortion that is close to $1$ means that $(K,\d)$ is quasi-isometric to $(K,\d')$, at scale $\delta$.
\end{rem}
If the two subspaces $K$ and $K'$ are too far apart, then it makes no sense to compare two distances $\d$ and $\d'$ defined on them, and one could expect the mutual distortion to explode. This is will typically the case when $\dh(K,K') \geq \delta$.

It is clear from the definition that using the notion of relative metric distortion defined above, the spherical distortion radius of $K$ may be expressed as
$$
\sdr_\delta(K,\d) = \sup\{r > 0~\middle|~\D_\delta(\d | \dsr) \leq 1\}.
$$
This point supports the idea that the relative metric distortion we defined is a suitable notion of proximity to assess stability of the spherical distortion radius, as exposed by the following proposition.

\begin{proposition} \label{prp:stab1}
Let $\delta_0 > 0$ and $\ve,\nu > 0$. Assume that both $\dh(K',K) \leq \ve$ and $\D_{\delta_0}(\d'|\d) \leq 1+\nu$.
Define
$$
 \xi(r) := 384 (1+\pi)\frac{r^4}{\delta_0^4}
 \text{~for all~} r \geq 0.
$$
Then, for all $\delta \geq \delta_0$, letting $\Upsilon := (\delta \nu)\vee \ve$ and $\mathrm{r}_{1} := \sdr_{\delta+2\ve}(K',\d')$, if $ \xi(\mathrm{r}_{1}) \Upsilon  < \mathrm{r}_{1}$, then
$$
\sdr_{\delta}(K,\d) \leq \sdr_{\delta+2\ve}(K',\d')+\xi(\mathrm{r}_{1}) \Upsilon. 
$$
\end{proposition}
A proof of \prpref{stab1} is given in Appendix \ref{sec:proof-sdr-properties}.
Note that the condition $\dh(K',K) \leq \ve$ may be relaxed via $\dh(K'|K) \leq \ve$, where $\dh(K'|K) := \sup_{x \in K'} \d(x,K)$. Also, under the assumptions of \prpref{stab1}, let us remark that if $ \sdr_{\delta+2\ve}(K',\d')$ is finite, then so is $\sdr_{\delta}(K,\d)$ with $\sdr_{\delta}(K,\d) \leq 2 \sdr_{\delta+2\ve}(K',\d')$. 

Proposition~\ref{prp:stab1} can be symmetrized to get the following two-sided control.

\begin{corollary} \label{cor:stab} Let $0 < \delta_0 < \delta_1$ and $\ve,\nu > 0$. Assume that both $\dh(K',K) \leq \ve$ and $\D_{\delta_0}(\d',\d) \leq 1+\nu$. Then, for any $\delta \in (\delta_0+2\ve,\delta_1 - 2\ve)$, it holds
$$
\sdr_{\delta-2\ve}(K,\d)-\xi_0 \Upsilon \leq \sdr_{\delta}(K',\d') \leq \sdr_{\delta+2\ve}(K,\d)+\xi_0 \Upsilon 
$$
with $\xi_0 := \xi(2\sdr_{\delta_1}(K,\d))$ and $\Upsilon := (\nu \delta) \vee \ve$, provided that $\xi_0 \Upsilon \leq  2\sdr_{\delta_1}(K,\d)$.
\end{corollary}

Corollary~\ref{cor:stab} is proven in Appendix \ref{sec:proof-sdr-properties}. 
It ensures that the spherical distortion radius enjoys an interleaving property. 
That is the SDR of $(K,\d)$ at scale $\delta$ may be framed by the SDR of an approximation $(K',\d')$ at scales $\delta\pm \ve$. This interleaving property is a common thread with the $\mu$-reach (see, e.g.,~\cite[Theorem~3.4]{Chazal06}) and the $\lambda$-reach (\cite[Theorem~3]{Chazal05}), that is not enough to ensure consistent estimation. In fact, for the two aforementioned quantities, consistency may be proved with the additional assumption of $\mu \mapsto \rch_\mu(K)$ (resp. $\lambda \mapsto \lambda$-reach) are continuous at the targeted $\mu$ (resp. $\lambda$). 

As opposed to the $\mu$-reach the $\lambda$-reach, the SDR is also stable with respect to its the scale parameter $\delta$.
Next, we prove that $\delta \mapsto \sdr_\delta(K,\d)$ is continuous over a fixed range $(0,\Delta^*)$ under mild structural assumptions on $(K,\d)$. These assumptions  will be easily checked in the model $\manifolds{k}{\rch_{\min}}{\mathbf{L}}$, hence ensuring consistency of the subsequent reach estimator.

\begin{ass}\label{ass:spread}
We say that $K \subset \R^D$ is \emph{spreadable} if there exist $\Delta_0 > 0$, $\ve_0 > 0$, and $C_0 > 0$ such that for all $x,y \in K$ such that $\|x-y\| \leq \Delta_0$ and all $\ve \leq \ve_0$, there exists a point $a \in K$ such that either 
\bitem
\item $\|a-y\| \leq \ve$ and $\|x-a\| \geq \|x-y\| + C_0 \ve$, or
\item $\|a-x\| \leq \ve$ and $\|y-a\| \geq \|x-y\| + C_0 \ve$.
\eitem
\end{ass}
Assumption~\ref{ass:spread} requires that every point $y$ of $K$ may be locally pushed away from any (close enough) point $x \in K$. In particular, this means that $K$ is nowhere discrete.
In the manifold case, this pushing may be carried out using the exponential map (see Proposition~\ref{prp:assreach}).
 
\begin{ass}\label{ass:subeuc} 
We say that $(K,d)$ is \emph{sub-Euclidean} if there exist $C_1 > 0$ and $\Delta_1 > 0$ such that for all $x,y \in K$ such that $\|x-y\| \leq \Delta_1$, we have $\d(x,y) \leq C_1 \|x-y\|$.
\end{ass}
Assumption~\ref{ass:subeuc} requires that the distance locally compares with the ambient Euclidean distance. 
This essentially means that the identity map $(K,\d) \to (K,\norm{\cdot})$ is locally Lipschitz.
Such an assumption is automatically fulfilled whenever $K$ has positive reach and $\d = \d_K$ (see~\cite{Federer59}), with explicit constants in the manifold case (see Proposition~\ref{prp:assreach}) 
Whenever these two conditions are met, the spherical distortion radius of $(K,\d)$ can be proved to be locally Lipschitz in $\delta$.

\begin{theorem} \label{thm:lip} Assume that the metric space $(K,\d)$ fulfills Assumptions~\ref{ass:spread} and~\ref{ass:subeuc}. Then $\delta \mapsto \sdr_\delta(K,\d)$ is locally Lipschitz on $(0,\Delta^*)$ where
$$
\Delta^* :=  \min\{ \Delta_0, \Delta_1, \sup\{\delta \geq 0~|~\sdr_{\delta}(K,\d) < \infty\} \}.
$$
More precisely, for all $0 < \delta_0 < \delta_1 < \Delta^*$,  the map $\delta \mapsto \sdr_\delta(K,\d)$ is $L_0$-Lipschitz on $[\delta_0,\delta_1]$ with
\[
L_0 :=  \frac{192 \mathrm{r}_1^3}{C_0 \delta_0^3}\(C_1 + \pi\frac{\mathrm{r}_1}{\delta_0}\),
\]
where $\mathrm{r}_1 := \sdr_{\delta_1}(K,\d)$.
\end{theorem}

A proof of \thmref{lip} can be found in Appendix \ref{sec:proof-sdr-properties}. 
Not only does it ensure that the spherical distortion radius at scale $\delta$ is continuous with respect to $\delta$, that is enough to guarantee consistency, but it also allows to control its variation via an explicit local Lipschitz constant. 
Combined with  \corref{stab}, this allows to convert a bound between $(K,\d)$ and $(K',\d')$ in terms of Hausdorff distance and metric distortion into a bound on the SDR's at scale $\delta$.

\begin{figure}[h!]
\centering
\includegraphics[width = 0.45\textwidth]{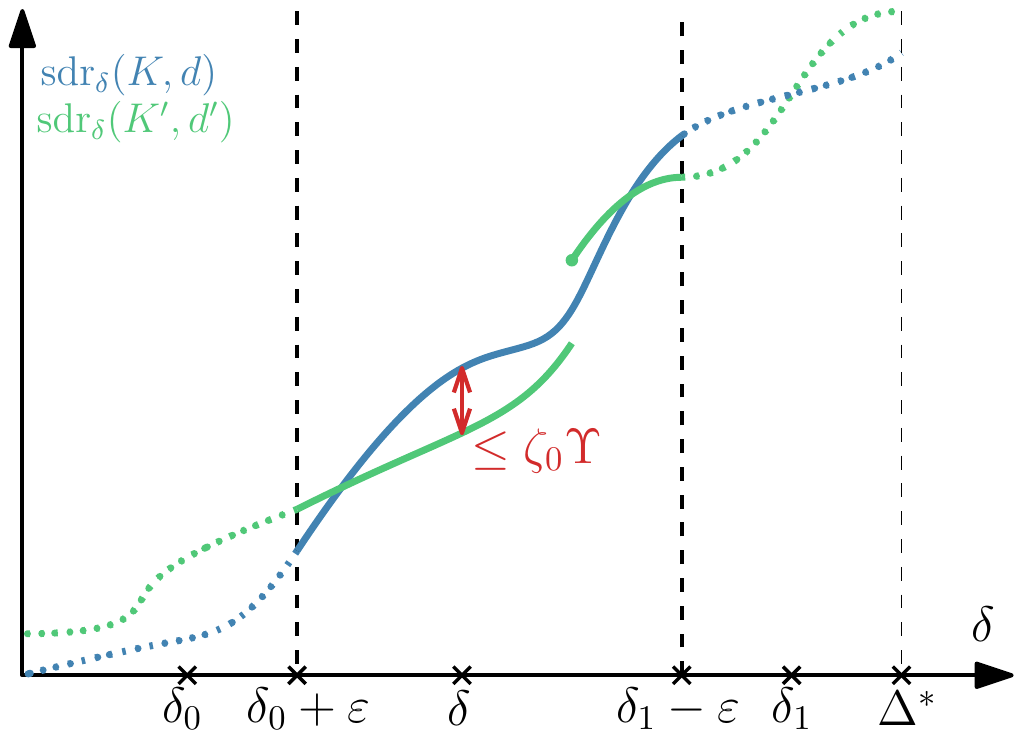}
\caption{Plot of $\delta \mapsto \sdr_\delta$ for $(K,d)$ and $(K',d')$ in the context of Theorem \ref{thm:stab}. On the interval $(\delta_0+\ve,\delta_1-\ve)$, the two functions do not differ of more than $\zeta_0 \Upsilon$. 
Even though $(K',d')$ might not be well-behaved, the regularity of $\delta \mapsto \sdr_\delta(K,d)$ (Theorem \ref{thm:lip}) is sufficient to insure stability.}
\label{fig:thm49}
\end{figure}

\begin{theorem} \label{thm:stab} Let $(K,\d)$ fulfill Assumptions~\ref{ass:spread} and~\ref{ass:subeuc}, and let $(K',\d')$ be such that $\dh(K,K') \leq \ve$ and $\D_{\delta_0}(\d,\d') \leq 1+\nu$ for some $\delta_0 < \Delta^*$. 
Then, for all $\delta_1 \in (\delta_0,\Delta^*)$ and $\delta \in (\delta_0+2\ve,\delta_1-2\ve)$, provided that $\xi_0 \Upsilon \leq 2\sdr_{\delta_1}(K,\d)$, we have
$$
\left|\sdr_{\delta}(K,\d) - \sdr_{\delta}(K',\d')\right| \leq \zeta_0 \Upsilon,
$$
with $\Upsilon = (\delta \nu) \vee \epsilon$ and $\zeta_0 = \xi_0 + 2L_0$, where $\xi_0$ is defined in \corref{stab}, $L_0$ is defined in \thmref{lip}.
\end{theorem}
We refer to Appendix \ref{sec:proof-sdr-properties} for a proof of this result and to Figure \ref{fig:thm49} for a diagram of the scales at play. Note that the constant $\zeta_0$ only depends on $\delta_0$ and features of $(K,\d)$, that the assumptions are required on $(K,\d)$ only, and that the constraint on $\ve$ depends only on $(K,\d)$ as well. 

The estimation of $K$ is a now well-understood in the manifold case (see~\cite{Aamari19b}). To obtain guarantees on the estimation of $\sdr_\delta(K,\d_K)$, it hence remains to investigate the estimation of $\d_K$. This is the aim of the following section.

\section{Optimal Metric Learning}
\label{sec:metric}
\subsection{Unsupervised Distance Metric Learning}

As explained in the introduction, various learning tasks lead to the problem of estimation the shortest-pat distance $\d_K$, via an estimator $\hat{d}$ on a sample of $K \subset \R^D$.
Though, there is no canonical choice of loss for measuring the proximity of $\wh \d$ to  $\d_K$. 
One could consider for instance  the empirical $\sup$-loss 
$$\ell_n(\wh \d | \d_K) := \sup_{x \neq y \in \bbX_n} \left|1- \frac{\wh \d(x,y) }{\d_K(x,y)}\right|,
$$
or the global $\sup$-loss
$$
\ell_\infty(\wh \d | \d_K) := \sup_{x \neq y \in K} \left|1- \frac{\wh \d(x,y) }{\d_K(x,y)}\right|.
$$

It might seem counter-intuitive to ask an estimator $\wh \d$ of $\d_K : K \times K \to \R_+$ to be defined on the whole set $K\times K$, while the this domain is unknown.
It actually is easy to extend any metric estimator to the whole space $\bbR^{D} \times \bbR^D$.
Indeed, given such a metric estimation procedure $\wh \d_n : \X_n \times \X_n \to \R_+$ that outputs a distance  $\wh \d_n[\bbX_n](x,y)$ between any pair of points of $\bbX_n$, we can define $\wt \d_n(x,y) := \wh \d_{n+2}[\bbX_n,x,y](x,y)$ for all $(x,y) \in \bbR^D \times \bbR^D$.
Informally this means that one can treat any given tuple of points $(x,y)$ as actual data points in the estimation process, and that we are only interested in the behavior of the later when $x$ and $y$ are in fact from $K$.

The losses $\ell_n$ and $\ell_\infty$ are naturally multiplicative, in particular because the usual notions of distortions are multiplicative by nature (see \secref{sdr}).
Indeed, the sup-loss $\ell_\infty(\wh \d | \d_K)$ being smaller than $\nu$ means that 
$$
\forall x,y \in K,~~~(1-\nu)\d_K(x,y) \leq \wh \d(x,y) \leq (1+\nu) \d_K(x,y),
$$
which is the usual way to quantify if the intrinsic metric is well-estimated. See for instance~\cite{tenenbaum2000global, Arias20}. 
When $\nu$ is small, it yields that $(K,\wh\d)$ is quasi-isometric to $(K,\d_K)$.

\begin{rem} We emphasize the fact that the global sup-loss $\ell_\infty$ and the mutual metric distortion $\D_\delta$ from Definition~\ref{def:metric_distortion} are different in essence. 
Indeed, while the mutual metric distortion $\D_\delta$ allows to compare different metrics on \emph{different} subsets of $\bbR^D$, the sup-loss $\ell_\infty$ compares two distances defined on the \emph{same} subset.

However, the global sup-loss and the mutual distortion metric may be related as follows. 
Consider $K$ endowed with either $\wh \d$ or $\d_K$. 
Denote by $\D_{0^+}(\d_K,\wh{\d}) := \lim_{\delta \rightarrow 0}\D_{\delta}(\d_K,\wh{\d})$.
Then, straightforward computation entails \begin{align*}
\ell_\infty(\wh{\d} | \d_K) +1 \leq \D_{0^+}(\d_K,\wh{\d}) \leq (1-\ell_\infty(\wh{\d}|\d_K))_+^{-1}.
\end{align*} 
Hence, the global sup-loss $\ell_\infty(\wh{\d}|\d_K)$ is somehow an additive counterpart to the mutual distortion $\D_{0^+}(\wh{\d},\d_K)$ in the case where $K=K'$. That is, when the support of the two metrics coincide in Definition~\ref{def:metric_distortion},  as already noticed in Remark \ref{rem:bilip}.
\end{rem}

When $K = M$ is a $\mathcal{C}^2$ submanifold of $\bbR^D$ of dimension $d$ with reach bounded below, methods using neighborhood graphs such as Isomap provably estimate $\d_M$ at rate $O(n^{-2/3d})$~\cite{arias2019unconstrained}.
As we will show in Theorem~\ref{thm:metric-minimax-ub}, this rate is far from being optimal. 
To date, the best minimax lower bound in this setting is due to~\cite{Arias20}, who obtain a rate of order $\Omega(n^{-2/d})$ in the particular case of a deterministic design on $\cC^2$ submanifolds.
Actually, we can extend the result of~\cite{Arias20} to our random design setting, and to general $\cC^k$ submanifolds with $k \geq 2$.

\begin{theorem}\label{thm:metriclb}
Assume that $f_{\min} \leq c_{d,k}/\rch_{\min}^d$ and $f_{\max} \geq C_{d,k}/\rch_{\min}^d$, and $L_j \geq C_{d,k}/\rch_{\min}^{j-1}$ for all $j \in \{2,\ldots,k\}$.
Then for $n$ large enough,
\begin{align*} 
%\label{metriclb}
 \inf_{\wh \d} 
 \sup_{P \in \distributions{k}{\rch_{\min}}{\mathbf{L}}{f_{\min}}{f_{\max}} 
} \bbE_{P^{\otimes n}}[\ell_\infty(\wh \d | \d_M)] 
\geq 
\tilde{c}_{d,k,\rch_{\min}}  \left(\frac{1}{n}\right)^{k/d},
\end{align*}
where the infimum is taken over all measurable estimator $\wh \d$ of $\d_M$ based on $n$ samples.
\end{theorem} 

This theorem is proved in Appendix \ref{sec:proof-lower-bound-metric}. As we shall prove shortly in Section~\ref{sec:optimal-metric-estimation}, this lower-bound can be provided with a matching upper-bound up to $\log n$ factors (\thmref{metricub}), and is thus optimal.

\subsection{An optimal Approach of Metric Estimation}
\label{sec:optimal-metric-estimation}

The existing unsupervised methods for metric learning are known to either have no theoretical guarantees, or to have a sub-optimal rate for estimating the intrinsic metric.
As stated before, Isomap reaches a rate of $n^{-2/3d}$, which is very far from the theoretical lower-bound $n^{-k/d}$ shown in \thmref{metriclb}. Other methods, such as taking the shortest path distance over a Delaunay triangulation~\cite{Arias20}, are shown to attain a precision of $n^{-2/d}$ which is optimal for $\cC^2$-model but not for $k \geq 3$. We propose here a fairly general approach that can output a family of minimax-optimal metric estimators. It relies on the following bound.

\begin{proposition} \label{prp:metric}
Let $K \subset \bbR^D$ be a set of positive reach $\rch(K) > 0$, and $K' \subset \bbR^D$ be any set such that $\dh(K',K) < \ve \leq \rch(K)/2$. Then,
$$
\ell_\infty(\d_{(K')^\ve} | \d_{K}) \leq \frac{2\ve}{\rch(K)}
,
$$
where 
we recall that
$(K')^\ve = \{ u \in \R^D \mid \d(u,K') \leq \ve\}$, so that $K \subset (K')^\ve$.

\end{proposition}
\prpref{metric} is proved in Appendix \ref{sec:proof-plugin-metric}.
It asserts that estimating geodesic distances of sets of positive reach is never harder than estimating the sets themselves in Hausdorff distance.
Beyond the framework of closed manifold developed here, note that for the convex case $\rch(K) = \infty$, $\d_K$ coincides with the Euclidean metric, so that estimating $\d_K$ becomes trivial.

A significant consequence of \prpref{metric} is that we can derive a consistent estimator of the intrinsic distance from any consistent estimator of the support, and with the same rate of convergence. 
In what follows, we write
\begin{align}
\label{eq:dmax}
\d_{\max} := \frac{5^d}{\omega_d f_{\min} \rch_{\min}^{d-1}},
\end{align}
where $\omega_d$ is the volume of the $d$-dimensional unit ball. In \lemref{boundgeo}, the length $\d_{\max}$ is proved to be an upper bound on the geodesic diameter of the supports of any distribution in the model $\distributions{k}{\rch_{\min}}{\mathbf{L}}{f_{\min}}{f_{\max}}$.

\begin{theorem} \label{thm:metricub} Let $k \geq 2$ and let $\wh M$ be an estimator satisfying
$$
\sup_{P \in \distributions{k}{\rch_{\min}}{\mathbf{L}}{f_{\min}}{f_{\max}}} P^{\otimes n}(\dh(\wh M, M) \geq \ve_n) \leq \eta_n, 
$$
for some positive sequences $\ve_n$ and $\eta_n$ converging to $0$.
Then the metric estimator
$$\wh \d(x,y) := \d_{\max} \wedge  \d_{(\wh M_{x,y})^{\ve_n}}(x,y)~~~\text{with}~~~\wh M_{x,y} := \wh M \cup \{x,y\}, 
$$ 
which is defined for all $x,y \in \bbR^D$, satisfies
$$
\sup_{P \in \distributions{k}{\rch_{\min}}{\mathbf{L}}{f_{\min}}{f_{\max}}} \bbE_{P^{\otimes n}}[\ell_\infty(\wh \d | \d_M)] \leq \frac{2}{\rch_{\min}}\ve_n+\left(1+\frac{\d_{\max}}{\ve_n}\right)\eta_n.
$$
\end{theorem}

\thmref{metricub} is proved in Appendix \ref{sec:proof-plugin-metric}. 
A particular advantage of this result is that it does not require the estimator $\wh M$ to have any geometric structure, nor to be regular in any sense.
This contrasts sharply with \cite{Arias20}, which extensively uses the structural properties of the intermediate estimator $\wh M$.
\thmref{metricub} is much more versatile, since here, $\wh M$ could just as easily be anything as a point cloud, a metric graph, a triangulation, or a union of polynomial patches. 
For instance, taking $\wh M = \{X_1,\dots,X_n\}$ to be the observed data, we can take $\ve_n = C(\log n/n)^{1/d}$ for $C$ large enough yields $\eta_n \leq \ve_n^2$ so that
$$
\sup_{P \in \distributions{2}{\rch_{\min}}{\mathbf{L}}{f_{\min}}{f_{\max}}} \bbE_{P^{\otimes n}} [\ell_\infty(\wh \d | \d_M)] \leq C_{\rch_{\min}, d, f_{\min}} \left(\frac{\log n}{n}\right)^{1/d},
$$ 
which is faster than the known rate of order $O(n^{-2/3d})$ for Isomap (see for instance \cite[Eq (1.2)]{Arias20}). 
Now, taking $\wh M$ to be a minimax optimal estimator of $M$ for the Hausdorff loss --- as that of~\cite{Aamari19}, for instance --- and $\ve_n = C (\log n/n)^{k/d}$  for some large constant $C > 0$ yields $\eta_n \leq \ve_n^2$ (see \lemref{hausdorff_and_covering}), and a metric estimator $\wh \d$ that achieves the following rate.
\begin{theorem} \label{thm:metric-minimax-ub}
Let $\wh \d$ be the estimator described in \thmref{metricub} built on top of $\wh M$ described in \lemref{hausdorff_and_covering}. Then for $n$ large enough,
$$
\sup_{P \in \distributions{k}{\rch_{\min}}{\mathbf{L}}{f_{\min}}{f_{\max}}} \bbE_{P^{\otimes n}} [\ell_\infty(\wh \d | \d_M)] \leq C_{\rch_{\min}, d, f_{\max}, f_{\min}, \bL,k} \left(\frac{\log n}{n}\right)^{k/d}.
$$ 
\end{theorem}
In virtue of \thmref{metriclb}, this rate is minimax optimal up to $\log n$ factors.

\section{Optimal Reach Estimation}
\label{sec:optimal_reach_estimation}
\subsection{Optimal Spherical Distortion Radius Estimation} 
Interesting as it is in its own right, we now investigate the estimation rates of the spherical distortion radius at scale $\delta > 0$.
To obtain a minimax lower bound, we simply note that $\sdr_\delta(M,\d_M)$ coincides with $\rch(M)$ whenever $\rch(M) = \wfs(M)$ (Proposition~\ref{prp:interpolate}). 
Hence, any lower bound for the estimation of $\rch(M)$ on a model over which $\rch(M) = \wfs(M)$ yields a lower bound for the estimation of $\sdr_\delta(M,\d_M)$.
In application of \thmref{lwr_bounds_clementb} with $\alpha \geq 0$, this immediately gives the following lower bound.
\begin{theorem} 
\label{thm:sdrlb} 
Assume that $f_{\min} \leq c_{d,k}/\rch_{\min}^d$ and $f_{\max} \geq C_{d,k}/\rch_{\min}^d$, and $L_j \geq C_{d,k}/\rch_{\min}^{j-1}$ for all $j \in \{2,\ldots,k\}$.
Then for $n$ large enough,
for all $\delta \in (0,\rch_{\min})$,
$$
\inf_{\wh\sdr_\delta} \sup_{P \in \distributions{k}{\rch_{\min}}{\mathbf{L}}{f_{\min}}{f_{\max}}} \bbE_{P^{\otimes}}[|\wh\sdr_\delta - \sdr_\delta(M,\d_M)|] \geq \tilde{c}_{\rch_{\min},d,k} n^{-k/d}.
$$
where the infimum is taken over all measurable estimators $\wh\sdr_\delta$ of $\sdr_\delta(M,\d_M)$ based on $n$ samples.
\end{theorem}
It turns out that this bound is optimal. 
To exhibit an estimator that achieves this rate, we take advantage of the Hausdorff and metric stability of the spherical distortion radius shown in \thmref{stab}. 
In order to apply it, we first need to check that Assumptions~\ref{ass:spread} and~\ref{ass:subeuc} are fulfilled for every manifolds in our models $\cC^k_{\rch_{\min},\bL}$.

\begin{proposition} \label{prp:assreach} Let $M \subset \bbR^D$ be a submanifold with bounded reach $\rch(M) > 0$. Then $M$ satisfies Assumptions \ref{ass:spread} and \ref{ass:subeuc} with parameters 
$$\ve_0 = \rch(M)/4, ~~~~\Delta_0 = \rch(M) ,~~~~C_0 = 3/16,~~~~\Delta_1 = \rch(M)/2 ~~~~\text{and} ~~~~C_1 = 2.
$$
\end{proposition}

\prpref{assreach} is proven in Appendix \ref{sec:proofoptireach}. In the vein of \thmref{metricub}, and using the stability of the spherical distortion radius with respect to the pair $(K,\d)$, we can now build an estimator of $\sdr_\delta(M,\d_M)$ in a plug-in fashion over $\cC^k$ submanifolds.
Recall that when $M$ is in $\cC^k_{\rch_{\min},\bL}$, and $\delta \in (0,\sqrt{2(D+1)/D}\wfs(M))$, then according to Propositions \ref{prp:interpolate} and \ref{prop:wfs_properties}, and to \lemref{boundgeo},
$$
0 < \rch_{\min} \leq \rch(M) \leq \sdr_\delta(M,\d_M) \leq \wfs(M) \leq \sqrt{\frac{D}{2(D+1)}} \diam(M) \leq \s_{\max} < \infty
,$$
where $\s_{\max} :=  \sqrt{D/(2(D+1))}  \d_{\max}$, with $\d_{\max}$ being the constant introduced in \eqref{eq:dmax}. 

\begin{theorem} \label{thm:sdrub} Given $k \geq 2$, let $\wh M$ be an estimator satisfying
$$
\sup_{P \in \distributions{k}{\rch_{\min}}{\mathbf{L}}{f_{\min}}{f_{\max}}} P^{\otimes n}(\dh(M,\wh M) \geq \ve_n) \leq \eta_n 
$$
for some positive sequences $\ve_n,\eta_n$ converging to $0$. Then, for any $\delta \in (0,\rch_{\min})$, the estimator $\wh\sdr_\delta := \sdr_\delta(\wh M,\wh\d) \wedge \s_{\max}$, where $\wh \d$ is defined in \thmref{metricub}, satisfies
$$
\sup_{P \in \distributions{k}{\rch_{\min}}{\mathbf{L}}{f_{\min}}{f_{\max}}} \bbE_{P^{\otimes n}} |\wh\sdr_\delta -\sdr_\delta(M,\d_M)| \leq C \left(\frac{\s_{\max}^4}{\delta^4} \ve_n +  \s_{\max}\eta_n\right).
$$
\end{theorem}

We refer to  Appendix \ref{sec:proofoptireach} for a proof of this result. 

\begin{rem}
In place of $\wh \d = \d_{\wh M^{\ve_n}}$, one could actually plug any estimator $\wh \d$ of the metric into\thmref{sdrub}. 
In light of the stability result of \thmref{stab}, as long as $\wh \d$ satisfies
\begin{align*}
\sup_{P \in \cP^k} P^{\otimes n}\(\D_\delta(\wh \d,\d_M) \geq 1+\frac{\ve_n}{\delta}\) \leq \eta_n, 
\end{align*}
the conclusion of \thmref{sdrub} would still hold. 
This comes in handy, especially if one wants to input a computationally efficient distance estimator, such as shortest-path distance on a neigbhorhood graph~\cite{tenenbaum2000global} or on Delaunay triangulations~\cite{Arias20}.
\end{rem}

Again, taking $\wh M$ to be a minimax optimal estimator for the Hausdorff loss~\cite{Aamari19} outputs an estimator $\wh\sdr_\delta$ of the spherical distortion radius satisfying \begin{theorem} 
For all $\delta \in (0,\rch_{\min})$, with the construction of $\wh \sdr_\delta$ above, we have that for $n$ large enough,
$$
\sup_{P \in  \distributions{k}{\rch_{\min}}{\mathbf{L}}{f_{\min}}{f_{\max}}} \bbE_{P^{\otimes n}} |\wh\sdr_\delta -\sdr_\delta(M,\d_M)| \leq C_{\rch_{\min}, d, f_{\max}, f_{\min}, \bL,k} \frac{1}{\delta^{4}} \left(\frac{\log n}{n}\right)^{k/d},
$$
and this rate is optimal in regard of \thmref{sdrlb}.  
\end{theorem}
Note the presence of the factor $1/\delta^{4}$ in the bound, which makes the rate diverge as $\delta \to 0$.
This blowup is to be expected for the following reason. As $\delta$ goes to $0$, the spherical distortion radius goes to the reach $\rch(M)$ (Proposition~\ref{prp:interpolate}). 
Since the estimation of $\rch(M)$ cannot be faster than $n^{-(k-2)/d}$ (\thmref{lwr_bounds_clementb}), the estimation rate of $\sdr_\delta(M,\d_M)$ must deteriorate in some way as $\delta\to 0$.

\subsection{Optimal Reach Estimation}

In light of Proposition~\ref{prop:reach_intermediate_local} and~\eqref{eq:double-plugin-precision}, it only remains to combine the maximal curvature estimator and the spherical distortion radius estimator to obtain an estimator of the reach. 
Naely, we let $\wh M$ be the minimax-Hausdorff estimator of \lemref{hausdorff_and_covering}. According to the very same \lemref{hausdorff_and_covering}, there exists $c_{\rch_{\min}, d, f_{\max}, f_{\min}, \bL,k} > 0$ such that denoting by 
\beq \label{eq:tune_ven}
\ve_n := c_{\rch_{\min}, d, f_{\max}, f_{\min}, \bL,k} \(\frac{\log n}{n}\)^{k/d},
\eeq
there holds
\beq \label{eq:proba_ven}
\sup_{P \in \distributions{k}{\rch_{\min}}{\mathbf{L}}{f_{\min}}{f_{\max}}}P^{\otimes n}(\dh(\wh M,M) \geq \ve_n) \leq \ve_n^2.
\eeq
We also let $\wh \d$ be the estimator of the intrinsic distance of \thmref{metricub} from $\wh M$ and $\ve_n$. 
We let $\wh\sdr_\delta := \sdr_\delta(\wh M,\wh \d) \wedge \s_{\max}$ for some $\delta \in (0,\rch_{\min})$ as in Theorem~\ref{thm:sdrub}. 
Finally, we write
$$
\wh \rch := \wh R_\ell \wedge \wh\sdr_\delta.
$$
The following \thmref{reachub} is a straightforward consequence of Theorems \ref{thm:cv_rates_curvature_max} and \ref{thm:sdrub}, inserted in the plugin strategy of Proposition~\ref{prop:reach_intermediate_local} and~\eqref{eq:double-plugin-precision}.

\begin{theorem} \label{thm:reachub} 
The estimator $\wh\rch$ described above with $\delta = \rch_{\min}/2$ satisfies
\begin{align*} 
\sup_{P \in \distributions{k}{\rch_{\min}}{\mathbf{L}}{f_{\min}}{f_{\max}}} \bbE_{P^{\otimes n}} | \wh \rch - \rch(M) | \leq C_{\rch_{\min}, d, f_{\max}, f_{\min}, \bL,k}  \left(\frac{\log n}{n}\right)^{(k-2)/d},
\end{align*} 
and, for all $\alpha > 0$,
\begin{align*} 
\sup_{P \in \distributions{k}{\rch_{\min}}{\mathbf{L}, \alpha}{f_{\min}}{f_{\max}}} \bbE_{P^{\otimes n}} | \wh \rch - \rch(M) | \leq C_{\rch_{\min}, d, f_{\max}, f_{\min}, \bL,k,\alpha} \left(\frac{\log n}{n}\right)^{k/d}. 
\end{align*} 
\end{theorem}

As a conclusion, Theorems~\ref{thm:lwr_bounds_clementb} and~\ref{thm:reachub} assert that $\wh \rch$ is minimax optimal, and that its rate of convergence adapts to whether $\rch(M)$ is attained by curvature (yielding the slower rate $O(n^{-(k-2)/d})$) or by a bottleneck (yielding the faster rate rate $O(n^{-k/d})$).

The computation of $\wh\rch$ depends explicitly on the parameters of the models at two levels. 
First, in tuning the value of $\ve_n$ as in \eqref{eq:tune_ven}. 
Second, in choosing $\delta \in (0,\rch_{\min})$. These two dependencies may be circumvented by picking
$$\wt \ve_n = {\log n}\left(\frac{\log n}{n}\right)^{k/d},$$ 
and $\delta_n = 1/\log n$. Then, for $n$ large enough, both \eqref{eq:proba_ven} and $\delta_n \in (0,\rch_{\min})$ will be fulfilled.
The price to pay for this way-around to calibration of constants limits to multiplicative $\log n$ factors in the upper-bound of \thmref{reachub}.

\section{Conclusion and Further Prospects}
\label{sec:conclusion}
We developed a general strategy for estimating the reach of a manifold $M$. 
It relies on two independent plugins, accountable for the estimation of the minimal curvature radius $R_\ell(M)$ and any another set-defined feature size $\theta(M)$ that lie between the reach and the weak feature size.
We then introduced and studied the spherical distortion radius, the estimation of which reduces to geodesic distance estimation, itself reducing to set estimation in Hausdorff distance.
All the derived results are minimax optimal, as testified by associated matching lower bounds up to $\log n$ factors.

Geometrically, one should note that this overall method relies heavily on the local/global dichotomy of the reach for closed submanifolds~\cite{Aamari19b}. Hence, it still remains unclear how to extend it to manifolds with boundary, even though their curvature and spherical distortion radius are likely to be estimated in a similar way~\cite{Aamari21}.

On the statistical side, a major extension of the results would consist in allowing for additive noise. 
Recent works obtained Hausdorff estimation rates for the support~\cite{Fefferman19,Aizenbud21,Puchkin22} in such a noisy setting, so that the estimation of the spherical distortion radius inherits the same rates straightforwardly.
In the same spirit as the iterated local polynomial fitting of~\cite{Aizenbud21}, we expect that the same method could likewise lead to maximal curvature estimation. 

Finally, since the main goal of this work was of minimax nature, we did not focus on the algorithmic properties of our estimators.
As they stand, $\wh R_\ell$ and $\wh \sdr$ both require to compute a supremum over the union of continuous patches $\wh M$, which is computationally prohibitive.
Actually, one can easily show that taking the same supremum over a discretization of $\wh M$ at scale $O\bigl(n^{-\beta/d}\bigr)$ -- i.e. $O(n^{\beta})$ points in total -- yields estimation rates of order $O\bigl(n^{-(\beta\wedge(k-2))/d}\bigr)$ for $R_\ell(M)$, and $O\bigl(n^{-(\beta\wedge k)/d}\bigr)$ for $\sdr_\delta(M,\d_M)$.
This suggests a possible estimation-computation tradeoff which one could take advantage of.
Yet, this is not a fully satisfactory solution, as $\sdr_\delta$ still requires to compute costly geodesic distances on a high-dimensional set.
More globally, the quest for computationally efficient -- yet optimal -- geometric estimators in high dimensions is still in its infancy.

\section*{Acknowledgments}
The authors would like to thank heartily \emph{Chez Adel} for its unconditional warmth and creative atmosphere, and Vincent Divol for helpful discussions.

\bibliographystyle{chicago}
\bibliography{biblio,biblio_cb}

\newpage

\appendix
\section{Proofs of \secref{reach}} \label{sec:proof3}

\subsection{Comparing Reaches, Weak Feature Size and Diameter} \label{sec:proof_of_prp_wfs_properties}
This Section is devoted to the Proof of Proposition~\ref{prop:wfs_properties}, which goes as follows.
\begin{proof}[Proof of Proposition~\ref{prop:wfs_properties}] 
For (i), recall that no closed compact submanifold can be contractible~\cite[Theorem~3.26]{hatcher2002algebraic}. Furthermore,~\cite[Theorem~4.8]{Federer59} and~\cite[Lemma~2.1]{Chazal06} combined together yield that $K^r$ is isotopic to $K$ for all $r < \wfs(K)$. 
On the other hand whenever $r > \rad(K)$ where $\rad(K)$ is the radius of the smallest ball enclosing $K$, $K^r$ is star-shaped with respect to any point of the non-empty intersection $\cap_{x \in M} \ball(x,r)$. 
We conclude that $\wfs(K) \leq \rad(K)$, Since $\rad(K) < \infty$ because $K$ is compact, we obtain $\wfs(K) < \infty$.

For (ii), the first two inequalities come from the definition of $\rch_\mu(K)$ (see~\eqref{eq:mu-reach}). 
The rightmost comes Jung's Theorem~\cite[Theorem~2.10.41]{Federer69}, which asserts that
$
 \rad(K) \leq \sqrt{\frac{D}{2(D+1)}}\diam(K)
$, 
and the fact that $\wfs(K) \leq \rad(K)$ whenever $\wfs(K)$ is finite (same argument as for (i)).
\end{proof}

\subsection{Minimax Lower Bound for $\mu$-Reach Estimation}
\label{sec:proof-of-lower-bound-mu-reach}
This Section is devoted to the proof of Theorem~\ref{thm:mu-reach-inconsistency}. It builds upon the possible discontinuities of the map $M \mapsto \Med_\mu(M)$ in Hausdorff distance.
The exhibition of such a discontinuity can be done in dimension $d=1$ and $D=2$, and can then be generalized to arbitrary $1 \leq d < D$ by using symmetry and rotation arguments.

\begin{figure}[h!]
\centering
\includegraphics[width=0.35\textwidth]{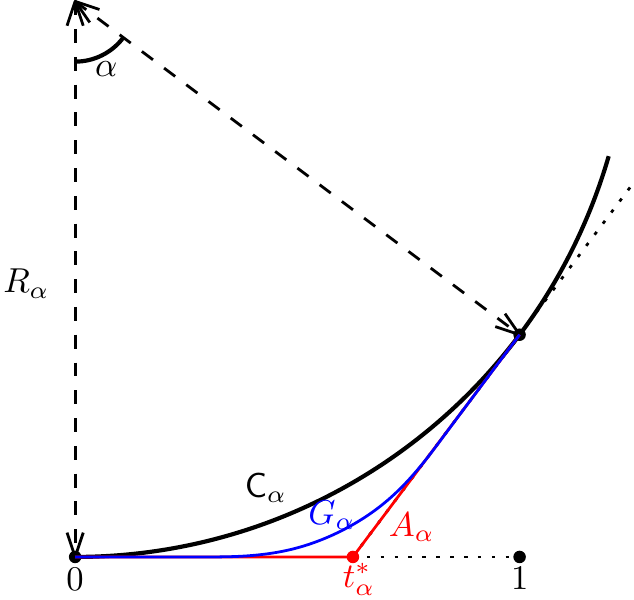}
\caption{Construction for Lemma~\ref{lem:turn-widget}: curves associated to $\circlett_\alpha$, $A_\alpha$, and $G_\alpha$.}
\label{fig:galpha}
\end{figure}

The building block of the construction is the following arc of curve.
For all $\alpha\in (0,\pi/4]$, write $R_\alpha := 1/\sin(\alpha)$. Let also $\circlett_\alpha : [0,1] \to \bbR_+$ be defined as $\circlett_\alpha(t) := R_\alpha - \sqrt{R_\alpha^2 - t^2}$, which graph is an arc of circle of radius $R_\alpha$ and aperture $\alpha$ (see Figure~\ref{fig:galpha}).
To be able to glue up smoothly $\alpha$-turns like $\circlett_\alpha$ with straight lines, we smooth it as follows.

\begin{lemma}
\label{lem:turn-widget}
There exists $G_\alpha : [0,1] \to \bbR_+$ infinitely differentiable such that:
\begin{enumerate}
\item $G_\alpha^{(\ell)}(0) = 0$ for all $\ell \geq 0$;
\item $G_\alpha(1) = \circlett_\alpha(1)$, $G_\alpha'(1) = \circlett_\alpha'(1)$ and $G_\alpha^{(\ell)}(1) = 0$ for all $\ell \geq 2$;
\item $\|G_\alpha^{(\ell)}\|_\infty \leq C_{\ell}/R_\alpha$ for all $\ell \geq 1$;
\item $G_\alpha(t) < \circlett_\alpha(t)$ for all $t \in (0,1)$; 
\item $G_\alpha$ is convex.
\end{enumerate}
\end{lemma}
See \figref{galpha} for a diagram of such a $G_\alpha$.
Let us first comment on the requirements on $G_\alpha$. Items 1 and 2 say that $G_\alpha$ is a $\cC^k$ interpolation between the two tangent lines of two points of $\circlett_\alpha$ who are $\alpha$-apart in term of polar coordinate. 
Item 3 says that the graph of $G_\alpha$, once rescaled by $1/R_\alpha$, will be bounded in $\cC^k$-norm for all $k$. 
Items 4 and 5 ensure well-behavior of the medial axes of our future construct (see Figure~\ref{fig:malpha}). 

\begin{proof}[Proof of Lemma~\ref{lem:turn-widget}]
The following construction applies to general convex functions, although we restrict it to $G_\alpha$ for simplicity.
Consider the piecewise linear map $A_\alpha$ given by  the tangent lines of $\circlett_\alpha$ at $t=0$ and $t=1$. That is, define $A_\alpha(t)$ for all $t \in \bbR$ by
\begin{align*}
A_\alpha(t)
:=&
\max\{ \circlett_\alpha(0) + \circlett_\alpha'(0)t , \circlett_\alpha(1) + (t-1)\circlett_\alpha'(1)\}
\\
=&
\max\{ 0 , \circlett_\alpha(1) + (t-1)\circlett_\alpha'(1) \}
.
\end{align*}
As $\circlett_\alpha$ is strictly convex, $A_\alpha < \circlett_\alpha$ on $\bbR\setminus \{0,1\}$.
We also denote by $t_\alpha^\ast$ the (unique) point of non-differentiability of $A_\alpha$, that is
\begin{align*}
t_\alpha^\ast
:=
1 - \frac{\circlett_\alpha(1)}{\circlett_\alpha'(1)}
=
R_\alpha \tan(\alpha/2)
.
\end{align*}
Note by now that for all $\alpha \in (0,\pi/4)$, $1/2 \leq t_\alpha^* \leq 2-\sqrt{2} \leq 6/10$.
Given $h>0$ to be chosen later, write $K_h(t) := h^{-1}K(t/h)$, where $K(t) := c_0 \exp(-1/(1-t^2)) \ind_{|t|< 1}$ is a non-negative $\mathcal{C}^\infty$ kernel, and $c_0$ is chosen so that $\int_{\bbR} K = 1$.
Finally, consider the convolution
\begin{align*}
G_\alpha(t) := \int_{\bbR} K_h(x) A_\alpha(t-x) \diff x
.
\end{align*}
By smoothness of $K_h$ and non-negativity of both $K_h$ and $A_\alpha$, $G_\alpha = K_h \ast A_\alpha$ is infinitely differentiable and non-negative.
Also, since $A_\alpha$ is convex and $K_h$ non-negative, $G_\alpha$ is convex (Item 5). Furthermore, one easily checks that outside the interval $[t_\alpha^\ast-h,t_\alpha^\ast+h]$, $G_\alpha$ coincides with $A_\alpha$. 
Hence, if $h \leq 1/4$, we have $[t_\alpha^\ast-h,t_\alpha^\ast+h] \subset [1/2-1/4,6/10+1/4]= [1/4,17/20]$, so that Items 1 and 2 holds directly.

To check that $G_\alpha < \circlett_\alpha$ on $(0,1)$, fix $t \in (0,1)$.
If $t \notin [t_\alpha^\ast-h,t_\alpha^\ast+h]$, $G_\alpha(t) = A_\alpha(t) < \circlett_\alpha(t)$ by construction.
If $t \in [t_\alpha^\ast-h,t_\alpha^\ast+h]$, we have $G_\alpha(t) \leq G_\alpha(t_\alpha^\ast + h) = h \circlett_\alpha'(1)$.
But on the other hand, $\circlett_\alpha(t) \geq \circlett_\alpha(t_\ast^\ast-h) > \circlett_\alpha(1/4)$. Hence, we do have $G_\alpha(t) < \circlett_\alpha(t)$ as soon as $h \leq 1/100$, since $\circlett_\alpha(1/4)/\circlett_\alpha'(1) > 1/100$ for all $\alpha \in (0,\pi/4)$. This yields Item 4.

Finally, letting $h = h_0 = 1/100$, we obtain for all $\ell \geq 1$ and $t \in [0,1]$,
\begin{align*}
|G_\alpha^{(\ell)}(t)|
&=
\bigl|
K_h^{(\ell)}
\ast
\circlett_\alpha(t)
\bigr|
\leq
\bigl\Vert
K_h^{(\ell)}
\bigr\Vert
_\infty
\bigl\Vert
\circlett_\alpha
\bigr\Vert_\infty
\leq
C_{\ell} \circlett_\alpha(1)
\leq
C_{\ell}/R_{\alpha},
\end{align*}
which yields Item 3 and concludes the proof.
\end{proof}

Given $R>0$, we now let $G_{\alpha,R}$ be the curve obtained by dilating homogeneously the graph of $G_\alpha$ by a scale factor $R/R_\alpha$. 
We extend the construction of these smooth $\alpha$-turns for $\alpha \in (\pi/4,\pi]$: for this, we glue two $G_{\alpha/2,R}$ or four $G_{\alpha/4,R}$ to define $G_{\alpha,R}$.

\begin{proposition}
\label{prop:hypotheses}
Assume that for all $j \in \{2,\ldots,k\}$, $L_j \geq C_{d,k}/\rch_{\min}^{j-1}$ for $C_{d,k}>0$ large enough.
Then for all $\mu \in [0,1)$ and $\varepsilon >0$ small enough, there exist $M,M' \in \manifolds{k}{\rch_{\min}}{\mathbf{L}}$ such that:
\begin{itemize}
\item
$
|\rch_\mu(M) - \rch_\mu(M')| \geq c_{d,k}\rch_{\min}
$
;
\item
$c'_{d,k} \rch_{\min}^d 
\leq 
\vol_d(M) \wedge \vol_d(M')
\leq
\vol_d(M) \vee \vol_d(M')
\leq C''_{d,k} \rch_{\min}^d$;
\item
$\vol_d(M \triangle M') \leq C'''_{d,k} \rch_{\min}^{d} \varepsilon$
.
\end{itemize}
\end{proposition}

\begin{proof}[Proof of Proposition~\ref{prop:hypotheses}]
For small enough (and arbitrarily small) $\ve >0$, we let $\alpha \in [0,\pi]$ be such that $\sin\bigl((\alpha+\ve)/2\bigr)^2 = 1-\mu^2$. Such an $\alpha$ always exists since $\mu^2 < 1$.
Given $\Delta,R_0,R_1> 0$ to be chosen later, we glue smooth turns from Lemma~\ref{lem:turn-widget} with straight lines to create a $\mathcal{C}^k$ closed curve in $\R^2$, as shown in Figure~\ref{fig:malpha}. 
Then, we obtain a $\cC^k$ closed $d$-dimensional submanifold $M_\alpha$ of $\bbR^{d+1}$, with a symmetry of revolution with respect to the horizontal axis of \figref{malpha}. 
\begin{figure}[h!]
\centering
\includegraphics[width=0.7\textwidth]{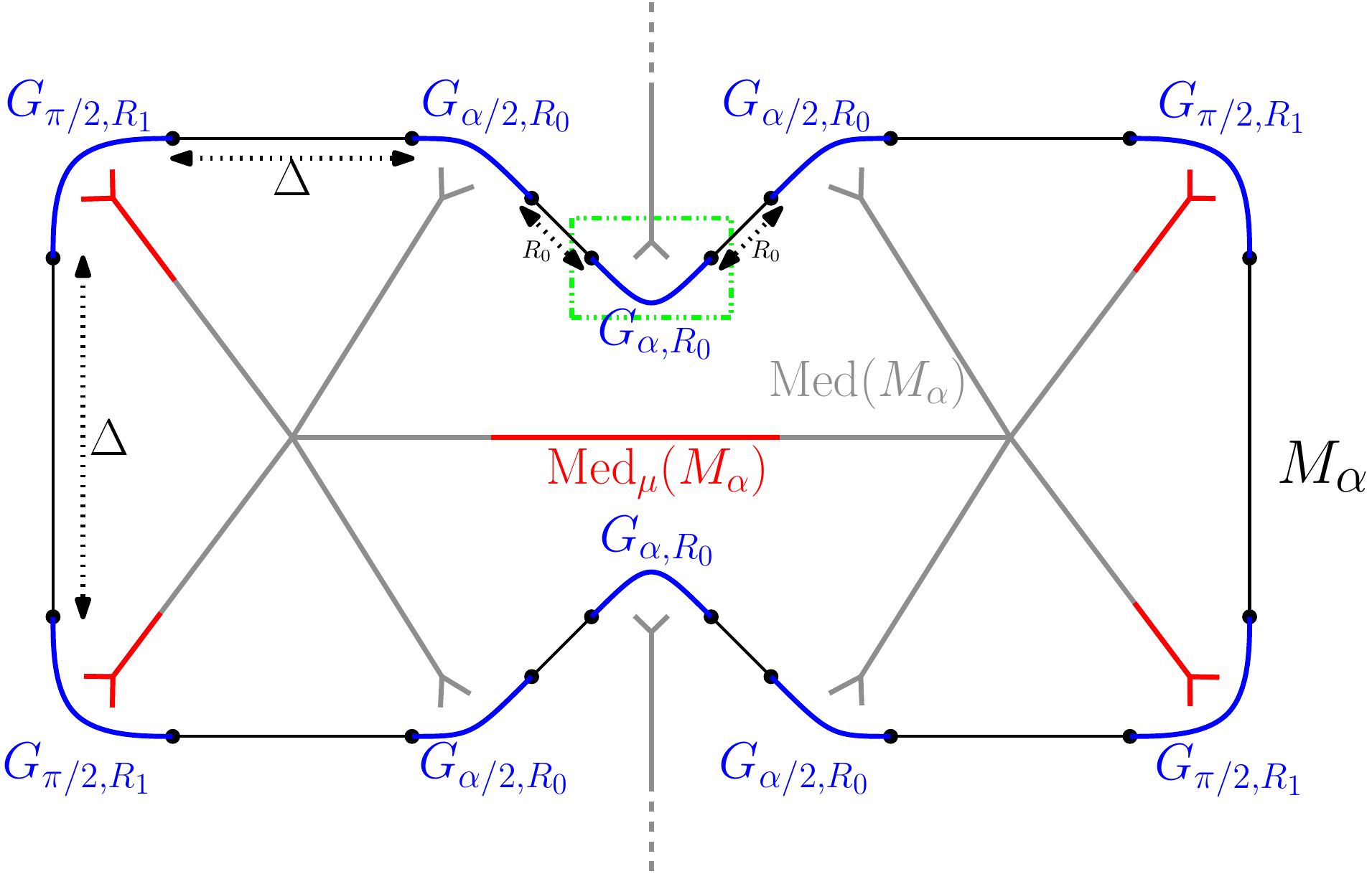}
\caption{Construction of $M_\alpha$ in the proof of Proposition~\ref{prop:hypotheses}.}
\label{fig:malpha}
\end{figure}

By construction, if $\Delta \geq 8R_0$, then $M_\alpha$ has local parametrizations on top of its tangent spaces (see Definition~\ref{def:geometric-model}) with $L_j \leq C_{d,k}/(\Delta \wedge R_0 \wedge R_1)^{j-1}$ for all $j \geq 2$, and has volume $\vol_d(M_\alpha) \leq C_{d,k} (\Delta \vee R_0 \vee R_1)^d$ and $\vol_d(M_\alpha) \geq c_{d,k}  (\Delta \wedge R_0 \wedge R_1)^d$.

We now examine the structure of the medial axis and the reach of $M_\alpha$.
If $u \in \Med(M_\alpha)$ is a point on the medial axis, rotational symmetry yields that two of its projections points must lie either:
\begin{itemize}
\item
In a plane containing its horizontal axis of symmetry (i.e. Figure~\ref{fig:malpha}).
As a result, its distance to $M_\alpha$ cannot be smaller than the smallest reach of each of its parts $G_{\pi/2,R_1}, G_{\alpha/2,R_0}$ and $G_{\alpha,R_0}$, so that $\d(u,M_\alpha) \geq c_{d,k} R_0 \wedge R_1$.
\item
In a $d$-plane orthogonal to the horizontal axis. By rotational invariance, this forces $u$ to be on this axis of symmetry. As a result, $\d(u,M_\alpha) \geq \Delta/2 - 3R_0 \geq c_{d,k} \Delta$ since $\Delta \geq 8R_0$.

\end{itemize}
In all, we get 
$\rch(M_\alpha) \geq c_{d,k} (\Delta \wedge R_0 \wedge R_1).$

We now examine the $\mu$-reach of $M_\alpha$. By definition, if $u \in \Med_\mu(M_\alpha)$ has two nearest neighbors $x,y \in M_\alpha$, the angle between $(u-x)$ and $(u-y)$ must be at most $2 \arcsin(\sqrt{1-\mu^2})$. As a result, a single branch of $M_\alpha$ between the two arcs of $G_{\alpha/2,R_0}$ cannot not generate any point of the $\mu$-medial axis, since $\alpha$ has been chosen so that $\alpha < 2 \arcsin(\sqrt{1-\mu^2})$.
Hence, for $\Delta,R_1$ large enough compared to $R_0$, we have $\rch_\mu(M_\alpha) \geq c'_{d,k} (\Delta \wedge R_1)$.

Finally, we build $M_\alpha'$ from $M_\alpha$ by bumping the curve near $G_{\alpha,R_0}$ as shown in Figure~\ref{fig:bump} (while still preserving the radial symmetry as before).
The manifold $M_\alpha'$ satisfies the same regularity conditions at $M_\alpha$.
Furthermore, $M_\alpha$ and $M_\alpha'$ only differ on a set of volume $\vol_d(M_\alpha \triangle M_\alpha') \leq C_{d,k} (\Delta\vee R_1)^{d-1} (R_0 \varepsilon)$.

\begin{figure}[h!]
\label{fig:bump}
\centering
\includegraphics[width=0.5\textwidth]{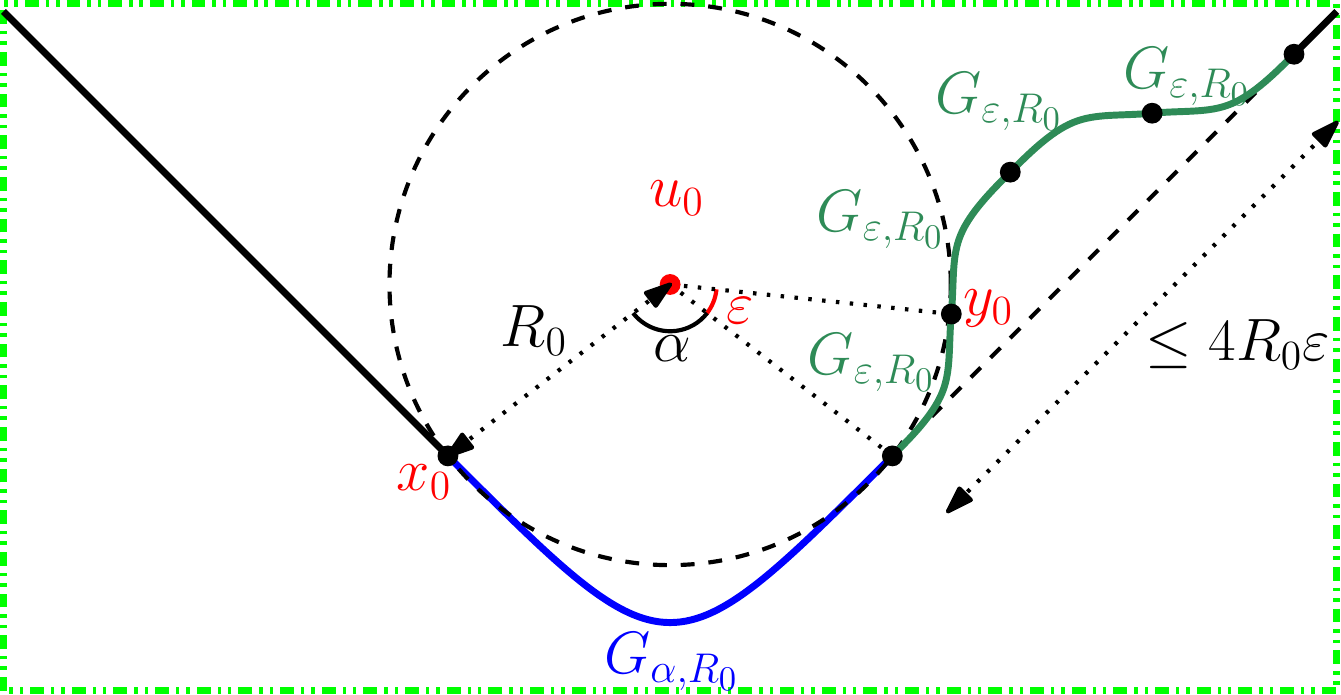}
\caption{
Local bump of $M_\alpha'$ for Proposition~\ref{prop:hypotheses},
in the boxed area of Figure~\ref{fig:malpha}.
}
\label{fig:bump}
\end{figure}

With this extra bump, we create a point $u_0 \in \Med(M_\alpha')$ that has two nearest neighbors $x_0,y_0 \in M_\alpha'$ at distance $R_0$, with angle between $(u_0-x_0)$ and $(u_0-y_0)$ equal to $\alpha' = \alpha + \varepsilon$, which satisfies $\sin(\alpha'/2)^2 = 1-\mu^2$. As a result, $u_0 \in \Med_\mu(M_\alpha')$, so that $\rch_\mu(M_\alpha') \leq \|u_0 - y_0\|=R_0$.
In particular, we have
\begin{align*}
|\rch_\mu(M_\alpha) - \rch_\mu(M_\alpha')| \geq c'_{d,k} (\Delta \wedge R_1) - R_0
.
\end{align*}

The proof is hence complete by setting $M = M_\alpha$ and $M'=M_\alpha'$, with $R_1 = \Delta = R_0 /c'_{d,k}$ and $R_0 = \rch_{\min} / c_{d,k}$ for small enough $c_{d,k},c'_{d,k}>0$.
\end{proof}

\begin{proof}[Proof of Theorem~\ref{thm:mu-reach-inconsistency}]
From Proposition~\ref{prop:hypotheses}, for $\varepsilon>0$ small enough, take $M,M' \in \manifolds{k}{\rch_{\min}}{\mathbf{L}}$ such that 
$
|\rch_\mu(M) - \rch_\mu(M')| \geq c_{d,k} \rch_{\min}
$
, $c'_{d,k} \rch_{\min}^d \leq \vol_d(M),\vol_d(M') \leq C_{d,k} \rch_{\min}^d$, and 
$\vol_d(M \triangle M') \leq C'_{d,k} \rch_{\min}^{d} \varepsilon$
.
Let us denote by $P$ and $P'$ the uniform distributions over $M$ and $M'$ respectively.
Elementary calculations directly yield that
$$
\TV(P,P') \leq \frac{\vol_d(M \triangle M')}{\vol_d(M) \vee \vol_d(M')} \leq C''_{d,k} \varepsilon
.
$$
Furthermore, since $c'_{d,k} \rch_{\min}^d \leq \vol_d(M) \wedge \vol_d(M') \leq \vol_d(M) \vee \vol_d(M') \leq C_{d,k} \rch_{\min}^d$, we obtain that $P,P' \in \distributions{k}{\rch_{\min}}{\mathbf{L}}{f_{\min}}{f_{\max}}$ as soon as $f_{\min} \leq 1/(C_{d,k} \rch_{\min}^d)$ and $f_{\max} \geq 1/(c'_{d,k} \rch_{\min}^d)$.
As a result, for all $n \geq 1$, Le Cam's Lemma \cite{yu1997assouad} yields
\begin{align*}
\inf_{\widehat{r}_\mu} \sup_{P \in \distributions{k}{\rch_{\min}}{\mathbf{L}}{f_{\min}}{f_{\max}}} 
&\E _{P^{\otimes n}}\left[
| \widehat{r}_\mu - \rch_\mu(M)| 
\right]
\\
&\geq
\frac{1}{2}
|\rch_\mu(M) - \rch_\mu(M')|
\bigl(1 - \TV(P,P'))^n
\\
&
\geq
c_{d,k} \rch_{\min} (1 - \varepsilon)^n
.
\end{align*}
As this construction is valid for all $\varepsilon>0$ small enough, we obtain the result by letting $\varepsilon$ tend to zero.
\end{proof}

\subsection{Maximal Curvature Estimation}\label{sec:proof_thm_cvrates_curvature_max}
This section is devoted to the proof of Theorem~\ref{thm:cv_rates_curvature_max}. It is based on a careful investigation of the local polynomial fitting procedure described in~\cite{Aamari19b}. First, recall that from~\cite[Lemma~2]{Aamari19b}, if $M \in \manifolds{k}{\rch_{\min}}{\mathbf{L}}$, $y \in M$ and $y' \in \ball \left (y, \frac{L_2 \wedge \rch_{\min}}{4} \right ) \cap M$, we may write
\begin{align}\label{eq:true_pol_decomposition}
y'-y = \pi^*_y&(y'-y) + \mathbb{T}_y^{(2),*}( \tens{\pi^*_y(y'-y)}{2}) + \hdots + \mathbb{T}_y^{(k-1),*}(\tens{\pi^*_y(y'-y)}{k-1}) \notag
		\\
		&+ R_{y}^{(k)}(y'-y),
\end{align}
where $\pi^*_y := \pi_{T_y M}$, $\mathbb{T}_y^{(j),*}$ are $j$-multilinear maps from $T_y M$ to $\R^D$, and $R_y^{(k)}$ satisfies
\begin{align*}
\norm{R_y^{(k)} (y'-y)} \leq C t_*^{k-1} \norm{y'-y}^{k},
\end{align*}
where $t_* = \max_{2 \leq j \leq k, y \in M}\|T^{(j),*}_y\|_{\mathrm{op}}^{\frac{1}{j-1}} \leq C_{k,d,\rch_{\min},\mathbf{L}}$. As assessed by~\cite[Lemma~2]{Aamari19b}, the polynomial decomposition expressed in~\eqref{eq:true_pol_decomposition} allows to recover the curvature tensor via $\II_y M = \mathbb{T}_{y}^{(2),*}$. Following~\cite{Aamari19b}, we estimate this curvature tensor via the second term of the polynomial decomposition provided by local fit to data points \eqref{eq:defi_pol_fit}.
To this aim, a slight adaptation of~\cite[Lemma~3]{Aamari19b} is needed, that allows to translate quality of approximation in terms of Hausdorff distance to guarantees on the monomial terms.
\begin{lemma}\label{lem:pol_expression}
          Set $h_0= (\tau_{min} \wedge L_2^{-1})/8$ and $h \leq h_0$. Let $M \in \mathcal{C}^{k}_{\tau_{min},\mathbf{L}}$, $x_0 = y_0 + z_0$, with $y_0 \in M$ and $\|z_0\| \leq \sigma \leq h/4$. Denote by $\pi^*_{y_0}$ the orthogonal projection onto $T_{y_0}M$, and by $\mathbb{T}_{y_0}^{(2),*}, \hdots, \mathbb{T}_{y_0}^{(k-1),*}$ the multilinear maps given by ~\eqref{eq:true_pol_decomposition}.
                    
          Let $x = y+z$ be such that    $y \in M$, $\|z\| \leq \sigma \leq h/4$ and $x \in \ball(x_0,h)$. We also let $\pi$ be an orthogonal projection,  and $\mathbb{T}^{(2)}, \hdots, \mathbb{T}^{(k-1)}$ be multilinear maps that satisfy 
          \begin{align*}
           \left (\max_{2 \leq j \leq k-1}\|\mathbb{T}^{(j)}\|_{\mathrm{op}}^{\frac{1}{j-1}} \right ) \vee t^* \leq t, \\
          th  \leq \frac{1}{4},
          \end{align*}
          for some $t \geq 0$. Then it holds 
           \begin{multline*}
        x-x_0 - \pi(x-x_0) - \sum_{j=2}^{k-1}{ \mathbb{T}^{(j)}(\tens{\pi(x-x_0)}{j})} = \sum_{j=1}^{k}\mathbb{T}^{(j),'}_{y_0}(\tens{\pi^*_{y_0}(y-y_0)}{j}) + R^{(k)}_{y_0}{(x-x_0)},    
\end{multline*}     
where $\mathbb{T}_{y_0}^{(j),'}$ are $j$-linear maps, and $\|R_{y_0}^{(k)}(x-x_0)\| \leq C \left ( \sigma + h^k (t_*^{k-1} + t^k h) \right )$, where $C$ depends on $d$, $k$, $\rch_{\min}$, $L_2$,$\ldots$,  $L_k$. Moreover, we have              
\[
    \begin{array}{@{}rl}
    \mathbb{T}^{(1),'}_{y_0}  &= (\pi^*_{y_0}-{\pi}),\\
    \mathbb{T}^{(2),'}_{y_0}  &= (\pi^*_{y_0}-{\pi})\circ \mathbb{T}_{y_0}^{(2),*} + (\mathbb{T}_{y_0}^{(2),*} \circ \pi_{y_0}^* - {\mathbb{T}}^{(2)} \circ {\pi}),
    \end{array}
    \]
    and, if $\pi = \pi^*_{y_0}$ and $\mathbb{T}^{(j)} = \mathbb{T}_{y_0}^{(j),*}$ for all $j \in \{2, \hdots, k-1\}$, then $\mathbb{T}_{y_0}^{(j),'}=0$ for all $j \in \{1, \hdots, k\}$.
          \end{lemma}
The proof of Lemma~\ref{lem:pol_expression} is deferred to Section~\ref{sec:proof_lem_pol_expression}. To ensure that our local curvature estimators allow to approximate the maximal curvature of $M$, we have to ensure that the sample covers $M$ well enough. That is the aim of the following Lemma.
\begin{lemma}[{\cite[Appendix, Lemma~B.7 \& Section 5.1.4]{Aamari19b}}]\label{lem:hausdorff_and_covering}
Let $P \in \distributions{k}{\rch_{\min}}{\mathbf{L}}{f_{\min}}{f_{\max}}$. Write $\mathbb{X}_n$ for an i.i.d. $n$-sample drawn from $P$. 
Let $h = \left( {C_{d,k}}\frac{f_{\max}^2}{f_{min}^3}  \frac{\log n}{n}\right)^{1/d}$, for $C_{d,k}$ large enough. Then, for $n$ large enough so that $h \leq \rch_{\min}/4$, with probability at least $1 - 2\left( \frac{1}{n}\right)^{2k/d}$, it holds
\begin{align*}
\dh\left(M,\mathbb{X}_n\right)
& \leq
h/4, \\
\dh (M,\hat{M}) & \leq C_{d,k,\rch_{\min}, \mathbf{L}} (t^*)^{k-1} \left ( \frac{f_{\max}^{2 + \frac{d}{2k}} \log n}{f_{\min}^{3 + \frac{d}{2k}} n} \right )^{k/d},
\end{align*}
where $\hat{M}$ denotes the union of local polynomial patches
\[
\hat{M} := \bigcup_{i=1}^n \hat{\Psi}_i \left ( \ball_{\hat{T}_i}(0,7h/8) \right )
\]
defined by \eqref{eq:defi_pol_fit} and \eqref{eq:polynomial_expansion},   and $t_* = \max_{y \in M, 2 \leq j \leq k}\|\mathbb{T}^{(j),*}_y\|_{\mathrm{op}}^{\frac{1}{j-1}} \leq C_{k,d,\rch_{\min},\mathbf{L}}$ as in Lemma~\ref{lem:pol_expression}.
\end{lemma}          
Equipped with these two lemmas, we are in position to prove Theorem~\ref{thm:cv_rates_curvature_max}.

\begin{proof}[Proof of Theorem~\ref{thm:cv_rates_curvature_max}]

Based on Lemma \ref{lem:hausdorff_and_covering}, for $h = \left ( C_{d,k} \frac{f^2_{\max} \log n}{f^3_{\min}n} \right )^{\frac{1}{d}}$, given $i \in \{1,\ldots,n\}$, we denote by $\hat{\Psi}_i$ the polynomial estimator around $X_i$ defined by
\begin{align*}
\hat{\Psi}_i(v) := X_i + v + \sum_{j=2}^{k-1} \hat{\mathbb{T}}_{i}^{(j)}(\tens{v}{j}),
\end{align*}
for all $v \in \hat{T}_i$. 
Setting
$$
\hat{M} := \bigcup_{i=1}^n \hat{\Psi}_i \left ( \ball_{\hat{T}_i}(0,7h/8) \right )
,
$$
we have that with probability larger than $1-2 \left( \frac{1}{n} \right )^{\frac{2k}{d}}$,
\begin{align}\label{eq:bound_hausdorff_patch}
\dh (\hat{M},M) \leq C_{d,k,\rch_{\min}, \mathbf{L}} (t_*)^{k-1} \left ( \frac{f_{\max}^{2 + \frac{d}{2k}} \log n}{f_{\min}^{3 + \frac{d}{2k}} n} \right )^{\frac{k}{d}} := \varepsilon_1,
\end{align}  
for $n$ large enough, according to Lemma~\ref{lem:hausdorff_and_covering}. In what follows we settle on the probability event of Lemma~\ref{lem:hausdorff_and_covering}. In particular, denoting by 
$$
\hat{t} = \max_{1 \leq i \leq n} \max_{2 \leq j \leq k-1} \|\hat{\mathbb{T}}_{i}^{(j)}\|_{\mathrm{op}}^\frac{1}{j-1},
$$
note that~\cite[Section 5.1.2]{Aamari19b} ensures that $\hat{t}\vee t_* \leq t \leq 1/(4h)$, for some fixed $t$, provided $n$ is large enough. 

We let $i \in \{1,\ldots,n\}$, $v \in \ball_{\hat{T}_i}(0,h/4)$, and intend to approximate $\II_{\pi_M(\hat{\Psi}_i(v))}$. To do so, we consider the following polynomial expansion centered at $v$: for $u \in \ball_{\hat{T}_i}(0,h/4)$,
\begin{align}\label{eq:pol_dec_1}
\hat{\Psi}_i(v+u) - \hat{\Psi}_i(v) & = u + \sum_{j=2}^{k-1} j \hat{\mathbb{T}}_{i}^{(j)}\left ( \tens{v}{j-1} \otimes u\right ) + \sum_{j=2}^{k-1} \sum_{r=j}^{k-1} \binom{r}{j} \hat{\mathbb{T}}_{i}^{(r)} \left ( \tens{v}{r-j} \otimes \tens{u}{j}\right ).
\end{align}
First we deduce from \eqref{eq:pol_dec_1} an estimate for the tangent space at $\pi_M(X_i+v)$, as well as a coordinate system. Namely, we let
\begin{align*}
  \hat{J}_{i,v}\colon \hat{T}_i &\longrightarrow \hat{J}_{i,v}(\hat{T}_i)
  \\
  u &\longmapsto u + \sum_{j=2}^{k-1} j \hat{\mathbb{T}}_{i}^{(j)} \left( \tens{v}{j-1} \otimes u \right).
\end{align*}
Note that since $th \leq 1/4$, we have 
\begin{align*}
\| \hat{J}_{i,v}(u) - u \| & \leq \sum_{j=2}^{k-1} j \left(\frac{th}{4} \right )^{j-1}\|u\| \\
                 & \leq \left (\sum_{j=1}^{\infty}j \left (\frac{th}{4} \right )^{j-1} -1 \right) \|u\| \\
                 & \leq \left ( \left ( \frac{1}{1-\frac{th}{4}} \right ) ^2 -1 \right ) \|u\| \leq \frac{\|u\|}{2},    
\end{align*}
so that $\hat{J}_{i,v}$ is full-rank. In what follows we write $\hat{T}_{i,v} := \Im(\hat{J}_{i,v})$  and $\hat{\pi}_{i,v} := \pi_{\hat{T}_{i,v}}$. We now may express \eqref{eq:pol_dec_1} in terms of the coordinate system given by $\hat{T}_{i,v}$:
\begin{align}\label{eq:pol_dec_2}
\hat{\Psi}_i(v+u) - \hat{\Psi}_i(v) = \hat{J}_{i,v}(u) + \sum_{j=2}^{k-1} \widetilde{\mathbb{T}}_{i,v}^{(j)}(\tens{\hat{J}_{i,v}(u)}{j}),
\end{align}
where the symmetric tensor  of order $j$ centered at $v$, $\widetilde{\mathbb{T}}_{i,v}^{(j)}$, is defined by
\begin{align*}
\widetilde{\mathbb{T}}_{i,v}^{(j)}(\tens{w}{j}) := \sum_{r=j}^{k-1} \binom{r}{j} \hat{\mathbb{T}}_{i}^{(r)} \left ( \tens{v}{r-j} \otimes \tens{\hat{J}_{i,v}^{-1}(w)}{j} \right ),
\end{align*}
for $w \in \hat{T}_{i,v}$. As well, since $th \leq \frac{1}{4}$, we may write
\begin{align*}
\norm{\widetilde{\mathbb{T}}_{i,v}^{(j)}}_{\mathrm{op}} & \leq \sum_{r=j}^{k-1} \binom{r}{j} (3/2)^j t^{r-1} \left ( \frac{h}{4} \right )^{r-j} \\ 
& \leq \left ( \sum_{r=j}^\infty \binom{r}{j} \left ( \frac{th}{4} \right )^{r-j} \right ) (3/2)^j t^{j-1} \\
&\leq \left ( \frac{1}{1-\frac{th}{4}} \right )^j (3/2)^j t^{j-1} \leq (3/2)^{2j} t^{j-1},
\end{align*}
so that $\max_{2 \leq j \leq k-1} \norm{\widetilde{\mathbb{T}}_{i,v}^{(j)}}_{\mathrm{op}}^{\frac{1}{j-1}} \leq \widetilde{t} \leq \left ( \frac{3}{2} \right )^4 t$.  

In particular, the bilinear form $\widetilde{\mathbb{T}}^{(2)}_{i,v} : \hat{T}_{i,v} \times \hat{T}_{i,v} \to \R^D$ may be expressed by
\begin{align*}
\widetilde{\mathbb{T}}_{i,v}^{(2)}(\tens{w}{2}) & := \sum_{j=2}^{k-1} \binom{j}{2} \hat{\mathbb{T}}_{i}^{(j)}\left ( \tens{v}{j-2} \otimes \tens{\hat{J}_{i,v}^{-1}(w)}{2} \right )
\end{align*}
for all $w \in \hat{T}_{i,v}$. Our second fundamental form estimator at $\pi_M(\hat{\Psi}_i(v))$ is then defined by
\begin{align*}
\hat{\mathbb{T}}_{i,v}^{(2)}: & = \widetilde{\mathbb{T}}_{i,v}^{(2)} \circ \hat{\pi}_{i,v} - \hat{\pi}_{i,v} \circ \widetilde{\mathbb{T}}_{i,v}^{(2)} \circ \hat{\pi}_{i,v},
\end{align*}
where with a slight abuse of notation, $\mathbb{T}\circ \pi (u) := \mathbb{T}\bigl(\tens{\pi(u)}{2}\bigr)$. Note that composition with $\hat{\pi}_{i,v}$ is performed to ensure that $\hat{\mathbb{T}}_{i,v}^{(2)}$ ranges into $\hat{T}_{i,v}^\perp$.

Our final max-curvature estimator can now be defined as
\begin{align*}
\hat{R}_\ell^{-1}
:=
\max_{1 \leq i \leq n} \max_{v \in \ball_{\hat{T}_i}(h/4)} \norm{\hat{\mathbb{T}}_{i,v}^{(2)}}_{\mathrm{op}}.
\end{align*}
First, we intend to show that, for a given  $v \in \ball_{\hat{T}_i}(h/4)$, $\hat{\mathbb{T}}_{i,v}^{(2)}$ is close to $\II_{y_0}$, for some $y_0 \in M$. To do so, we let $u \in \ball_{\hat{T}_i}(0,h/4)$, $x := \hat{\Psi}_i(v+u), x_0 := \hat{\Psi}_i(v)$, and $\widetilde{P}^{(r:k-1)}_{i,v} := \sum_{j=r}^{k-1} \widetilde{\mathbb{T}}^{(j)}_{i,v}$. Then, we have the decomposition
\begin{align*}
\hat{J}_{i,v}(u) & =  \hat{\pi}_{i,v}(x-x_0) - \sum_{j=2}^{k-1} \hat{\pi}_{i,v} \circ \widetilde{\mathbb{T}}^{(j)}_{i,v}(\tens{\hat{J}_{i,v}(u)}{j}) \\
& = \hat{\pi}_{i,v}(x-x_0) - \sum_{j=2}^{k-1} \hat{\pi}_{i,v} \circ \widetilde{\mathbb{T}}^{(j)}_{i,v}\left [\tens{\left (\hat{\pi}_{i,v}(x-x_0) -  \hat{\pi}_{i,v} \circ \widetilde{P}^{(2:k-1)}_{i,v}(\hat{J}_{i,v}(u)) \right )}{j} \right ] \\
& = \hat{\pi}_{i,v}(x-x_0) + \sum_{j=2}^{k} \mathbb{T}^{(j),''}_{i,v}( \tens{\hat{\pi}_{i,v}(x-x_0)}{j}) + R_{i,v}^{(k)}(x-x_0),
\end{align*}
with $\mathbb{T}^{(2),''}_{i,v} = - \hat{\pi}_{i,v} \circ \widetilde{\mathbb{T}}^{(2)}_{i,v}$, higher order tensors satisfying $\norm{\mathbb{T}^{(j),''}_{i,v}}_{\mathrm{op}} \leq C_k \widetilde{t}^{j-1} \leq C_k t^{j-1}$, and remainder term $\|R^{(k)}_{i,v}\| \leq C_k t^k h^{k+1}$. 
Plugging the above inequalities into~\eqref{eq:pol_dec_2} yields
\begin{align}\label{eq:pol_dec_decentre}
x-x_0 = \hat{\pi}_{i,v}(x-x_0) + \mathbb{T}^{(2)}_{i,v}\left(\tens{\hat{\pi}_{i,v}(x-x_0)}{2} \right ) + \sum_{j=3}^k \mathbb{T}^{(j)}_{i,v}(\tens{\hat{\pi}_{i,v}(x-x_0)}{j}) + R^{(k),'}_{i,v}(x-x_0),
\end{align} 
with $\mathbb{T}^{(2)}_{i,v} = \widetilde{\mathbb{T}}^{(2)}_{i,v} - \hat{\pi}_{i,v} \circ \widetilde{\mathbb{T}}^{(2)}_{i,v}$, $\norm{\mathbb{T}^{(j)}_{i,v}}_{\mathrm{op}} \leq C_k t^{j-1}$, and $\norm{R^{(k),'}_{i,v}(x-x_0)} \leq C_k t^k h^{k+1}$.

Then, according to Lemma \ref{lem:hausdorff_and_covering}, there exists $y_0 \in \ball(X_i, \frac{8}{7 \times 4}h) \cap M$ such that $\norm{y_0-x_0} \leq \varepsilon_1$, where $\varepsilon_1$ is defined by~\eqref{eq:bound_hausdorff_patch}. 
We further have
\begin{align*}
\norm{v - \hat{\pi}_{i}(y_0-X_i)} & \leq \varepsilon_1 + \norm{\hat{\Psi}_i(v) - (X_i + v)} \\
& \leq \varepsilon_1 + \norm{\sum_{j=2}^{k-1} \hat{\mathbb{T}}_{i}^{(j)}(\tens{v}{j})} \\ 
& \leq \varepsilon_1 + h/16 \\
& \leq h/8,
\end{align*}
since $\hat{t}h \leq \frac{1}{4}$, provided that $\varepsilon_1 \leq h/16$ (satisfied for $n$ large enough). 
Next, if $z \in \ball \left (y_0,\frac{h}{8} \right ) \cap M$, we have 
\begin{align*}
\norm{\hat{\pi}_{i}(z-X_i) - v } \leq & \norm{\hat{\pi}_{i}(z-y_0)} + \norm{v - \hat{\pi}_{i}(y_0-X_i)} \\
& \leq h/4,
\end{align*}
so that, writing $x_z:= \hat{\Psi}_i(\hat{\pi}_{i}(z-X_i))$, it holds $\norm{z-x_z} \leq \varepsilon_1$ and~\eqref{eq:pol_dec_decentre} applies.  
Next, provided $C_kth < 1/4$ and $C_kt \geq t_*$ (satisfied whenever $n$ is large enough), Lemma~\ref{lem:pol_expression} yields that 
\begin{multline*}
x_z - x_0 - \left ( \hat{\pi}_{i,v}(x_z-x_0) + \mathbb{T}^{(2)}_{i,v}(\tens{\hat{\pi}_{i,v}(x_z-x_0)}{2}) + \sum_{j=3}^k \mathbb{T}^{(j)}_{i,v}(\tens{\hat{\pi}_{i,v}(x_z-x_0)}{j}) \right ) \\
 = \sum_{j=1}^k \mathbb{T}^{(j),'}_{i,v}(\tens{\pi_y^*(z-y_0)}{j}) + R^{(k)}_{y_0}(x_z-x_0), 
\end{multline*}
so that 
\begin{align*}
\norm{\sum_{j=1}^k \mathbb{T}^{(j),'}_{i,v}(\tens{\pi_y^*(z-y_0)}{j})} = \norm{R^{(k),'}_{i,v}(x_z-x_0) - R^{(k)}_{y_0}(x_z-x_0)} \leq C_{k,d,\rch_{\min},\mathbf{L}} \varepsilon_1,
\end{align*}
according to \eqref{eq:pol_dec_decentre} and Lemma~\ref{lem:pol_expression}, since $t^kh \leq C_{k,d,\rch_{\min},\mathbf{L}}$.
Using the development~\eqref{eq:pol_dec_1} and the inclusion $\ball_{T_{y_0}M}(0,h/16) \subset \pi^*_{y_0} \left ( \ball(y_0,h/8) \cap M - y_0 \right )$ from~\cite[Lemma~2]{Aamari19b} then entails
\begin{align*}
\norm{\sum_{j=1}^k \mathbb{T}^{(j),'}_{i,v}(\tens{\pi^*_{y_0}(w)}{j})}  \leq C_{k,d,\rch_{\min},\mathbf{L}} \varepsilon_1,
\end{align*} 
for all $w \in \ball_{T_{y_0}M}(0,h/16)$. Proceeding as in~\cite[Proof of Theorem~2]{Aamari19b}, we get
\begin{align*}
\norm{\mathbb{T}^{(1),'}_{i,v}}_{\mathrm{op}} & \leq C_{k,d,\rch_{\min},\mathbf{L}} \varepsilon_1 h^{-1}, \\
\text{and~}\norm{\mathbb{T}^{(2),'}_{i,v}}_{\mathrm{op}} & \leq C_{k,d,\rch_{\min},\mathbf{L}} \varepsilon_1 h^{-2}.
\end{align*}
In turn, following~\cite[Proof of Theorem~4]{Aamari19b} entails
\begin{align*}
\norm{\hat{\mathbb{T}}^{(2)}_{i,v} \circ \hat{\pi}_{i,v} - \mathbb{T}^{(2),*}_{y_0} \circ \pi^*_{y_0}}_{\mathrm{op}} \leq C_{k,d,\rch_{\min},\mathbf{L}} \varepsilon_1 h^{-2}.
\end{align*}
Since $\II_{y_0} = \mathbb{T}^{(2),*}_{{y_0}}$ (\cite[Lemma~2]{Aamari19b}), we deduce that 
\begin{align}\label{eq:sup_curv_1}
\max_{1 \leq i \leq n} \max_{v \in \ball_{\hat{T}_j}(0,h/4)} \norm{\hat{\mathbb{T}}^{(2)}_{i,v} \circ \hat{\pi}_{i,v}}_{\mathrm{op}} \leq \max_{y \in M} \norm{\II_y}_{\mathrm{op}} + C_{k,d,\rch_{\min},\mathbf{L}} \varepsilon_1 h^{-2}. 
\end{align}
Conversely, since $X_1, \hdots, X_n$ is a $(h/4)$-covering of $M$ onto the probability event described in Lemma~\ref{lem:hausdorff_and_covering}, we deduce that for all $y \in M$, there exists $i_0 \in \{1,\ldots,n\}$ such that $\|X_{i_0} - y \| \leq h/4$. In particular, we have 
\begin{align*}
v:= \hat{\pi}_{i_0,v}(y - X_{i_0}) \in \ball_{\hat{T}_{i_0}}(0,h/4).
\end{align*}
Proceeding as above similarly leads to 
\begin{align*}
\norm{\hat{\mathbb{T}}^{(2)}_{i_0,v} \circ \hat{\pi}_{i_0,v} - \II_y \circ \pi^*_y}_{\mathrm{op}} \leq C_{k,d,\rch_{\min},\mathbf{L}} \varepsilon_1 h^{-2},
\end{align*}
so that 
\begin{align}\label{eq:sup_curv_2}
\max_{y \in M} \norm{\II_y}_{\mathrm{op}} \leq \max_{1 \leq i \leq n} \max_{v \in \ball_{\hat{T}_i}(0,h/4)} \norm{\hat{\mathbb{T}}^{(2)}_{i,v} \circ \hat{\pi}_{i,v}}_{\mathrm{op}}  + C_{k,d,\rch_{\min},\mathbf{L}} \varepsilon_1 h^{-2}.
\end{align}
Combining~\eqref{eq:sup_curv_1} and~\eqref{eq:sup_curv_2} yields that for $n$ large enough, 
\begin{align*}
\bigl| \hat{R}_\ell - R_\ell(M) \bigr| \leq R_\ell(M)^2  C_{k,d,\rch_{\min},\mathbf{L}} \varepsilon_1 h^{-2}
,
\end{align*}
which concludes the proof.
\end{proof}

\subsection{Proof of Lemma~\ref{lem:pol_expression}}\label{sec:proof_lem_pol_expression}                   
          \begin{proof}[Proof of Lemma~\ref{lem:pol_expression}] We follow the proof of~\cite[Lemma~3]{Aamari19b}. Without loss of generality we take $y_0=0$, so that $\|y\| \leq 3h/2$. Let $z' = z-z_0$, so that $\|z'\| \leq h/2$. We write
          \begin{align*}
          & x-x_0 - \pi(x-x_0) - \sum_{j=2}^{k}{ \mathbb{T}^{(j)}(\tens{\pi(x-x_0)}{j})} \\ 
          & \qquad  = y + z' - \pi(y+z') - \sum_{j=2}^{k}\mathbb{T}^{(j)} (\tens{ (\pi(y) + \pi(z')) }{j}) \\
          & \qquad  = y + z' - \pi(y+z') - \sum_{j=2}^{k} \left [ \mathbb{T}^{(j)}(\tens{\pi(y)}{j}) + \sum_{r=0}^{j-1} \binom{j}{r} \mathbb{T}^{(j)} \left ( \tens{\pi(y)}{r} \otimes \tens{\pi(z')}{j-r} \right ) \right ].
          \end{align*}
          Since, for any $j \geq 2$ and $r \in \{0,\ldots,j-1\}$, 
          \begin{align*}
          \norm{ \mathbb{T}^{(j)} \left ( \tens{\pi(y)}{r} \otimes \tens{\pi(z')}{j-r} \right )} & \leq t^{j-1} (3h/2)^{r} (2 \sigma)^{j-r} \\ 
          & \leq C_k \sigma t^{j-1}h^{j-1} \leq C_k \sigma,  
          \end{align*}
          we may write
          \begin{multline}\label{eq:dec_pol_noise}
          x-x_0 - \pi(x-x_0) - \sum_{j=2}^{k}{ \mathbb{T}^{(j)}(\tens{\pi(x-x_0)}{j})} \\ 
          = y  - \pi(y) - \sum_{j=2}^{k} \mathbb{T}^{(j)}(\tens{\pi(y)}{j}) + R^{(k),'}(x-x_0),
          \end{multline}
          where $\norm{R^{(k),'}(x-x_0)} \leq C_k \sigma$. Next,~\eqref{eq:true_pol_decomposition} entails
           \begin{align*}
y = \pi_{y_0}^*&(y) + \mathbb{T}^{(2),*}_{y_0}( \tens{\pi_{y_0}^*(y)}{2}) + \hdots + \mathbb{T}_{y_0}^{(k-1),*}(\tens{\pi^*_{y_0}(y)}{k-1}) 
		\\
		&+ R^{(k),''}_{y_0}(y), 
\end{align*}
           with $\| R^{(k),''}_{y_0}(y) \| \leq C_{k,d,\rch_{\min}, \mathbf{L}} t_*^{k-1} h^k$. Denoting by 
           \[
           {P}^{*,(1:k-1)}_{y_0}(\pi^*_{y_0}(y)) := \pi^*_{y_0}(y) + \sum_{r=2}^{k-1} \mathbb{T}_{y_0}^{(r),*}(\tens{\pi^*_{y_0}(y)}{r}),
           \]
            we deduce that 
          \begin{multline*}
          y - \pi(y) - \sum_{j=2}^{k}{ \mathbb{T}^{(j)}(\tens{\pi(y)}{j})} = \\
          {P}^{*,(1:k-1)}_{y_0}(\pi^*_{y_0}(y)) + R^{(k),''}_{y_0}(y) - \pi \left ( {P}^{*,(1:k-1)}_{y_0}(\pi^*_{y_0}(y)) + R^{(k),''}_{y_0}(y) \right ) \\
          - \sum_{j=2}^{k} \mathbb{T}^{(j)} \left ( \tens{ \pi \left [ \left ({P}^{*,(1:k-1)}_{y_0}(\pi^*_{y_0}(y)) + R^{(k),''}_{y_0}(y) \right )\right ]}{j} \right )
          .
          \end{multline*}
Note that 
$$
\norm{\pi(R^{(k),''}_{y_0}(y))} \leq \norm{R^{(k),''}_{y_0}(y)} \leq C_{k,d,\rch_{\min}, \mathbf{L}} t_*^{k-1} h^k
.
$$ 
Next, since $\|y\| \leq 3h/2$, it holds
\begin{align*}
\norm{ {P}^{*,(1:k-1)}_{y_0}(\pi^*_{y_0}(y)) } & \leq \sum_{r=1}^{k-1} t_*^{r-1} \left ( \frac{3h}{2} \right )^r  \\
& \leq \frac{3h}{2} \frac{1}{1-\frac{3t_*h}{2}} \leq 3h,
\end{align*}          
so that, for all $j \in \{2,\ldots,k\}$, 
\begin{multline*}
\left\Vert
\mathbb{T}^{(j)} \left ( \tens{ \pi \left [ \left ( {P}^{*,(1:k-1)}_{y_0}(\pi^*_{y_0}(y)) + R^{(k),''}_{y_0}(y) \right ) \right ]}{j} \right) 
- 
\mathbb{T}^{(j)} \left ( \tens{ \pi \left [ \left ( {P}^{*,(1:k-1)}_{y_0}(\pi^*_{y_0}(y)) \right ) \right ]}{j} \right)
\right\Vert \\
\begin{aligned}[t]
& \qquad \leq  t^{j-1} \sum_{r=1}^j \binom{j}{r} \norm{R^{(k),''}_{y_0}}^r (3h)^{j-r} \\
& \qquad \leq C_{k,d,\rch_{\min}, \mathbf{L}} t^{j-1}h^j\max_{1 \leq r \leq j} t_*^{(k-1)r}h^{(k-1)r} \\
& \qquad \leq C_{k,d,\rch_{\min}, \mathbf{L}} t^k h^{k+1}. 
\end{aligned}
\end{multline*}
Thus, we may write
\begin{multline*}
y - \pi(y) - \sum_{j=2}^{k}{ \mathbb{T}^{(j)}(\tens{\pi(y)}{j})}  = {P}^{*,(1:k-1)}_{y_0}(\pi^*_{y_0}(y)) - \pi \left ( {P}^{*,(1:k-1)}_{y_0}(\pi^*_{y_0}(y)) \right ) \\
 - \sum_{j=2}^{k} \mathbb{T}^{(j)} \left ( \tens{ \pi \left [ \left ( {P}^{*,(1:k-1)}_{y_0}(\pi^*_{y_0}(y))\right ) \right ]}{j} \right ) + R^{(k),'''}_{y_0}(y),
\end{multline*}
where $\norm{R^{(k),'''}_{y_0}(y)} \leq C_{k,d,\rch_{\min}, \mathbf{L}} h^k (t_*^{k-1} + t^k h)$. At last, for $j \in \{2,\ldots,k\}$, and $r_1, \hdots, r_j \in \{1,\ldots,k-1\}$ such that $\sum_{s=1}^{j} r_s \geq k+1$, we have 
\begin{align*}
\norm{\mathbb{T}^{(j)} \left ( \bigotimes_{s=1}^j \pi \left ( \mathbb{T}_{y_0}^{(r_s),*} \left( \tens{\pi^*_{y_0}(y)}{r_s}\right ) \right ) \right )} & \leq t^{j-1} \prod_{s=1}^j t_*^{r_s-1}h^{r_s} \\
& \leq (th)^{\left ( \sum_{s=1}^j r_s \right )-1} h \leq t^k h^{k+1}, 
\end{align*}
where $\mathbb{T}_{y_0}^{(1),*} = \pi_{y_0}^*$, with a slight abuse of notation. Hence, it holds 
\begin{multline*}
y - \pi(y) - \sum_{j=2}^{k}{ \mathbb{T}^{(j)}(\tens{\pi(y)}{j})} = 
( \pi^*_{y_0} - \pi\circ \pi^*_{y_0})(y) + \mathbb{T}_{y_0}^{(2),*}(\tens{\pi_{y_0}^*(y)}{2}) -  \pi \left ( \mathbb{T}_{y_0}^{(2),*}(\tens{\pi^*(y)}{2}) \right ) \\ 
- \mathbb{T}^{(2)} \left ( \tens{\left ( \pi \circ \pi_{y_0}^*(y) \right )}{2} \right )  + \sum_{j=3}^k \mathbb{T}^{(j),'}_{y_0} \left ( \tens{\pi_{y_0}^*(y)}{j} \right ) + R^{(k),''''}_{y_0}(y),
\end{multline*}
where $\norm{R^{(k),''''}_{y_0}(y)} \leq C_{k,d,\rch_{\min}, \mathbf{L}} h^k (t_*^{k-1} + t^k h)$. Plugging the above equation into \eqref{eq:dec_pol_noise} gives the result.
          \end{proof}

\section{Proofs of \secref{sdr}} \label{sec:proofsdr}

\subsection{Comparing Reach, Weak Feature Size and Spherical Distortion Radius}
\label{sec:proof-compare-reach-sdr}

Let us prove \prpref{interpolate}.

\begin{proof}[Proof of \prpref{interpolate}]
The monotonicity follows trivially from the definition, and since by~\cite[Theorem~1]{boissonnat2019reach},  $\sdr_0(K,\d_K) = \rch(K,\d_K)$, there holds immediately that $\sdr_\delta(K,\d_K) \geq \rch(K)$ for any $\delta \geq 0$. 
Now take $\delta \leq \sqrt{2(D+1)/D}\wfs(K)$, and take $z$ a critical point of $K$, so that $z \in \conv \Gamma$ where $\Gamma := \{x \in K~|~\|x-z\| = \d(z,K)\}$. Using Jung's theorem~\cite[Theorem~2.10.41]{Federer69}, there holds
$$
\diam(\Gamma) \geq \sqrt{\frac{2(D+1)}{D}} \rad(\Gamma) = \sqrt{\frac{2(D+1)}{D}} \d(z,K) \geq \sqrt{\frac{2(D+1)}{D}} \wfs(K) \geq \delta
$$
so that there exists two points $x,y \in \Gamma$ such that $\|x-y\| \geq \delta$. Furthermore, since the interior of $\ball(z,\wfs(K))$ contains no point of $K$, there holds
$$
\d_K(x,y) \geq \d_{\cS(\wfs(K))}(x,y) > \dsr(x,y),\text{~for all~} r > \wfs(K),
$$
so that indeed $\sdr_\delta(K,\d_K) \leq \wfs(K)$.
\end{proof}

\subsection{Stability Properties of the Spherical Distortion Radius}
\label{sec:proof-sdr-properties}
We now move to the proofs of the stability properties of the SDR. As a first step, we will need the following lemma on geodesic distances over spheres.

\begin{lemma} \label{lem:dstech} Let $r,\ve > 0$ and take $x,y,a,b \in K$ such that $\|x-y\| < 2r$ and 
$$
\|a-b\| \leq \left(1+\frac{A\ve}{r}\right) \|x-y\|
$$
for some $A > 0$. For all $\lambda > 0$, define
$$
\zeta_\lambda = \max\{\frac{192 r^3}{\|a-b\|^3}(\lambda + A\pi), 4A \}.
$$
Then, for all $\zeta \geq \zeta_\lambda$ such that $\zeta \ve \leq r$, there holds
\[
\d_{\cS(r+ \zeta \ve)}(a,b) \leq \dsr(x,y) - \lambda \ve.
\]
\end{lemma}
\begin{proof}[Proof of Lemma~\ref{lem:dstech}] Notice that, denoting by $\rho = \|x-y\|$,
\[
\dsr(x,y) = 2 r \arcsin\(\frac{\rho}{2r}\) = \rho \times \vp(2r/\rho)~~\text{with}~~\vp(u) := u \arcsin(1/u).
\]
The map $\vp$ is decreasing on $[1,\infty)$ and, using the development of 
\[
\arcsin(u) = \sum_{n = 0}^\infty (2n)! u^{2n+1}/(2^{2n}n!^2(2n+1)),
\] 
we find that
$$
\vp'(u) = - \sum_{n = 1}^\infty  \frac{(2n)! \times 2n }{2^{2n}n!^2(2n+1)}\frac{1}{u^{2n+1}} \leq - \frac{1}{3u^3}.
$$
Notice furthermore that, by assumption
\begin{align*} 
\frac{2(r+\zeta\ve)}{\|a-b\|} &\geq  \frac{2(r+\zeta\ve)}{(1+A\ve/r)\|x-y\|} = \frac{1+\zeta\ve/r}{1+A\ve/r} \frac{2r}{\|x-y\|} \\
&\geq \left(1+\frac{\zeta\ve}{2r}\right) \frac{2r}{\|x-y\|}
\end{align*} 
where we used that $A \leq \zeta/4$, $\zeta \ve \leq r$, and that $(1+u)/(1+u/4) \geq 1+u/2$ for $|u| \leq 2$. 
Now, as $\vp \leq \pi/2$ and that $|\vp'|$ is decreasing, we can write 
\begin{align*} 
\d_{\cS(r+\zeta\ve)}(a,b) 
&\leq 
\frac{A\ve}{r}\|x-y\| \vp\left(2(r+\zeta\ve)/\|a-b\|\right)+ \|x-y\|\vp\left(2(r+\zeta\ve)/\|a-b\|\right) \\
&\leq A\pi\ve+ \dsr(x,y)
\\
&\hspace{3em}
-\|x-y\| \times | \vp' | \(\frac{2(r+\zeta\ve)}{\|a-b\|}\)\times \left(\frac{2(r+\zeta\ve)}{\|a-b\|} -\frac{2r}{\|x-y\|} \right) \\
&\leq \dsr(x,y)+A\pi\ve-\frac{\|a-b\|^3}{3(2(r+\zeta\ve))^3} \zeta\ve\\
&\leq  \dsr(x,y)+\left(A\pi-\frac{\|a-b\|^3}{192 r^3} \zeta\right)\ve,
\end{align*} 
and using $\zeta \geq \zeta_\lambda$ ends the proof.
\end{proof}

We are now in position to prove \prpref{stab1} and  \thmref{lip}.

\begin{proof}[Proof of \prpref{stab1}] If $\mathrm{r}_1 = \infty$ there is nothing to show. Otherwise, notice that because $\mathrm{r}_1 \geq \delta_0/2$ by definition, there holds that
\[
\xi(R) \geq \max\{\frac{192 R^3}{\delta_0^3}\(1+\pi\frac{2R}{\delta_0}\), \frac{8 R}{\delta_0}\} \geq 1 
\]
for all $R > \mathrm{r}_1$. Now, since $\xi(\mathrm{r}_1) \Upsilon < \mathrm{r}_1$, one can find $R > \mathrm{r}_1$ such that $\xi(R) \Upsilon < R$. By definition of $\mathrm{r}_1$, there exist $x,y \in K'$ such that $\delta+2\ve \leq \|x-y\| < 2R$ and $\d_{\cS(R)}(x,y) < \d'(x,y)$. 
Now, let $a,b \in K$ be two closest points (in Euclidean distance) from $x$ and $y$ such that $\d(a,b) = \d\(\pr_K(\{x\}), \pr_K(\{y\}\)$. Then 
\[
\delta \leq \|a-b\| \leq \|x-y\|+2\ve < 2R + 2\Upsilon \leq 2(R+\xi(R)\Upsilon) 
\]
and 
\[
\|a-b\| \leq  \|x-y\|+2\ve \leq \left(1 + \frac{2\Upsilon R}{\delta_0 R}\right) \|x-y\|.
\]
We now can apply \lemref{dstech} with $A = 2R/\delta_0$ and $\lambda = 1$ to find that
\begin{align*} 
\d(a,b) 
&\geq 
\frac{1}{1+\nu} \d'(x,y) 
\\
&> \frac{1}{1+\nu}\d_{\cS(R)}(x,y)  
\\
&\geq 
\frac{1}{1+\nu}\left(\d_{\cS(R+\xi(R)\Upsilon)}(a,b) + \Upsilon\right)
\\
&\geq 
\frac{1+\Upsilon/\delta}{1+\nu} \d_{\cS(R+\xi(R)\Upsilon)}(a,b),
\end{align*} 
where the last inequality uses that $\d_{\cS(R+\xi(R)\Upsilon)}(a,b) \geq \|a-b\| \geq \delta$. At the end of the day, since $\Upsilon \geq \delta\nu$, we have $\d(a,b)  >  \d_{\cS(R+\xi(R)\Upsilon)}(a,b)$, so that $\sdr_{\delta}(K,\delta) < R + \xi(R)\Upsilon$. Taking $R$ to $\mathrm{r}_1$ yields the result.
\end{proof}

\begin{proof}[Proof of \thmref{lip}] We take $\ve > 0$ such that
$$\ve <  C_0\ve_0, ~~\ve < (\delta_1-\delta_0)/2,~~\text{and}~~\ve < r_0/L_0,$$ 
and take $\delta \in [\delta_0,\delta_1-\ve)$. We write $r_\delta := \sdr_\delta(K,\d)$ and $r_{\delta+\ve} := \sdr_{\delta+\ve}(K,\d)$ for short. Recall that $r_\delta \leq r_{\delta+\ve}$. Now take $r \leq r_{\delta+\ve} - L_0\ve$, and two points $x,y \in K$ such that $\delta \leq \|x-y\| < 2r$ (if there are none, then $r \leq r_\delta$ automatically). If $\|x-y\| \geq \delta + \ve$, then $\d(x,y) \leq \d_{\cS(r)}(x,y)$ because $r \leq r_{\delta+\ve}$. If now $\|x-y\| < \delta+\ve$, since $\|x-y\| \leq \Delta_0$, we can use \assref{spread} and find a point $a \in K$ such that $\|a-y\| \leq \ve/C_0$ and $\|x-a\| \geq \|x-y\|+ \ve \geq \delta+\ve$. Now, since $r + L_0 \ve \leq r_{\delta+\ve}$, it holds $\d(x,a) \leq \d_{\cS(r+L_0\ve)}(x,a)$. Furthermore, notice that
$$
\|x-a\| \leq \|x-y\| + \frac1{C_0}\ve \leq \left(1 + \frac{r_1 \ve}{C_0 \delta_0 r}\right) \|x-y\| 
.
$$
Using \assref{subeuc} and \lemref{dstech} with $A = r_1 /(C_0 \delta_0)$ and $\lambda = C_1/C_0$, we find
\begin{align*} 
\d(x,y) \leq \d(x,a) + \d(a,y) \leq  \d_{\cS(r+L_0\ve)}(x,a)+ \frac{C_1}{C_0} \ve \leq \dsr(x,y)
,
\end{align*} 
so that in the end $r \leq r_\delta$. Taking $r$ to $r_{\delta+\ve} -L_0\ve$ yields that $r_{\delta+\ve} \leq r_{\delta} + L_0\ve$, ending the proof. 
\end{proof}

Finally, \corref{stab} follows as a direct corollary of \prpref{stab1}.

\begin{proof}[Proof of \corref{stab}] Since $\xi_0 \ve \leq  2\sdr_{\delta_1}(K,\d)$, the radius $\sdr_{\delta_1}(K,\d)$ is in particular finite so that, according to \prpref{stab1}, $\sdr_{\delta}(K',\d') \leq 2 \sdr_{\delta_1}(K,\d)$ and, consequently, $\xi_1 \Upsilon \leq  \sdr_{\delta}(K',\d')$ and $\xi_2 \Upsilon \leq  \sdr_{\delta+2\ve}(K,\d)$, where $\xi_1 = \xi(\sdr_{\delta}(K',\d'))$ and $\xi_2 = \xi( \sdr_{\delta+2\ve}(K,\d))$. Applying \prpref{stab1} twice -- which is possible, since $\Upsilon \geq ((\delta-2\ve)\nu) \vee \ve)$ --, we thus find
$$
\sdr_{\delta-2\ve}(K,\d)-\xi_1 \Upsilon \leq \sdr_{\delta}(K',\d') \leq \sdr_{\delta+2\ve}(K,\d)+\xi_2 \Upsilon, 
$$
and we conclude by noticing that both $\xi_1$ and $\xi_2$ are less than $\xi_0$.
\end{proof}

\begin{proof}[Proof of \thmref{stab}] Using \corref{stab} and \thmref{lip}, one find that
\begin{align*} 
\sdr_{\delta}(K',\d') 
&\leq 
\sdr_{\delta+2\ve}(K,\d)+\xi_0 \Upsilon 
\\
&\leq  
\sdr_{\delta}(K,\d)+2L_0 \ve + \xi_0 \Upsilon 
\\
&\leq  
\sdr_{\delta}(K,\d)+ \zeta_0 \Upsilon,
\end{align*} 
and likewise for the lower bound.
\end{proof}

\section{Proofs of \secref{metric}} \label{sec:proofmetric}

\subsection{Minimax Lower Bound for Metric Learning}
\label{sec:proof-lower-bound-metric}
We now turn towards the proof of \thmref{metriclb}. 
It relies on an adaptation of the classical Le Cam's argument~\cite{yu1997assouad} to the asymmetric loss $\ell_\infty$.

\begin{lemma} \label{lem:lecammet} Let $x,y \in \bbR^D$ and let $M_0$ and $M_1$ be two submanifolds of $\bbR^D$ such that $x, y \in M_0 \cap M_1$ and the uniform distribution $P_0$ (resp. $P_1$) on $M_0$ (resp. $M_1$) is in $\distributions{k}{\rch_{\min}}{\mathbf{L}}{f_{\min}}{f_{\max}}$. Then if $\d_{M_0}(x,y) \leq \d_{M_1}(x,y)$, 
\beq
 \inf_{\wh\d} \sup_{P \in \cP^k} \bbE_{P^{\otimes n}}[\ell_\infty(\wh \d|\d_M)] \geq \frac12 \times \left|1- \frac{\d_{M_0}(x,y) }{\d_{M_1}(x,y)}\right| \times (1 - \tv(P_0^{\otimes n},P_1^{\otimes n})), \label{eq:lecammet}
\eeq
\end{lemma}

\begin{proof}[Proof of \lemref{lecammet}] For brevity, we write $\cR_n$ be the minimax risk appearing in the left-hand side of~\eqref{eq:lecammet}. First, we write
\begin{align*} 
\cR_n 
&\geq  
\inf_{\wh \d} \sup_{P \in \{P_0, P_1\}} \bbE_{P^{\otimes n}}[\ell_\infty(\wh \d | \d_M)]  \\
&\geq
\inf_{\wh \d} \sup_{P \in \{P_0, P_1\}} \bbE_{P^{\otimes n}}  \left[\left|1- \frac{\wh \d(x,y) }{\d_M(x,y)}\right|\right] \\
&\geq \frac12 \inf_{\wh \d}\{\bbE_{P_0^{\otimes n}} \left[ \left|1- \frac{\wh \d(x,y) }{\d_{M_0}(x,y)}\right|\right] +\bbE_{P_1^{\otimes n}} \left[ \left|1- \frac{\wh \d(x,y) }{\d_{M_1}(x,y)}\right|\right]  \}\\
&\geq \frac12 \inf_{\wh \d} \bbE_{P_0^{\otimes n}}\left[\left(\left|1- \frac{\wh \d(x,y) }{\d_{M_0}(x,y)}\right| + \left|1- \frac{\wh \d(x,y) }{\d_{M_1}(x,y)}\right|\right) \times \left(1 \wedge \frac{\diff P_1^{\otimes n}}{\diff P_0^{\otimes n}}\right)\right].
\end{align*} 
But now, using that $\d_{M_0}(x,y) \leq \d_{M_1}(x,y)$, a simple computation shows that the functional 
$$
\delta \mapsto \left|1- \frac{\delta}{\d_{M_0}(x,y)}\right| + \left|1- \frac{\delta}{\d_{M_1}(x,y)}\right|
$$
is minimal for $\delta = \d_{M_0}(x,y)$ so that
\begin{align*} 
\cR_n  &\geq \frac12  \bbE_{P_0^{\otimes n}}\left[\left|1- \frac{\d_{M_0}(x,y) }{\d_{M_1}(x,y)}\right| \times \left(1 \wedge \frac{\diff P_1^{\otimes n}}{\diff P_0^{\otimes n}}\right)\right] \nonumber\\
&= \frac12 \times \left|1- \frac{\d_{M_0}(x,y) }{\d_{M_1}(x,y)}\right| \times (1 - \tv(P_0^{\otimes n},P_1^{\otimes n})),
\end{align*} 
which ends the proof. 
\end{proof}

\begin{proof}[Proof of \thmref{metriclb}] 
Without loss of generality, we set the analysis in $\bbR^{d+1} \simeq \bbR^{d+1} \times \{0\}^{D-(d+1)} \subset \bbR^D$. 

\subsubsection{Submanifolds Construction} We let $M_0 \subset \bbR^{d+1}$ be a submanifold of $\cC^k_{2\rch_{\min},\mathbf{L}/2}$ such that it contains the cylinder 
$$
\{(s,z) \in \bbR^{2} \times \bbR^{d-1}~|~\|s\| = R \text{~and~} \|z\| \leq 3R\}
.
$$
Such a manifold always exists as soon as $R \geq 2 \rch_{\min}$ and $L_j$ is large enough compared to $1/R^{j-1}$.
For instance, one can design $M_0$ as a hypersurface of revolution obtained based on patches the interpolating curves of \lemref{turn-widget}.

In what follows, we denote any $x \in \bbR^{d+1} = \bbR^{d} \times \bbR$ as $x = (w,h) \in \bbR^{d}\times \bbR$. 
With this notation, we define,  for $\ve>0$ and $c>0$ to be chosen later,
$$
\Phi_{\ve}(x) := x + c \ve^{k} K(w/\ve) e_{d+1}~~~\text{where}~~~ e_{d+1} = (0,\dots,0,1) \in \bbR^{d+1},
$$
where $K(w)$ equals $\exp(-1/(1-\|w\|^2)_{+})$ for $\|w\|<1$ and $0$ otherwise.

For $\ve \leq 1$ and $c$ small enough, $\Phi_{\ve}$ is a diffeomorphism of $\bbR^{d+1}$ with derivative bounded up to the order $k$. 
Using~\cite[Proposition~A.4]{Aamari19}, we get that $M_\ve := \Phi_{\ve}(M)$, the image of $M_0$ by $\Phi_{\ve}$, belongs to $\cC^k_{\rch_{\min},\mathbf{L}}$ provided that $c$ is small enough (depending on $R$) and $\ve \leq c R$.

Locally around the apex $(0,R) \in \bbR^{d+1}$, $M_0$ can be seen as the graph of $
\Psi_0(w) := \sqrt{R^2-w_1^2}$, defined on $(-R,R) \times \ball_{\bbR^{d-1}}(0,3R)$, while $M_\ve$ is the graph of
$$
\Psi_\ve(w) := \Psi_0(w) + c\ve^{k} K(w/\ve).
$$
Finally, we let $\bar\Psi_\ve(w) := (w,\Psi_\ve(w))$ and similarly define $\bar\Psi_0$. We refer to \figref{cyl} for a diagram of the situation.

\begin{figure}[h!]
\centering
\begin{subfigure}{0.49\textwidth}
	\centering
	\includegraphics[height = 0.7\linewidth,page=1]{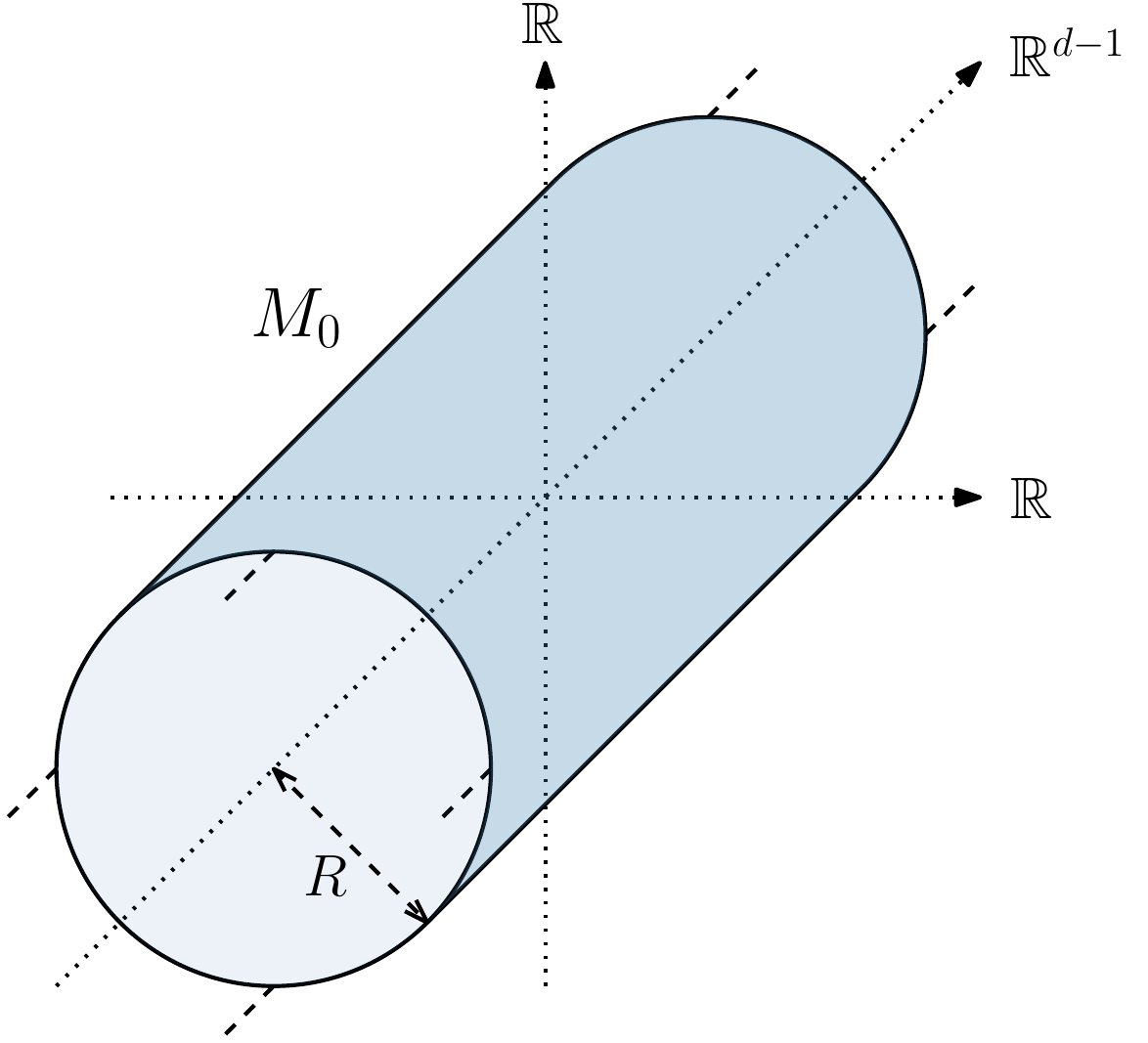}
	\caption{}
	\label{fig:Mgeo0}
\end{subfigure}
\begin{subfigure}{0.49\textwidth}
	\centering
	\includegraphics[height = 0.7\linewidth,page=2]{Mgeo}
	\caption{}
	\label{fig:Mgeo1}
\end{subfigure}
\caption{(\subref{fig:Mgeo0}) The cylindrical section of $M_0$ used in the proof of \thmref{metriclb}, and (\subref{fig:Mgeo1}) the perturbed submanifold $M_\ve$.}
\label{fig:cyl}
\end{figure}

\subsubsection{Shortest-Path Properties}
In this section, we seek to derive a lower bound on $\left|1 - \d_{M_0}(x,y)/\d_{M_\ve}(x,y)\right|$, so as to apply Lemma~\ref{lem:lecammet}. For this, we will consider well-chosen $x,y \in M_0 \cap M_\ve$ and derive a lower bound on $\d_{M_\ve}(x,y) - \d_{M_0}(x,y)$.

We let $\ell < R$, and we pick $x := \bar\Psi_0(-\ell e_1)$ and $y := \bar\Psi_0(\ell e_1)$ where $e_1 = (1,0,\dots,0) \in \bbR^d$. By construction, $x$ and $y$ belong to $M_0$.
Furthermore, provided that $\ell \geq \ve$, there holds that $x =  \bar\Psi_\ve(-\ell e_1)$ and $y = \bar\Psi_\ve(\ell e_1)$ so that $x$ and $y$ are also in $M_\ve$. 
We let $\gamma_\ve : [-1,1] \to M_{\ve}$ be a shortest path in $M_\ve$ between $x$ and $y$, parametrized at constant speed. We denote paths
$$
w_\ve := a_\ve e_1 + b_\ve := \pr_{\bbR^{d}\times\{0\}}(\gamma_\ve)
,
$$ 
where $b_\ve \in \{0\} \times \bbR^{d-1}$. We refer to \figref{ab} for a diagram of the situation. Several observations are in order.

\begin{figure}[h!]
\centering
\begin{subfigure}[b]{0.75\textwidth}
	\includegraphics[width = 1\linewidth,page=1]{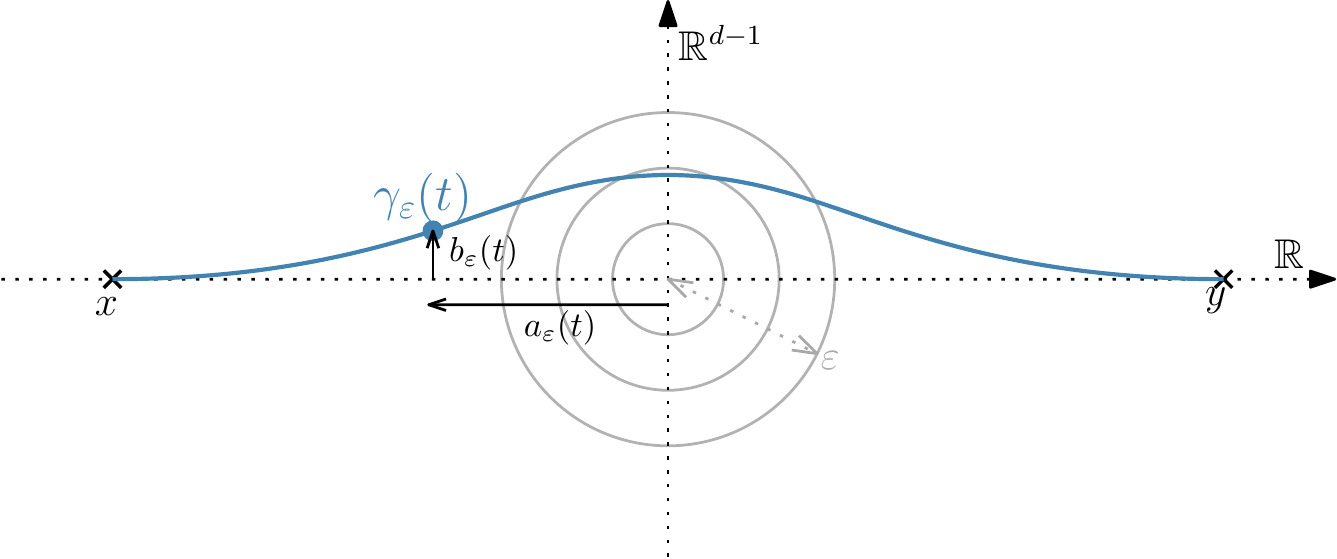}
	\caption{}
	\label{fig:topMve0}
\end{subfigure}

\begin{subfigure}[b]{0.75\textwidth}
	\includegraphics[width = 1\linewidth,page=2]{topMve}
	\caption{}
	\label{fig:topMve1}
\end{subfigure}

\caption{(\subref{fig:topMve0}) Top view of $M_\ve$ and of one of the shortest path between $x$ and $y$, in blue. In light grey is represented the bump of size $\ve$. 
(\subref{fig:topMve1}) Same view of $M_\ve$ as (\subref{fig:topMve0}), illustrating the fact that any shortest path must go from left to right (otherwise one can construct a shorter path, through $s_1$ in the figure) and cannot go outside the shaded area (otherwise one can construct a shorter path, through $s_2$ in the figure).}
\label{fig:ab}
\end{figure}

\begin{itemize}
\item Since $w_\ve(\pm 1) = \pm\ell e_1$, we have $a_\ve(\pm1) = \pm \ell$ and $b_\ve(\pm1) = 0$. Also, because $\gamma_\ve$ is a minimizing path, $a_\ve$ is nondecreasing, and $\|b_\ve\|_\infty \leq \ve$ (see \figref{ab}).
\item Because $\gamma_\ve$ has constant speed on $[-1,1]$, there holds 
\beq \label{eq:gamma}
\|\gamma'_\ve(t)\| = \frac12 \d_{M_\ve}(x,y) \in [A_1 \ell, A_2 \ell],~~~\text{for all}~~~t \in [-1,1],
\eeq
with $A_1,A_2$ depending on $R$ only, uniformly on small $\ve$.
\item $a_\ve$ and $b_\ve$ are smooth and $\gamma_\ve = \bar\Psi_\ve(w_\ve)$.
\item Since $M_\ve$ is symmetric with respect to $\{0\} \times \bbR^d$, so should be the shortest path between $x$ and $y$. This entails in particular that $b_\ve$ is even and that $a_\ve$ is odd;
\item As $\gamma_\ve$ has constant speed and has a curvature bounded from above (as a shortest path in a bounded-curvature space), the ratio $\|\gamma''_\ve \| / \|\gamma_\ve'\|^2$ is bounded in sup-norm by a constant depending on $R$ only. Therefore, there exists a constant $B > 0$ depending on $R$ only such that, uniformly on $\ve$ small enough,
\beq \label{eq:abder}
\max\{\|a_\ve'\|_\infty / \ell, \|a_\ve''\|_\infty/ \ell^2,\|b_\ve'\|_\infty/\ell,\|b_\ve''\|_\infty/\ell^2\} \leq B.
\eeq
\item By symmetry also, $\gamma_\ve$ crosses the hyperplane $\{0\} \times \bbR^d$ orthogonally. As a consequence $\inner{\gamma'_\ve(0)}{e_{d+1}} = 0$, $b_\ve'(0) = 0$ and
$$
a_\ve'(0) = \|w_\ve'(0)\| = \| \gamma'_\ve(0) \| \in \left[ A_1 \ell, A_2\ell\right], 
$$
where $A_1$ and $A_2$ were introduced in~\eqref{eq:gamma}.  
\item Finally, using~\eqref{eq:abder}, we deduce that there exists $C > 0$ depending on $R$ only such that for all $t \in [-1,1]$,
\beq \label{eq:dera}
\begin{cases}
|a_\ve(t) - a_\ve'(0) t | &\leq C \ell^2 t^2, \\
|a_\ve'(t) - a_\ve'(0) | &\leq C \ell^2 t, \\
|a_\ve'(t)a_\ve(t) - a_\ve'(0)^2 t | &\leq C \ell^3 t^2, \\
|b_\ve(t) - b_\ve(0)| &\leq C \ell t.
\end{cases}
\eeq
\end{itemize}

\subsubsection{Perturbative Expansion of the Geodesic Length} We let $\gamma_0(t) := \bar\Psi_0(a_\ve(t) e_{1})$. Although not constant-speed, monotonicity of $a_\ve$ implies that $\gamma_0$ is the shortest path in $M_0$ between $x$ and $y$, and we get, using~\eqref{eq:abder} and~\eqref{eq:dera}, that for some constant $A_3$ depending on $R$,
\beq \label{eq:gamma0}
  \frac12 A_1\ell  \leq \|\gamma'_0(t)\| \leq A_3 \ell ~~~~\text{if}~~\ell \leq \frac{A_1}{2C}, 
\eeq
which we will assume henceforth. Furthermore, the velocity of $\gamma_\ve$ writes 
\begin{align*} 
\gamma_\ve' &= \d\bar\Psi_\ve(w_\ve)[w_\ve'] = \d\bar\Psi_0(w_\ve)[w_\ve'] +c \ve^{k-1}    \inner{\nabla K(w_\ve/\ve)}{w_\ve'} e_{d+1} \\
&=  w_\ve' + \inner{\nabla \Psi_0(w_\ve)}{w_\ve'}e_{d+1}+c\ve^{k-1}    \inner{\nabla K(w_\ve/\ve)}{w_\ve'} e_{d+1} \\
&= 
a_\ve' e_1 +b_\ve'
+
\bigl( 
\underbrace{\inner{\nabla \Psi_0(a_\ve)}{a_\ve'}}_{:=\nabla_0}
+ 
\underbrace{c\ve^{k-1}\inner{\nabla K(w_\ve/\ve)}{w_\ve'}
}_{:=\nabla_1}
\bigr)
e_{d+1},
\end{align*} 
where we used the fact that $\Psi_0$ depends only on its first variable. 
We write the last term as $(\nabla_0 + \nabla_1)e_{d+1}$. 
Using that each three terms in the preceding development are orthogonal, we obtain
\beq \label{eq:devgv}
\|\gamma_\ve'\|^2 = a_\ve'^2 + \|b_\ve'\|^2+(\nabla_0 + \nabla_1)^2 = \underbrace{a_\ve'^2 +\nabla_0^2}_{= \|\gamma_0'\|^2}+ \underbrace{\|b_\ve'\|^2+2 \nabla_0 \nabla_1 + \nabla_1^2}_{:=Q_\ve},
\eeq
and it only remains to study the last three terms, denoted by $Q_\ve$. First, notice that using~\eqref{eq:abder}, one can find two constants $D_0$ depending on $R$ such that  $Q_\ve \geq -D_0 \ve^2 \ell^2$. Together with~\eqref{eq:gamma0}, this yields that $Q_\ve/\|\gamma_0'\|^2 \geq -1$ for $\ve$ small enough (depending on $R$). Likewise, we can show that $Q_\ve \leq D_1 (\ell^2+\ell^2\ve^2+\ve^4)$, for some constant $D_1$ depending on $R$. This again yields  
\beq \label{eq:qe}
\frac{Q_\ve}{\|\gamma_0'\|^2} \leq D_2 ~~~~\text{if}~~~~\ve \leq D_3 \ell,
\eeq
for some constants $D_2$ and $D_3$ depending on $R$ only. All in all, we have that $Q_\ve/\|\gamma_0'\| \in [-1,D_2]$. Using that
$$
\sqrt{1+z} \geq \begin{cases}
1 + z~~~~~&\text{if}~~~z \in [-1,0], \\
1 + D_4 z~~~~~&\text{if}~~~ z\in [0,D_2],~\text{with}~~D_4 = \frac1{D_2}(\sqrt{1+D_2}-1),
\end{cases}
$$
we can finally derive from~\eqref{eq:devgv} and~\eqref{eq:qe} the following bound
\begin{align} \label{eq:newdev} 
\|\gamma_\ve'\| = \|\gamma_0'\| \sqrt{1+\frac{Q_\ve}{\|\gamma_0'\|^2}} \geq   \|\gamma_0'\|  + \tau(Q_\ve) Q_\ve,
\end{align}  
where
$$\tau(z) := \frac{2}{A_1\ell} \ind_{z < 0} + \frac{D_4}{A_3\ell}  \ind_{z \geq 0},
$$
and where we also used~\eqref{eq:gamma0} to bound $1/\|\gamma_0'\|$.
In particular, integrating \eqref{eq:newdev} over $[-1,1]$ yields that
\begin{align*}
\d_{M_\ve}(x,y) \geq \d_{M_0}(x,y) +
\int_{-1}^1 \tau(Q_\ve) Q_\ve
.
\end{align*}

To obtain a more explicit bound, let us now study $Q_\ve$.
For this, first rewrite $\nabla_0$ and $\nabla_1$ more explicitly as
$$
\nabla_0 = -\frac{a_\ve a'_\ve}{\sqrt{R^2-a_\ve^2}}~~~~\text{and}~~~~\nabla_1 = -2c \ve^{k-2}\frac{K(w_\ve/\ve)}{(1-\|w_\ve/\ve\|^2)^2} \inner{w_\ve}{w_\ve'}
.
$$
Hence, noticing that $\inner{w_\ve}{w_\ve'} = a_\ve a_\ve'+\inner{b_\ve}{b_\ve'}$, one can write $2\nabla_0 \nabla_1$ as $P_0+P_1$ with
\begin{align*} 
\begin{cases}
P_0 &=   \ve^{k-2} (a_\ve a'_\ve)^2 T_\ve  \\ 
P_1 &=   \ve^{k-2} T_\ve a_\ve a'_\ve \inner{b_\ve}{b'_\ve}
\end{cases}
~~~~ \text{with}~~~T_\ve:=  \frac{4c K(w_\ve/\ve)}{\sqrt{R^2-a_\ve^2} (1-\|w_\ve/\ve\|^2)^2}. 
\end{align*} 
For $\ell \leq A_1/4C$, condition~\eqref{eq:dera} together with $a'_\ve(0) \geq A_1 \ell/2$ imply that 
$$
\|w_\ve(t) \| \geq |a_\ve(t)| \geq \ve~~~\text{for all~}|t| \geq t_\ve~~\text{with}~~t_\ve := \frac{4\ve}{A_1 \ell},
$$
so that in particular, $T_\ve(t) = 0$ for $|t| \geq t_\ve$. Furthermore, notice that, provided that $\ell$ is small before $R$, $T_\ve$ is bounded by some constant $E > 0$ depending on $R$ only. Using again~\eqref{eq:dera}, we find that for $\ell \leq A_1^2/8C$, there holds
\begin{align} 
&(a_\ve' a_\ve)^2(t) \geq \frac{1}{2} a_\ve'(0)^4 t^2 - C^2 \ell^6 t^4 \geq \frac{1}{32} A_1^4 \ell^4 t^2, \label{eq:aaprime}\\
\text{and}~~~~~~~~
&|a_\ve' a_\ve |(t) \leq a_\ve'(0)^2 |t| + C \ell^3 t^2 \leq 5 A_1^2 \ell^2 |t|~~~~~\text{for all}~~t\in [-1,1]. \nonumber
\end{align} 
In particular, we find that
\begin{align*} 
\int_{-1}^1 |P_1(t)| \diff t 
&\leq 
5   \ve^{k-2}  E A_1^2  \ell^2 \|b_\ve\|_\infty \|b'_\ve\|_\infty  \int_{-t_\ve}^{t_\ve} |t| \diff t 
\\
&= 5   \ve^{k-2}  E A_1^2   \ell^2  \|b_\ve\|_\infty \|b'_\ve\|_\infty  t_\ve^2 
\\
&\leq 
80 B E \|b_\ve\|_\infty  \ell^2 \ve^{k}
,
\end{align*} 
where we used~\eqref{eq:abder} in the last inequality. On the other hand, letting $t_0 \in (-1,1)$ be a time at which $\|b_\ve(t_0)\| =\|b_\ve\|_\infty$,
notice that
\begin{align*} 
\int_{-1}^1  \|b_\ve'\|^2 
&= 
\int_{-1}^{t_0}  \|b_\ve'\|^2 + \int_{t_0}^1  \|b_\ve'\|^2 
\\
&\geq 
\frac{1}{1+t_0} \left\|\int_{-1}^{t_0}  b_\ve' \right\|^2 + \frac{1}{1-t_0} \left\| \int_{t_0}^1  b_\ve'\right\|^2 
\\
&= 
\left(\frac{1}{1+t_0}+\frac{1}{1-t_0}\right) \|b_\ve(t_0)\|^2 \geq 2 \|b_\ve\|^2_\infty.
\end{align*} 
Integrating~\eqref{eq:devgv} and using that $\nabla_1^2 \geq 0$ thus yields 
\beq \label{eq:finaldev}
\int_{-1}^1 \tau(Q_\ve) Q_\ve 
\geq 
2 \|b_\ve\|_\infty\( \tau_1 \|b_\ve\|_\infty -40  \tau_2  B E  \ell^2 \ve^{k}\) + \tau_1  \ve^{k-2} \int_{-1}^1 (a_\ve a'_\ve)^2 T_\ve.
\eeq
where $\tau_1$ is the smallest value of $\tau$, and $\tau_2$ its greatest value.
Now we distinguish on the value of $\|b_\ve(0)\|$:
\begin{itemize}
\item
If $\|b_\ve(0)\| \geq \ve/2$, then $\|b_\ve\|_\infty \geq \ve/2$ and for $\ve$ small enough, we get, noticing that the last term in~\eqref{eq:finaldev} is non-negative, 
$$
\int_{-1}^1 \tau(Q_\ve) Q_\ve  \geq c_R\ve(\ve/2 - \ve/4)/\ell \geq c_R \ve^2/\ell.
$$
\item
Otherwise, if $\|b_\ve(0)\| \leq \ve/2$, then, using~\eqref{eq:dera}, we find that
$$
\begin{cases}\| b_\ve(t) \| &\leq 3\ve/4, \\
|a_\ve(t)| &\leq \ve/2, \\
\end{cases}
~~~\text{for all~} |t| \leq t^*_\ve~~~\text{with}~~ t^*_\ve := \min\{\frac{\ve}{4C\ell},\frac{\ve}{8A_2\ell},\frac{2A_2}{C\ell}\}. 
$$
For $\ve$ small before $R$, $t^*_\ve$ is of the form $t^*_\ve = G\ve/\ell$ with $G$ depending on $R$ only. Furthermore, notice that for $|t| \leq t_\ve^*$, there holds $\|w_\ve(t)\|^2 = |a_\ve(t)|^2 + \|b_\ve(t)\|^2 \leq 13 \ve^2/16$. 
In particular, $T_\ve$ is lower-bounded on $[-t_\ve^*,t_\ve^*]$ by a constant $H$ depending on $R$ only. 
Noticing that $t^*_\ve \leq t_\ve$, we can use the inequality in~\eqref{eq:aaprime} to obtain
\begin{align*} 
 \int_{-1}^1 (a_\ve a'_\ve)^2 T_\ve 
 \geq  
 \frac1{32} A_1^4 \ell^4 H \int_{-t_\ve^*}^{t_\ve^*} t^2 \diff t
 = 
 \frac{1}{48} A_1^4 H G^3 \ell \ve^3.
\end{align*} 
Finally, since $z \mapsto z(z-\nu)$ is minimal on $\bbR_+$ at $z = \nu/2$ with minimal value $-\nu^2/4$, we find the bound
\begin{align*}
\int_{-1}^1 \tau(Q_\ve) Q_\ve \geq c_R \ve^{k+1} - c'_R  \ell^3 \ve^{2k} \geq c_R \ve^{k+1},
\end{align*}  
provided that $\ve$ is small enough before $R$.
\end{itemize}
In both cases, we find that $\int_{-1}^1 \tau(Q_\ve) Q_\ve  \geq c_R \ve^{k+1}$. Now integrating~\eqref{eq:newdev} gives
$$
\d_{M_\ve}(x,y) \geq \d_{M_0}(x,y) + c'_R \ve^{k+1} > \d_{M_0}(x,y).
$$
Finally,~\eqref{eq:gamma} yields $\d_{M_\ve}(x,y)  \geq 2A_2 \ell$ and letting $\ell := (1 \vee D_3^{-1}) \ve$, which we can from~\eqref{eq:qe}, finally gives
\beq 
\label{eq:vek}
\left|1 - \frac{\d_{M_0}(x,y)}{\d_{M_\ve}(x,y)} \right| \geq c_R \ve^{k}.
\eeq
\subsubsection{Concluding with Le Cam's lemma} We apply \lemref{lecammet} with $M_0$ and $M_1 := M_\ve$ for $\ve$ properly chosen. 
Their volumes are bounded from above and below by something depending on $R$ and $d$ only, so that the uniform distribution on $M_0$ and $M_\ve$ are in $\distributions{k}{\rch_{\min}}{\mathbf{L}}{f_{\min}}{f_{\max}}$ provided that $f_{\min}$ and $f_{\max}$ are respectively small enough and large enough compared to $1/R^d$. Finally, we set $R = 2\rch_{\min}$ and $\ve = (C_{\rch_{\min},d} n)^{-1/d}$. For $n$ large enough so that all previous controls are verified, \lemref{lecammet} finally yields
$$
 \inf_{\wh\d} \sup_{P \in \cP^k} \bbE_{P^{\otimes n}}[\ell_\infty(\wh \d|\d_M)] \geq \frac12 c_{\rch_{\min}} \ve^{k}(1- C_{\rch_{\min},d} n \ve^d) \geq c_{\rch_{\min},d,k} n^{-k/d},
$$
where the total variation was bounded using~\cite[Lemma~7]{Berenfeld22}.
\end{proof}

\subsection{Plug-in Estimation for Metric Learning}
\label{sec:proof-plugin-metric}
We start by giving the proof of \prpref{metric}.
\begin{proof}[Proof of \prpref{metric}]
Let $x,y \in K$. Notice that, since $K \subset (K')^\ve$, there holds trivially that $\d_{(K')^\ve}(x,y) \leq d_{K}(x,y)$. For the converse inequality, let $\gamma : [0,1] \to \bbR^D$ be a continuous path in $(K')^\ve$ between $x$ and $y$. Since $\ve < \rch(K)/2$ the closest-point projection on $K$ is well-defined on $(K')^\ve \subset K^{2\ve}$ and we can consider $\gamma_0 = \pr_K \circ \gamma$, which is a continuous path in $K$. For any subdivision $0 = t_0 < t_1 < \dots < t_k = 1$, there holds
$$
\sum_{i=0}^{k-1} \|\gamma_0(t_{i+1}) - \gamma_0(t_i) \| \leq  \frac{\rch(K)}{\rch(K)-2\ve} \sum_{i=0}^{k-1} \|\gamma(t_{i+1}) - \gamma(t_i) \| 
$$
where we used the fact that $\pr_{K}$ is $ \rch(K)/(\rch(K)-2\ve)$-Lipschitz on $K^{2\ve}$~\citep[Theorem~4.8~(8)]{Federer59}. Taking the supremum over all subdivision yields 
$$
\d_{K}(x,y) \leq L(\gamma_0) \leq   \frac{\rch(K)}{\rch(K)-2\ve} L(\gamma)
$$
and then taking the infimum on all continuous path $\gamma$ finally gives
$$
\d_{(K')^{\ve}}(x,y) \geq \left(1 - \frac{2\ve}{\rch(K)}\right) \d_{K}(x,y)
$$
ending the proof.
\end{proof}
To prove \thmref{metricub}, an intermediate result  that bounds the intrinsic diameters of the supports in our statistical model is needed.

\begin{lemma} \label{lem:boundgeo} For any $P \in \distributions{k}{\rch_{\min}}{\mathbf{L}}{f_{\min}}{f_{\max}}$, if $M = \support(P)$, then 
$$
\sup_{x,y \in M} \d_M(x,y) \leq \d_{\max}.
$$
where $\d_{\max}$ is defined in \thmref{metricub}.
\end{lemma}
\begin{proof}[Proof of Lemma~\ref{lem:boundgeo}] We let $x_1,\dots,x_N$ be a $\rch_{\min}/4$-packing of $M$. We let $x,y \in M$, and $G$ be the neighborhood graph built on top of $x,y,x_1,\dots,x_N$ with connectivity radius $\rch_{\min} / 2$. Using~\cite[Theorem~6.3]{niyogi2008finding}, denoting $z_0 = x, z_1,\dots, z_k = y$ the shortest path between $x$ and $y$ in $G$, there holds
\[
\d_M(x,y) \leq \sum_{i=0}^{k-1} \d_M(z_i,z_{i+1}) \leq  \sum_{i=0}^{k-1} 2 \|z_i - z_{i+1}\| \leq  k \rch_{\min}.
\]
But now $k \leq N-1$ and 
\[
N \leq \frac{\vol_d(M)}{\min_{x \in M} \vol_d(M \cap \ball(x,\rch_{\min}/4))} \leq  \frac{\vol_d(M)}{\(1-1/8^2\)^{d/2} \omega_d \rch_{\min}^d/4^d}, 
\]
where we used~\cite[Lemma~5.3]{niyogi2008finding}. Noticing that $\vol_d(M) \leq 1/f_{\min}$, we easily conclude.
\end{proof}

We are now in position to prove \thmref{metricub}.
\begin{proof}[Proof of \thmref{metricub}] We let $\cA_n := \{\dh(\wh M, M) \leq \ve_n\}$ denote the event where $\wh M$ is $\ve_n$-precise in Hausdorff distance, and we take $x, y \in M$. On the event $\cA_n$, for $n$ large enough such that $\ve_n \leq \rch_{\min}/2$, \prpref{metric} applies to $K' = \wh M \cup \{x,y\}$ and, together with \lemref{boundgeo}, yields 
$$\left| 1 - \frac{\wh \d(x,y)}{\d_M(x,y)}\right| \leq \frac{2\ve_n}{\rch_{\min}}.$$ 
On $\cA_n^c$, we distinguish whether $\|x- y\| \leq \ve_n$ or not. If so, then $\wh\d(x,y) = \|x-y\| \leq \d_M(x,y)$. In the other case, $\d_M(x,y) \geq \|x-y\| \geq \ve_n$ and $\d(x,y) \leq \d_{\max}$ so that, in any case
$$
\left| 1 - \frac{\wh \d(x,y)}{\d_M(x,y)}\right| \leq 1 + \frac{\wh \d(x,y)}{\d_M(x,y)} \leq 1 + \frac{\d_{\max}}{\ve_n},
$$ 
for $n$ large enough such that $\ve_n \leq \d_{\max}$. Patching these two bounds together yields
$$
 \bbE_{P^{\otimes n}}[\ell_\infty(\wh \d| \d_M)] \leq \frac{2\ve_n}{\rch_{\min}} P^{\otimes n}(\cA_n)+ \left(1 + \frac{\d_{\max}}{\ve_n}\right) P^{\otimes n}(\cA_n^c),
$$
ending the proof.
\end{proof}

\section{Proofs of \secref{optimal_reach_estimation}} \label{sec:proofoptireach}

We first prove that submanifolds of the model do fulfill \assref{spread} and \assref{subeuc}.

\begin{proof}[Proof of \prpref{assreach}] 
\assref{subeuc} is a simple consequence of~\cite[Proposition~6.3]{niyogi2008finding} which yields fulfillment for $\Delta_1 = \rch(M)/2$ and $C_1 = 2$. 
For \assref{spread}, take $x,y \in M$ such that $\|x-y\| \leq \rch(M)$ and take $\ve < \rch(M)/4$.  
We consider $a = \exp_y(v)$, where 
$$
v = - \ve \frac{\pr_{T_y M}(x-y)}{\|\pr_{T_y M}(x-y)\|}.
$$
Thanks to~\cite[Theorem~4.8 (7)]{Federer59}, there holds
\begin{align*} 
\|\pr_{T_y M}(x-y)\|^2 
&= 
\|x-y\|^2 - \d^2(x-y,T_y M) 
\\
&\geq \|x-y\|^2- \frac{\|x-y\|^4}{4 \rch^2(M)} \\
&\geq \frac34 \|x-y\|^2,
\end{align*} 
and
$$
\inner{v}{y-x} =  \ve \frac{\inner{x-y}{\pr_{T_y M}(x-y)}}{\|\pr_{T_y M}(x-y)\|} = \ve \|\pr_{T_y M}(x-y)\| \geq \frac12 \ve \|x-y\|,
$$
so that
\[
\|x-y-v\|^2 \geq \|x-y\|^2+ \ve \|x-y\| +\ve^2 \geq \(\|x-y\|+\frac12\ve\)^2,
\]
and thus $\|x-y-v\| \geq \|x-y\|+\ve/2$. But now $\|x-a\| \geq \|x-y-v\|-\|a-y-v\|$ and $\|a-y-v\| \leq 5\ve^2/4\rch(M)$ according to~\cite[Lemma~1]{Aamari19}. All in all, we get that
$$
\|x-a\| \geq \|x-y\|+\frac12 \ve-\frac{5}{4\rch(K)}\ve^2 \geq \|x-y\|+\frac{3}{16} \ve,
$$
ending the proof.
\end{proof}

To prove \thmref{sdrub}, a bound on the metric distortion between our distance estimator and $\d_M$ is needed, that easily follows from  Proposition~\ref{prp:metric}.

\begin{proposition} \label{prp:metdisto} In the context of \prpref{metric}, we have that for all $\delta > 4\ve$,
$$
\D_\delta(\d_{K}, \d_{(K')^\ve}) \leq 1 + \frac{4\ve}{(\delta-4\ve) \wedge \rch(K)}.
$$
\end{proposition}
\begin{proof}[Proof of \prpref{metdisto}] \prpref{metric} already gives that $\D_\delta(\d_{K}|\d_{(K')^{\ve}}) \leq 1 + 2\ve/\rch(K)$. For the other control, notice that for any two $x,y \in (K')^\ve$ that are $\delta$-apart for the Euclidean distance, there holds denoting $x_0 = \pr_{K}(x)$ and $y_0 = \pr_{K}(y)$, 
$$\d_{(K')^\ve}(x,y) \leq 4\ve + \d_{K}(x_0,y_0)
$$
because the piecewise-defined path consisting of the segment $[x,x_0]$ of the (or a near-minimizing) shortest-path between $x_0$ and $y_0$ in $K$, and of the segment $[y_0,y]$, is a continuous path in $(K')^\ve$ between $x$ and $y$ of length the RHS of the display above. Now notice that
$$
 \d_{K}(x_0,y_0) \geq \|x_0 - y_0\| \geq \delta - 4\ve,
$$  
which immediately yields $ \D_\delta(\d_{(K')^\ve}|\d_{K}) \leq 1 + \frac{4\ve}{\delta-4\ve}$.
\end{proof}
The rate of the plug-in SDR estimator follows straightforwardly. 

\begin{proof}[Proof of \thmref{sdrub}] Let $\cA_n := \{\dh(M,\wh M) \leq \ve_n\}$. On this event, we have $\D_{\delta}(\wh\d,\d_M) \leq 1 + 8\ve_n/\delta$ according to \prpref{metdisto}, so that applying \thmref{stab} with $\delta_0 = \delta/2$, $\ve= \ve_n$ and $\nu = 8\ve_n/\delta$ yields $|\wh\sdr_\delta -\sdr_\delta(M,\d_M)| \leq  \zeta_0 \ve_n$ with $
\zeta_0 \leq C \s_{\max}^4/\delta^4$. We conclude that
$$
 \bbE_{P^{\otimes n}} |\wh\sdr_\delta -\sdr_\delta(M,\d_M)| \leq \zeta_0 \ve_n P^{\otimes n }(\cA_n) + 2 \s_{\max} P^{\otimes n }(\cA_n^c),
$$
which ends the proof.
\end{proof}

\end{document}